\documentclass[a4paper,oneside,11pt,final]{article}

\usepackage[
  a4paper,
  left=20mm,
  right=20mm,
  top=20mm,
  bottom=20mm
]{geometry}

\usepackage[T1]{fontenc}
\usepackage[utf8]{inputenc}

\usepackage{lmodern}

\usepackage{amsmath}
\usepackage{amssymb}
\usepackage{amsthm}
\usepackage{mathtools}

\usepackage{titlesec}

\titleformat{\section}
  {\normalfont\sffamily\Large\bfseries}
  {\thesection}{1em}{}

\titleformat{\subsection}
  {\normalfont\sffamily\large\bfseries}
  {\thesubsection}{1em}{}

\titleformat{\subsubsection}
  {\normalfont\sffamily\normalsize\bfseries}
  {\thesubsubsection}{1em}{}

\titleformat{\paragraph}
  {\normalfont\sffamily\normalsize\bfseries}
  {\theparagraph}{}{}

\usepackage{setspace}
\usepackage[backend=biber,style=numeric,giveninits=true,noerroretextools,maxbibnames=99]{biblatex}
\usepackage{csquotes}

\DeclareSourcemap{
  \maps[datatype=bibtex]{
     \map{
        \step[fieldsource=doi,final]
        \step[fieldset=isbn,null]
        \step[fieldset=issn,null]
        }
      }
}

\AtEveryBibitem{\clearfield{url}} 
\AtEveryCitekey{\clearfield{url}} 

\usepackage{hyperref} 
\usepackage[english]{babel}
\usepackage{enumerate}
\usepackage{aliascnt}
\usepackage{verbatim}
\usepackage[normalem]{ulem}
\usepackage{esint} 

\usepackage[noabbrev,capitalize]{cleveref}

\usepackage[right]{showlabels}

\usepackage[nottoc]{tocbibind} 
\usepackage[title]{appendix}
\usepackage{xcolor}
\usepackage{bm} 
\usepackage{tikz}

\numberwithin{equation}{section} 

\usepackage{hyperref} 
\hypersetup{final,
  colorlinks,
  linkcolor={blue!70!black},
  citecolor={blue!40!black},
  urlcolor={blue!40!black}} 

\let\etoolboxforlistloop\forlistloop 
\usepackage{autonum}
\let\forlistloop\etoolboxforlistloop 

\allowdisplaybreaks[4] 

\newtheoremstyle{sansplain}
  {3pt}   
  {3pt}   
  {\itshape} 
  {}      
  {\sffamily\bfseries} 
  {.}     
  {.5em}  
  {}      

\newtheoremstyle{sansdefinition}
  {3pt}
  {3pt}
  {\itshape}
  {}
  {\sffamily\bfseries}
  {.}
  {.5em}
  {}

\newtheoremstyle{sansremark}
  {3pt}
  {3pt}
  {\normalfont}
  {}
  {\sffamily}
  {.}
  {.5em}
  {}

\theoremstyle{sansplain}

\newtheorem{proposition}{Proposition}[section]

\newaliascnt{lemma}{proposition} 

\newtheorem{lemma}[lemma]{Lemma}

\aliascntresetthe{lemma} 
\Crefname{lemma}{Lemma}{Lemmas}

\newaliascnt{theorem}{proposition} 
\newtheorem{theorem}[theorem]{Theorem}

\aliascntresetthe{theorem} 

\newaliascnt{corollary}{proposition} 
\newtheorem{corollary}[corollary]{Corollary}

\aliascntresetthe{corollary}

\newaliascnt{hypothesis}{proposition} 
\newtheorem{hypothesis}[hypothesis]{Hypothesis}

\aliascntresetthe{hypothesis} 

\theoremstyle{sansdefinition}

\newaliascnt{definition}{proposition} 
\newtheorem{definition}[definition]{Definition}

\aliascntresetthe{definition} 
\Crefname{definition}{Definition}{Definitions}

\newaliascnt{problem}{proposition} 

\aliascntresetthe{problem} 

\newaliascnt{example}{proposition} 

\aliascntresetthe{example} 

\theoremstyle{sansremark}

\newaliascnt{remark}{proposition} 
\newtheorem{remark}[remark]{Remark}

\aliascntresetthe{remark} 

\def\equationautorefname~#1\null{%
	(#1)\null
}

\newcommand{\R}{\mathbb{R}}
\newcommand{\C}{\mathbb{C}}

\newcommand{\N}{\mathbb{N}}

\renewcommand{\S}{\mathbb{S}}

\newcommand{\W}{\mathcal{W}}
\newcommand{\tr}{\operatorname{tr}}

\newcommand*{\dd}{\mathop{}\!\mathrm{d}}

\newcommand{\vol}{\mathrm{Vol}}

\title{
    A new localization approach for Willmore flows:\\ Asymptotic analysis for the Dirichlet problem
}

\author{
    \large{Simon Blatt}\thanks{University of Salzburg, Hellbrunner Str.~34, 5020 Salzburg, Austria. \texttt{simon.blatt@plus.ac.at}} 
    \and 
    \large{Manuel Schlierf}\thanks{University of Salzburg, Hellbrunner Str.~34, 5020 Salzburg, Austria. \texttt{math@manuelschlierf.info}}
}

\date{\today}
\date{}

\begin{document}
\maketitle
\vspace{-0.5cm}
\begin{abstract}
    We study the Willmore flow with Dirichlet boundary conditions by means of a new geometric localization approach. In contrast to previous works, we do not rely on localized energy estimates, but on uniformly controlled isothermal coordinates and localized parabolic a priori estimates. For genus-zero surfaces satisfying an explicit energy condition depending only on the prescribed boundary data, we prove global existence together with a dichotomy between smooth convergence and subconvergence to a complete noncompact Willmore surface. If the initial Willmore energy is at most $4\pi$, this dichotomy reduces to smooth convergence after reparametrization, and we show that this additional threshold is optimal. As an application, we obtain an existence criterion for the stationary Dirichlet Willmore problem.
\end{abstract}

{\small

}

\section{Introduction and main results}

Let $\Sigma$ be a compact oriented surface, possibly with boundary, and let $f\colon\Sigma\to\R^3$ be a smooth immersion. The Willmore energy of $f$ is defined by
\begin{equation}\label{eq:def-Willmore-energy}
  \W(f)\vcentcolon=\frac14\int_\Sigma\mathbf H_f^2\dd\vol_{g_f},
\end{equation}
where $g_f=f^*g_{\R^3}$ is the induced metric and $\dd\vol_{g_f}$ its area measure. In positively oriented local coordinates, the \emph{unit normal}, the \emph{second fundamental form}, the \emph{mean curvature}, and the \emph{traceless second fundamental form} are given by
\begin{equation}\label{eq:geometric-conventions}
  \mathbf n_f=\frac{\partial_1f\times\partial_2f}{|\partial_1f\times\partial_2f|},\quad (\mathbf A_f)_{ij}=\langle\partial_{ij}^2f,\mathbf n_f\rangle,\quad \mathbf H_f=\mathrm{tr}_{g_f}\mathbf A_f,\quad\text{and}\quad \mathbf A_f^0=\mathbf A_f-\frac12\mathbf H_fg_f,\quad\text{respectively}.
\end{equation}
In particular, $\mathbf H_f$ is the sum of the principal curvatures. We denote by $\nabla_{g_f}$ the covariant derivative and use the sign convention $\Delta_{g_f}=\mathrm{div}_{g_f}\nabla_{g_f}$, so that $\Delta_{g_f}f=\mathbf H_f\mathbf n_f$. 

The Willmore energy goes back to the study of classical conformal geometry, see \cite{blaschkethomsen1929}. First results on its minimizers are due to Willmore \cite{willmore1965,willmore1993}: for closed surfaces, one has $\W(f)\geq4\pi$, with equality precisely for round spheres. More generally, the Li--Yau inequality \cite{liyau1982} yields that $\W(f)\geq4\pi m$ if some point of $\R^3$ has at least $m$ preimages under $f$. In particular, an immersion of a closed surface with energy below $8\pi$ is embedded.

For a variation $V\colon\Sigma\to\R^3$ supported away from $\partial\Sigma$, the first variation is given by
\begin{equation}\label{eq:first-variation-Willmore-interior}
  \partial_{t,0}\W(f+tV)=\int_\Sigma\langle\nabla\W(f),V\rangle\dd\vol_{g_f}\quad\text{where}\quad \nabla\W(f)=\frac12\big(\Delta_{g_f}\mathbf H_f+|\mathbf A_f^0|^2\mathbf H_f\big)\mathbf n_f.
\end{equation}
An immersion satisfying the PDE $\nabla\W(f)=0$ is called a \emph{Willmore immersion}. Existence of smooth closed minimizers of prescribed genus was established in \cite{simon1993,bauerkuwert2003}. In particular, the monotonicity formula Simon obtained in \cite{simon1993} is a fundamental tool in the compactness theory, and is also crucial for our main result on the Willmore flow with Dirichlet boundary conditions below. A regularity theory and removability results for point-singularities of Willmore immersions were developed in \cite{kuwertschaetzle2004,riviere2008}. For spherical topology, the classification in \cite{bryant1984} implies that the energy of a Willmore immersion is an integer multiple of $4\pi$, and that every non-round Willmore sphere has energy at least $16\pi$.

In this work, we study the Willmore flow with Dirichlet boundary conditions. Our main result implies global existence and a uniform curvature bound for immersed genus-zero surfaces below an explicit energy threshold depending only on the prescribed boundary data. The proof uses a new localization scheme,  which allows proving a lifespan result and parabolic smoothing estimates up to the boundary. A second ingredient is a rigidity result for stationary boundary blow-ups. 

From now on, let $\Sigma$ be connected with $\partial\Sigma\neq\emptyset$, let $\Gamma\subset\R^3$ be a smooth not necessarily connected compact one-dimensional submanifold, and let $\mathbf N\colon\Gamma\to\S^2$ be smooth with $\mathbf N(p)\perp T_p\Gamma$. Fix a smooth embedding $\gamma\colon\partial\Sigma\to\Gamma$ onto $\Gamma$. We impose the \emph{Dirichlet}, or \emph{clamped}, boundary conditions
\begin{equation}\label{eq:Dirichlet-boundary-conditions}
  f=\gamma\quad\text{and}\quad\partial_{\eta_f}f=\mathbf N\circ\gamma\quad\text{on $\partial\Sigma$},
\end{equation}
where $\eta_f$ is the outer unit conormal to $\partial\Sigma$ in $(\Sigma,g_f)$. Thus, both the boundary position and the ambient outer conormal are prescribed. For variations $V$ preserving \eqref{eq:Dirichlet-boundary-conditions}, the boundary terms in the first variation vanish, and \eqref{eq:first-variation-Willmore-interior} remains valid.

The stationary boundary problem has been studied using several approaches. Uniqueness results for Willmore disks with special boundary data were obtained in \cite{palmer2000,dallacqua2012}. A general variational construction of branched solutions, with smooth embedded solutions under additional energy assumptions, was given in \cite{schaetzle2010}. A different approach was used in \cite{novagapozzetta2020} to establish existence of minimizers in the class of curvature varifolds with boundary and connected supports. Interior regularity for these varifold-minimizers under suitable energy assumptions was established in \cite{ruppscharrer2025}. For surfaces of revolution generated by graphical profile curves, existence for symmetric Dirichlet data was proved in \cite{dallacquafroehlichgrunauschieweck2011}. Moreover, \cite{eichmann2016} shows that uniqueness fails in general, already among surfaces of revolution. The minimization problem for cylinders of revolution was resolved in \cite{eichmannschaetzle2024}, including a geometric obstruction to the existence of a minimizer. More recently, existence for small graphical boundary data with low regularity was established in \cite{gulyak2026}.

\subsection{The Dirichlet Willmore flow: main results}

We consider a maximal smooth solution $f\colon[0,T)\times\Sigma\to\R^3$ of
\begin{equation}\label{eq:Dirichlet-Willmore-flow}
  \begin{cases}
    \partial_t f=-\frac12\big(\Delta_{g_f}\mathbf H_f+|\mathbf A_f^0|^2\mathbf H_f\big)\mathbf n_f&\text{on $[0,T)\times\Sigma$}\\
    f(0,\cdot)=f_0&\text{on $\Sigma$}\\
    f=\gamma\text{ and }\partial_{\eta_f}f=\mathbf N\circ\gamma&\text{on $[0,T)\times\partial\Sigma$}.
  \end{cases}
\end{equation}
Here, $f_0$ is a smooth immersion satisfying the boundary conditions and the compatibility conditions required for a solution smooth up to $t=0$. Maximality means that the solution cannot be smoothly extended beyond $T$ while still satisfying \eqref{eq:Dirichlet-Willmore-flow}. We prescribe the full velocity, which is normal to the evolving surface. The Willmore flow is the $L^2(\dd\vol_{g_f})$-gradient flow of $\W$ in the sense that it decreases Willmore energy according to
\begin{equation}\label{eq:energy-dissipation-Dirichlet}
  \W(f(t))+\int_s^t\int_\Sigma|\nabla\W(f(\tau))|^2\dd\vol_{g_{f(\tau)}}\dd\tau=\W(f(s))\quad\text{for $0\leq s\leq t<T$}.
\end{equation}

For closed surfaces, short-time existence and stability of round spheres were established in \cite{simonett2001}. Initial data of $C^{1+\alpha}$ regularity, $0<\alpha<1$, were later treated in \cite{lecroneshaosimonett2020}. A more general asymptotic analysis of the flow was developed by Kuwert--Sch\"atzle in \cite{kuwertschaetzle2001,kuwertschaetzle2002,kuwertschaetzle2004}. Their energy estimates control the flow's lifespan in terms of local concentration of $\mathbf A_f$ in $L^2$, and yield suitable parabolic smoothing estimates. For spherical initial data with $\W(f_0)\leq8\pi$, \cite{kuwertschaetzle2004} showed that the flow exists for all times and converges to a round sphere. Moreover, they constructed a blow-up limit $\hat f$ along a sequence $t_j\nearrow T$ by rescaling and translating the flow in $\mathbb R^3$. The limit is a nonflat Willmore immersion. By \cite{chillfasangovaschaetzle2009}, the blow-up $\hat f$ is noncompact if and only if a singularity occurs, i.e.\ the flow does not both exist globally and smoothly converge to a Willmore immersion up to reparametrization. This blow-up analysis was used in \cite{dallacquamullerschatzlespener2020} to prove that the Willmore flow exists globally and converges to the Clifford torus up to reparametrization when starting at a torus of revolution with energy at most $8\pi$.

The analysis in the closed case relies on localized $L^2$-type energy estimates and repeated integration by parts. For surfaces with boundary, these computations produce additional boundary terms. Under our Dirichlet boundary conditions, these terms are not directly controlled by the interior estimates. Instead, we obtain the required parabolic smoothing estimates by applying a new localization technique which does not rely on such energy estimates. The recent works \cite{palmurellariviere2022,palmurellariviere2024} introduced a different approach for studying the Willmore flow of closed surfaces, which is based on choosing a tangential speed in the PDE enforcing \emph{conformal parametrization} at each time. 

For the initial energy threshold in our main result, we introduce the boundary functional
\begin{equation}\label{eq:def-Boundary-integral}
  \mathbf B_{\Gamma,\mathbf N}(z)\vcentcolon=\int_\Gamma\frac{\langle\mathbf N(p),p-z\rangle}{|p-z|^2}\dd\mathcal H^1(p)\quad\text{for $z\in\R^3$},
\end{equation}
where $\mathcal H^1$ denotes one-dimensional Hausdorff measure. The functional is well-defined for all $z\in\R^3$, but in general has a directional jump across $\Gamma$, for more details see \Cref{lem:boundary-integral-jump}.

We work with a concept of local smooth convergence up to the boundary, after reparametrization, for limits whose domains may be noncompact. More precisely, the reparametrizations are defined on an exhaustion of the limit surface, preserve its boundary, and their images eventually contain every intrinsic ball of fixed radius about a chosen basepoint. For more details, see \Cref{def:local-smooth-pointed-convergence}. Completeness always refers to the induced intrinsic distance with the boundary included. For noncompact surfaces, their genus is the supremum of $\mathrm{genus}(K)$ where $K$ is a compact subsurface.

\begin{theorem}\label{thm:main-result}
  Let $f$ be a maximal solution of \eqref{eq:Dirichlet-Willmore-flow}. If $\mathrm{genus}(\Sigma)=0$ and
  \begin{equation}\label{eq:1-thm:main-result}
    \W(f_0)\leq\min\Big\{8\pi-2\sup_{z\in\R^3}\mathbf B_{\Gamma,\mathbf N}(z),4\pi-2\sup_{z\in\Gamma}\mathbf B_{\Gamma,\mathbf N}(z)\Big\},
  \end{equation}
  then $T=\infty$, $\sup_{t\in[0,\infty)}\Vert\mathbf A_{f(t)}\Vert_{L^\infty(\Sigma)}<\infty$, and exactly one of the following cases occurs.
  \begin{enumerate}
    \item[(a)] (Infinite-time blow-up) There exist $t_j\nearrow\infty$, $p_j\in\Sigma$, a complete noncompact genus-zero surface $\hat\Sigma$ with $\partial\hat\Sigma\neq\emptyset$, and a proper smooth Willmore immersion $\hat f\colon\hat\Sigma\to\R^3$ and a point $\hat p\in\hat\Sigma$ such that $(f(t_j),p_j)$ converges to $(\hat f,\hat p)$ in the sense of \Cref{def:local-smooth-pointed-convergence}. Moreover, $\Vert\mathbf A_{\hat f}\Vert_{L^\infty(\hat\Sigma)}<\infty$, the restriction $\hat f\vert_{\partial\hat\Sigma}\colon\partial\hat\Sigma\to\Gamma$ is a not necessarily surjective embedding, and $\partial_{\eta_{\hat f}}\hat f=\mathbf N\circ\hat f$ on $\partial\hat\Sigma$.
    \item[(b)] (Smooth convergence) As $t\to\infty$, the flow $f(t)$ converges smoothly after reparametrization to a Willmore immersion $f_\infty\colon\Sigma\to\R^3$, such that $f_\infty\vert_{\partial\Sigma}$ is a reparametrization of $\gamma$, and $\partial_{\eta_{f_\infty}}f_\infty=\mathbf N\circ f_\infty$ on $\partial\Sigma$.
  \end{enumerate}
  Moreover, if $\W(f_0)\leq4\pi$, then case (a) cannot occur.
\end{theorem}

It is important to note that $L^\infty$-control of $\mathbf A_f$ does not by itself imply convergence to a compact surface. In case (a), the subsequence is not rescaled, and the limit may retain only some of the prescribed boundary components. In case (b), smooth convergence means convergence in $C^\infty(\Sigma,\R^3)$, including the boundary, after composition with diffeomorphisms of $\Sigma$.

The smooth-convergence dichotomy in \Cref{thm:main-result} yields the following existence criterion for the stationary Dirichlet problem.

\begin{corollary}\label{cor:existence-clamped-Willmore}
  Let $\mathrm{genus}(\Sigma)=0$ and let $f_0$ be a smooth initial immersion as above. If
  \begin{equation}\label{eq:existence-clamped-Willmore}
    \W(f_0)\leq\min\Big\{4\pi,8\pi-2\sup_{z\in\R^3}\mathbf B_{\Gamma,\mathbf N}(z),4\pi-2\sup_{z\in\Gamma}\mathbf B_{\Gamma,\mathbf N}(z)\Big\},
  \end{equation}
  then there exists a smooth Willmore immersion $f_\infty\colon\Sigma\to\R^3$ with $\W(f_\infty)\leq\W(f_0)$. Moreover, $f_\infty\vert_{\partial\Sigma}$ is a reparametrization of $\gamma$ and $\partial_{\eta_{f_\infty}}f_\infty=\mathbf N\circ f_\infty$ on $\partial\Sigma$.
\end{corollary}
\begin{proof}
  By \Cref{thm:main-result}, the flow converges smoothly after reparametrization. The energy inequality follows from \eqref{eq:energy-dissipation-Dirichlet}.
\end{proof}

\begin{remark}[Lower regularity for competitor $f_0$]
	We sketch an argument to show that, in \Cref{cor:existence-clamped-Willmore}, it is sufficient to only require that, for some $\alpha\in(0,1)$, $f_0\in C^{1+\alpha}(\Sigma)\cap W^{2,2}(\Sigma)$ is an immersion satisfying \eqref{eq:Dirichlet-boundary-conditions} and \eqref{eq:existence-clamped-Willmore}. As a by-product of our lifespan result \Cref{thm:lifespan-and-smoothing}, for some $\delta>0$ there exists a family of immersions $f\in C^{0}([0,\delta],C^1(\Sigma))\cap C^\infty((0,\delta]\times\Sigma)$ satisfying \eqref{eq:energy-dissipation-Dirichlet} for $0<s\leq t\leq\delta$, $f(0)=f_0$, and \eqref{eq:Dirichlet-boundary-conditions} up to reparametrization. If one additionally proved continuous dependence on the data of solutions to the underlying DeTurck-Willmore flow \eqref{eq:DeTurck-Willmore-flow}, the argument used in \cite[Lemma~5.9]{dallacquamuellerruppschlierf2025} could be adapted to show $\liminf_{t\searrow 0}\W(f(t)) = \W(f_0)$. One could then apply \Cref{cor:existence-clamped-Willmore} with $f_0$ replaced by $f(\delta)$.
\end{remark}

Condition \eqref{eq:1-thm:main-result} plays the role of the energy threshold $\W(f_0)\leq 8\pi$ for the Willmore flow of spheres or tori of revolution, as we explain after \Cref{thm:intro-stationary-blow-ups}. Its two terms correspond to excluding curvature concentration in the interior and at the boundary, respectively. Unlike for the Willmore flow of spheres or tori of revolution, however, the noncompact case (a) in \Cref{thm:main-result} can occur under \eqref{eq:1-thm:main-result}. We next give an example for a class of initial data $f_0$ for which case (a) must occur. This also shows that the additional $4\pi$ threshold in \Cref{thm:main-result} excluding (a) is optimal. 

\begin{proposition}\label{prop:infinite-time-blow-up-disk}
  Let $D=\{z\in\R^2\times\{0\}\mid|z|\leq1\}$ be the unit disk, $\Gamma=\partial D$, and $\mathbf N(x)=-x$ for $x\in\Gamma$. Let $f_0\colon D\to\R^3$ be a smooth admissible initial immersion with $f_0(x)=x$ and $\partial_{\eta_{f_0}}f_0=\mathbf N$ on $\Gamma$. If $\W(f_0)\leq6\pi$, then the solution of \eqref{eq:Dirichlet-Willmore-flow} starting at $f_0$ satisfies case (a) of \Cref{thm:main-result}. In particular, $T=\infty$ and
  \begin{equation}\label{eq:diameter-infinite-time-blow-up}
    \limsup_{t\to\infty}\mathrm{diam}(f(t,D))=\infty.
  \end{equation}
\end{proposition}
\begin{proof}
  One computes that $\sup_{z\in\R^3}\mathbf B_{\Gamma,\mathbf N}(z)=0$ and $\mathbf B_{\Gamma,\mathbf N}(z)=-\pi$ for $z\in\Gamma$. Thus, \eqref{eq:1-thm:main-result} reduces to $\W(f_0)\leq6\pi$. If case (b) occurred, the limit would be a Willmore disk with boundary $\Gamma$ and tangent plane $\R^2\times\{0\}$ along $\Gamma$. By \cite[Theorem~1.1]{dallacqua2012}, its image would be $D$. Its outer conormal would then be $x$ on $\Gamma$, contradicting the prescribed value $-x$. Hence case (a) occurs. The complete noncompact proper limit has unbounded image, which in particular yields \eqref{eq:diameter-infinite-time-blow-up}.
\end{proof}

\begin{figure}[htb!]
  \centering
  \begin{tikzpicture}[x=1cm,y=1cm,line cap=round,line join=round]
    \begin{scope}[shift={(-3.7,0)}]
      \draw[gray!50] (-2.25,0)--(-1,0) (1,0)--(2.25,0);
      \draw[thick] (-1,0)..controls(-1.4,0)and(-2,0.2)..(-2,1.45)..controls(-2,3.65)and(2,3.65)..(2,1.45)..controls(2,0.2)and(1.4,0)..(1,0);
      \draw[thick,dashed] (1,0) arc[start angle=0,end angle=180,x radius=1,y radius=0.28];
      \draw[thick] (-1,0) arc[start angle=180,end angle=360,x radius=1,y radius=0.28];
      \draw[thick,->] (-1,0)--(-0.4,0);
      \draw[thick,->] (1,0)--(0.4,0);
      \node at (0,-0.55) {$\Gamma$};
      \node at (0,0) {$\mathbf N$};
    \end{scope}
    \draw[->] (-1.1,1.8)--(0.55,1.8);
    \node at (-0.275,2.15) {increasing $R$};
    \begin{scope}[shift={(3.6,0)}]
      \draw[gray!50] (-2.65,0)--(-1,0) (1,0)--(2.65,0);
      \draw[thick] (-1,0)..controls(-1.5,0)and(-2.4,0.12)..(-2.4,1.9)..controls(-2.4,4.6)and(2.4,4.6)..(2.4,1.9)..controls(2.4,0.12)and(1.5,0)..(1,0);
      \draw[thick,dashed] (1,0) arc[start angle=0,end angle=180,x radius=1,y radius=0.28];
      \draw[thick] (-1,0) arc[start angle=180,end angle=360,x radius=1,y radius=0.28];
      \draw[thick,->] (-1,0)--(-0.4,0);
      \draw[thick,->] (1,0)--(0.4,0);
      \node at (0,-0.55) {$\Gamma$};
      \node at (0,0) {$\mathbf N$};
    \end{scope}
  \end{tikzpicture}
  \caption{Disk-type initial data $f_{0,R}$ which, away from their boundary, are given by round spheres of radius $R$. Their energies tend to $4\pi$ from above.}
  \label{fig:infinite-time-blow-up-initial-data}
\end{figure}
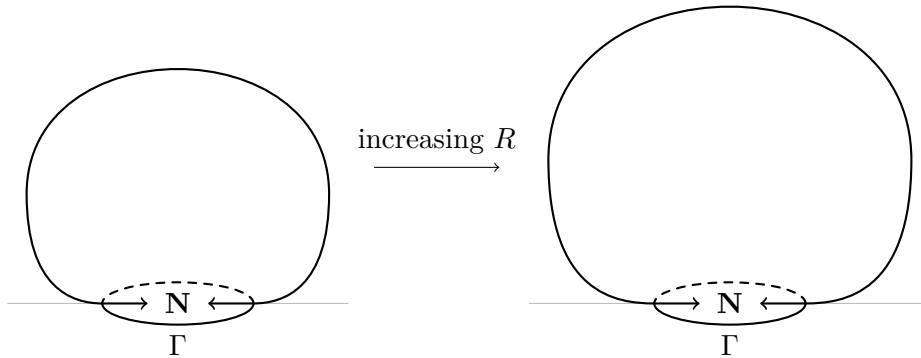

Moreover, initial data $f_0$ as in \Cref{prop:infinite-time-blow-up-disk} exist with energy arbitrarily close to $4\pi$, see \Cref{fig:infinite-time-blow-up-initial-data}. More precisely, we consider a round sphere of radius $R$ centered at $(0,0,\sqrt{R^2-1})$ and remove the small lower cap bounded by $\Gamma$. Near $\Gamma$, modify the remaining surface on a fixed annulus so that it agrees with the plane $\R^2\times\{0\}$ in a collar of $\Gamma$. On this annulus, the spherical graph and its derivatives are of order $R^{-1}$. Applying a fixed smooth cutoff, we can therefore construct embedded disks $f_{0,R}$ with the prescribed boundary and inward conormal, such that
\begin{equation}\label{eq:sharpness-four-pi}
  \W(f_{0,R})=4\pi+O(R^{-2})\quad\text{as $R\to\infty$}.
\end{equation}
The fact that $\W(f_{0,R})>4\pi$ follows immediately from \Cref{thm:main-result}, since case (a) occurs by \Cref{prop:infinite-time-blow-up-disk}.

For comparison, global existence and convergence for cylindrical surfaces of revolution with clamped boundary were proved by the second author in \cite{schlierf2024na,schlierf2025acv}, with an equivalent improved energy condition obtained independently in \cite{eichmann2024}. In that setting, the threshold in \cite{schlierf2025acv} is the first term in \eqref{eq:1-thm:main-result}, with the supremum restricted to the axis of rotation. Starting at cylinders of revolution simplifies the Willmore flow to a one-dimensional problem for profile curves, excludes case (a) in \Cref{thm:main-result} already below energy $8\pi$, and restricts the possible curvature concentration to the axis. A different approach to simplify the geometry a priori is to consider graphs $f_0(x)=(x,u_0(x))$ with $x$ in a sufficiently regular planar domain. For this problem, a time-weighted parabolic theory with homogeneous clamped boundary data was developed in \cite{lecronepruesswilke2014}. The recent work \cite{gulyak2026flow} obtained a stability-type result for the non-homogeneous problem, i.e.\ global existence and exponential convergence for small Lipschitz initial data and small graphical boundary data. 

Our result instead treats arbitrary genus-zero immersions.

\subsection{Proof strategy: New localization approach and classification of boundary blow-ups}

We now describe the two main ingredients of the proof. First, as an application of our new localization approach, we get the following blow-up construction which applies to $\Sigma$ with arbitrary genus, and does not require \eqref{eq:1-thm:main-result}.

\begin{theorem}[Stationary blow-ups]\label{thm:intro-stationary-blow-ups}
  Let $f$ be a maximal solution of \eqref{eq:Dirichlet-Willmore-flow}. If 
  \begin{equation}
  	  T<\infty \quad\text{or}\quad \sup_{t\in[0,T)}\Vert\mathbf A_{f(t)}\Vert_{L^\infty(\Sigma)}=\infty,
  \end{equation}
  then there exist $t_j\nearrow T$, $r_j\searrow0$, $p_j\in\Sigma$, and a smooth proper Willmore immersion $\hat f\colon\hat\Sigma\to\R^3$ of a complete connected oriented surface such that, for some $\hat p\in\hat\Sigma$,
  \begin{equation}\label{eq:intro-pointed-convergence-blow-ups}
    \Big(\frac{f(t_j)-f(t_j,p_j)}{r_j},p_j\Big)\longrightarrow(\hat f,\hat p)
  \end{equation}
  in the sense of \Cref{def:local-smooth-pointed-convergence}, with $|\mathbf A_{\hat f}(\hat p)|>0$. Moreover,
  \begin{equation}\label{eq:intro-properties-stationary-blow-up}
    \int_{\hat\Sigma}|\mathbf A_{\hat f}|^2\dd\vol_{g_{\hat f}}<\infty\quad\text{and}\quad\mathrm{genus}(\hat\Sigma)\leq\mathrm{genus}(\Sigma).
  \end{equation}
  Either $\partial\hat\Sigma=\emptyset$, or, after a Euclidean isometry, $\hat f\vert_{\partial\hat\Sigma}$ is an embedding onto $\R e_1$ and $\partial_{\eta_{\hat f}}\hat f=-e_2$ on $\partial\hat\Sigma$, where $e_1,e_2,e_3$ are the standard coordinate vectors. The density at infinity
  \begin{equation}\label{eq:intro-density-at-infinity}
    \theta^2(\hat f,\infty)\vcentcolon=\lim_{R\to\infty}\frac{\int_{\hat f^{-1}(B_R(0))}\dd\vol_{g_{\hat f}}}{\pi R^2}
  \end{equation}
  exists and satisfies
  \begin{equation}\label{eq:intro-density-bound-stationary-blow-up}
    \theta^2(\hat f,\infty)+\frac{\W(\hat f)}{4\pi}\leq\frac{\W(f_0)}{4\pi}+\begin{cases}
      \displaystyle\frac1{2\pi}\sup_{z\in\R^3}\mathbf B_{\Gamma,\mathbf N}(z)&\text{if $\partial\hat\Sigma=\emptyset$},\\
      \displaystyle\frac1{2\pi}\sup_{z\in\Gamma}\mathbf B_{\Gamma,\mathbf N}(z)&\text{if $\partial\hat\Sigma\neq\emptyset$}.
    \end{cases}
  \end{equation}
\end{theorem}

The estimate \eqref{eq:intro-density-bound-stationary-blow-up} explains the energy threshold \eqref{eq:1-thm:main-result}. Assuming for the moment that \eqref{eq:1-thm:main-result} holds strictly, we see that noncompact blow-ups $\hat f$ without boundary have density at infinity strictly below two, and those with boundary have density at infinity strictly below one. These bounds are crucial for classifying (and thus excluding) them. Both estimates follow from Simon's monotonicity formula. The estimate is not new in the case without boundary (it is a direct analog of \cite[Equation~(A.23)]{kuwertschaetzle2004}). However, to show that one only needs to take the supremum over all $z\in\Gamma$, rather than over all $z\in\R^3$ in the case $\partial\hat\Sigma\neq\emptyset$, we carefully analyze the jump behavior of $\mathbf B_{\Gamma,\mathbf N}$, see \Cref{lem:boundary-integral-jump}.

The main analytic step in the proof of \Cref{thm:main-result} is to obtain a lifespan estimate and parabolic smoothing estimates that are uniform under a smallness assumption on
\begin{equation}
  \Vert\mathbf A_{f_0}\Vert_{L^\infty(\Sigma)}+\Vert\bm\kappa_\Gamma\Vert_{L^\infty(\Gamma)}
\end{equation}
and a positive lower bound for $\mathfrak r_\Gamma$. Here $\kappa_\Gamma$ is the curvature vector of $\Gamma$, and $\mathfrak r_\Gamma$ relates its intrinsic distance to the ambient Euclidean distance, see \Cref{def:rGamma}. These assumptions allow us to represent $f_0$ locally as a graph over its tangent plane on universally chosen domains. Near the boundary, the domains of these graphs have uniformly controlled curved edges. We then apply classical theory for quasiconformal maps \cite{ahlforsbers1960} to construct isothermal coordinates on fixed interior and boundary domains. Since the graphs are uniformly small in $W^{2,\infty}$, classical Schauder theory and elliptic regularity yield uniform smooth control of the harmonic coordinate changes. This produces an atlas and a subordinate localization system whose constants depend only on these geometric bounds, rather than on the particular parametrization of the immersion, or on its induced metric.

In these coordinates, the localized flow is treated as a fourth-order parabolic boundary value problem, which is only a small perturbation of the biharmonic heat flow. We use time-weighted maximal $L^p$-regularity theory \cite{pruesssimonett2004,denkhieberpruess2007} for homogeneous boundary problems, following the approach in \cite{lecronepruesswilke2014}, to obtain a uniform existence time and smoothing estimates from small $C^{1+\alpha}$ graph norms, with $0<\alpha<1$. The uniform atlas allows us to combine the local estimates into the geometric lifespan and smoothing theorem, see \Cref{thm:lifespan-and-smoothing}. 

Our new geometric localization approach is in fact rather flexible. We expect that it can also be used for other types of boundary conditions, or in the closed case when additional constraints are imposed (see \cite{rupp2024,ruppscharrerschlierf2026,scharrerwest2025}), so long as the resulting localized PDEs in the coordinates of the isothermal atlas are again well-posed in the setting of time-weighted maximal $L^p$-regularity.

The localization approach can also be viewed as a quantitative analog of the theory of \emph{uniformly regular Riemannian manifolds}, for instance see \cite[Definition~3.3]{amann2025}. This theory was applied in \cite{lecroneshaosimonett2020} to the Willmore flow of surfaces without boundary, including noncompact surfaces. In contrast to previous applications, our focus lies on obtaining uniform quantitative estimates which only depend on geometric data, such as $\mathbf A_f$, $\bm\kappa_\Gamma$ and $\mathfrak r_\Gamma$. 

The second ingredient concerns classifying noncompact blow-ups with boundary. 

\begin{proposition}[Rigidity of boundary blow-ups]\label{prop:intro-boundary-rigidity}
  Let $\hat f\colon\hat\Sigma\to\R^3$ be a blow-up obtained by the construction in \Cref{thm:intro-stationary-blow-ups} with $\partial\hat\Sigma\neq\emptyset$. Then $\mathbf H_{\hat f}=0$ on $\partial\hat\Sigma$, and reflection about the line $\R e_1$ extends $\hat f$ to a smooth Willmore immersion of a complete surface without boundary. If, in addition, $\mathrm{genus}(\hat\Sigma)=0$ and $\theta^2(\hat f,\infty)<1$, then $\hat f$ parametrizes a flat half-plane.
\end{proposition}

Here reflection about $\R e_1$ means the map $(x_1,x_2,x_3)\mapsto(x_1,-x_2,-x_3)$, applied to a second copy of the surface and joined along its boundary. 

For \Cref{prop:intro-boundary-rigidity}, we first prove that the blow-ups satisfy the additional boundary condition that $\mathbf H_{\hat f}=0$ on $\partial\hat\Sigma$. This follows from a Pohozaev-type identity, related to the identities used in \cite{palmer2000,dallacqua2012}, together with boundary $\varepsilon$-regularity estimates proved in this work. The decay obtained from finite total curvature controls the terms at infinity. Then \cite{kuwertlamm2021} yields that the reflection is a smooth complete Willmore immersion without boundary, to which the classification and removability theory \cite{kuwertschaetzle2004,bryant1984} can be applied.

\begin{remark}
  For the sake of a cleaner presentation, the geometric smallness assumption in the above sketch of the construction of the isothermal atlas with uniform estimates is given in terms of the second-order data $\mathbf A_f$ and $\bm\kappa_\Gamma$. However, we actually require only $C^{1+\alpha}$-type geometric control. More precisely, instead of imposing smallness of $\mathbf A_f$, we control $[\mathbf n_f]_{\alpha,1}$, a suitably defined Hölder seminorm of the unit normal. For more details, see \Cref{app:graph-rep}. Similarly, instead of $\bm\kappa_\Gamma$, we use $[\bm\tau_\Gamma]_{C^\alpha(\Gamma)}$ where $\bm\tau_\Gamma$ is a unit tangent field along $\Gamma$.
\end{remark}

\subsection{Structure of the article}

In \Cref{sec:localization-approach}, we establish the main tools necessary to carry out the new localization approach sketched above. That is, for an immersion $f\in C^{1+\alpha}(\Sigma)$ with small Hölder seminorm of the unit normal $\mathbf n_f$ and boundary unit tangent $\bm\tau_\Gamma$, and positive lower bound for $\mathfrak r_\Gamma$, we first construct an isothermal atlas with uniform estimates depending only on these quantities. The main result is \Cref{prop:ex-good-atlas}. In \Cref{prop:smooth-clamped-immersion}, we construct a smooth reference surface $G_f$ and a diffeomorphism $\Xi_f\colon\Sigma\to\Sigma$ such that 
\begin{equation}
  G_f=f\circ\Xi_f \quad\text{and}\quad d(G_f-f\circ\Xi_f)=0\quad\text{on $\partial\Sigma$}.
\end{equation}
This is stronger than merely requiring that $G_f$ and $f\circ\Xi_f$ both satisfy \eqref{eq:Dirichlet-boundary-conditions}. Moreover, $\Vert \nabla^m_{g_{G_f}}\mathbf A_{G_f}\Vert_{L^\infty(\Sigma)}$ is uniformly controlled in terms of $m$, the Hölder seminorm bound for $\mathbf n_f$, $\mathfrak r_\Gamma$, and higher-order norms of the boundary data $\bm\kappa_\Gamma$ and $D_s\mathbf N$.

Then, \Cref{sec:lifespan-estimate} contains the proof of the lifespan and parabolic smoothing result. The argument combines ideas of \cite{lecronepruesswilke2014} with the localization approach based on the isothermal atlas constructed in \Cref{sec:localization-approach}. For an introduction to the time-weighted parabolic maximal $L^p$-regularity spaces and linear theory we use, see \Cref{app:max-reg-unif-surfaces}.

Finally, in \Cref{sec:pf-intro-results} we prove \Cref{thm:main-result,thm:intro-stationary-blow-ups,prop:intro-boundary-rigidity}, and in particular the $\varepsilon$-regularity result for noncompact blow-ups in \Cref{prop:eps-bdry-reg}.

\section{Isothermal atlas and approximation of $C^{1+\alpha}$-immersions with estimates}
\label{sec:localization-approach}

Fix $\alpha\in(0,1)$. Throughout this section, constants are called \emph{universal} if they depend at most on $\alpha$. Moreover, $\Sigma$ denotes a compact orientable surface with non-empty boundary, and we write $Q_r=(-r,r)^2$ and $Q_r^+=(-r,r)\times[0,r)$ for $r>0$.

\subsection{Construction of an isothermal atlas with uniform estimates}

\begin{proposition}\label{prop:ex-good-atlas}
  There exist universal constants $\varepsilon_{\ref{prop:ex-good-atlas}}>0$ and $C_{\ref{prop:ex-good-atlas}}>0$ with the following property. Let $f\in C^{1+\alpha}(\Sigma,\R^3)$ be an immersion with $\Gamma=f(\partial\Sigma)$ a smooth submanifold of $\R^3$ and, using the notation of \Cref{def:rGamma} and \eqref{eq:def-nf-alpha-seminorm}, assume that 
  \begin{equation}\label{eq:1-prop:ex-good-atlas}
    \mathfrak r_\Gamma>16\quad\text{and}\quad[\mathbf n_f]_{\alpha,1} + [\bm\tau_\Gamma]_{C^\alpha(\Gamma)} \leq\varepsilon \leq \varepsilon_{\ref{prop:ex-good-atlas}}.
  \end{equation}
  There exist a universal constant $\rho^\ast\in(0,1/16)$, $\mathbf I_{\mathrm{int}}=\{1,\dots,N_1\}$ and $\mathbf I_{\mathrm{bdr}}=\{N_1+1,\dots,N_1+N_2\}$, and charts
  \begin{equation}
    \alpha_j\colon D_j\vcentcolon=Q_{\rho^\ast}\to \Omega_j\text{ for $j\in \mathbf I_{\mathrm{int}}$},\quad \alpha_j\colon D_j\vcentcolon=Q_{16\rho^\ast}^+ \to\Omega_j\text{ for $j\in \mathbf I_{\mathrm{bdr}}$}
  \end{equation}
  with $\varphi_j$ denoting their inverses such that 
  \begin{equation}\label{eq:2-prop:ex-good-atlas}
    1\leq \sum_{j\in \mathbf I_{\mathrm{int}}}\chi_{\alpha_j(Q_{\rho^\ast/2})} + \sum_{j\in \mathbf I_{\mathrm{bdr}}}\chi_{\alpha_j(Q_{12\rho^\ast}^+)} \leq \sum_{j\in \mathbf I_{\mathrm{int}}\cup \mathbf I_{\mathrm{bdr}}} \chi_{\Omega_j}\leq Z^\ast,\quad \sup_{\varphi_k(\Omega_j\cap\Omega_k)} \vert \partial^\beta(\varphi_j\circ\alpha_k)\vert\leq C_{\vert\beta\vert}
  \end{equation}
  for a universal constant $Z^\ast>0$ and constants $C_{\vert\beta\vert}$ only depending on the order of the multi-index $\beta\in\N_0^2$. Moreover, there exists $\lambda_j\in C^{\alpha}(D_j)$ with 
  \begin{equation}\label{eq:3-prop:ex-good-atlas}
    \alpha_j^*g_f = e^{2\lambda_j}g_{\R^2},\quad \Vert \lambda_j\Vert_{C^\alpha(D_j)}\leq C_{\ref{prop:ex-good-atlas}}\varepsilon.
  \end{equation}
  Finally, there exist $C^{1+\alpha}$-maps $\omega_j\colon D_j\to\R^2$, $h_j\colon D_j\to\R$ such that
  \begin{equation}\label{eq:almost-graph-charts}
    f\circ\alpha_j=A_j(\omega_j,h_j)\text{ on $D_j$}\quad\text{and}\quad \Vert \omega_j-\operatorname{id}\Vert_{C^{1+\alpha}(D_j)} + \Vert h_j\Vert_{C^{1+\alpha}(D_j)} \leq C_{\ref{prop:ex-good-atlas}}\varepsilon.
  \end{equation}
  Here, $A_j$ are Euclidean isometries with $A_j(0)=f(p_j)$, $DA_j.e_3=\mathbf n_f(p_j)$ for $j\in \mathbf I_{\mathrm{int}}\cup \mathbf I_{\mathrm{bdr}}$ where $p_j\in\Sigma\setminus\partial\Sigma$ for $j\in \mathbf I_{\mathrm{int}}$, and $p_j\in\partial\Sigma$ for $j\in \mathbf I_{\mathrm{bdr}}$, and additionally $DA_j.e_2=-\partial_\eta f(p_j)$ for $j\in \mathbf I_{\mathrm{bdr}}$.
\end{proposition}

The fundamental strategy in constructing the isothermal atlas uses two ingredients: First we apply \Cref{prop:ex-of-graph-charts} to get a covering of $\Sigma$ by $C^{1+\alpha}
$-small graph charts, and then apply classical existence theory for quasiconformal maps by solving \emph{Beltrami's equation} in these charts. Each step is contained in a separate lemma below.

\begin{lemma}\label{lem:greedy-selection-of-Urps}
  Suppose that the assumptions of \Cref{prop:ex-good-atlas} hold. For $\varepsilon_{\ref{prop:ex-good-atlas}}$ sufficiently small, there exist $\mathbf I_{\mathrm{int}}=\{1,\dots,N_1\}$, $\mathbf I_{\mathrm{bdr}}=\{N_1+1,\dots,N_1+N_2\}$, and points $p_j\in \Sigma$ with $d_{g_f}(p_j,\partial\Sigma)>4\sqrt 2\rho^\ast$ for $j\in \mathbf I_{\mathrm{int}}$ and $p_j\in\partial\Sigma$ for $j\in \mathbf I_{\mathrm{bdr}}$ such that
  \begin{equation}\label{eq:1-lem:greedy-selection-of-Urps}
    1\leq \sum_{j\in \mathbf I_{\mathrm{int}}}\chi_{U_{{\rho^\ast}/{4},p_j}} + \sum_{j\in \mathbf I_{\mathrm{bdr}}}\chi_{U_{8\sqrt 2\rho^\ast,p_j}} \leq \sum_{j\in \mathbf I_{\mathrm{int}}}\chi_{U_{2\rho^\ast,p_j}} + \sum_{j\in \mathbf I_{\mathrm{bdr}}}\chi_{U_{32\rho^\ast,p_j}} \leq Z^\ast.
  \end{equation}
\end{lemma}
\begin{proof}
  Suppose that $\varepsilon_{\ref{prop:ex-good-atlas}}$ is chosen sufficiently small such that
  $\varepsilon_{\ref{prop:ex-good-atlas}}\leq\varepsilon_{\ref{prop:ex-of-graph-charts}}$ and $C_{\ref{prop:ex-of-graph-charts}}\varepsilon_{\ref{prop:ex-good-atlas}}\leq 2-\sqrt2$.
  Since $U_{r,p}\subset U_{1,p}$ for $r\leq1$, the assumptions of \Cref{prop:ex-of-graph-charts} are satisfied at each $r$ if $p\in\partial\Sigma$ or $d_{g_f}(p,\partial\Sigma)\geq2\sqrt2r$.
  In these cases, \eqref{eq:4-prop:ex-of-graph-charts} and \eqref{eq:5-prop:ex-of-graph-charts} yield that
  \begin{equation}\label{eq:1-pf-lem:greedy-selection-of-Urps}
    B_{d_{g_f}}(p,r)\subset U_{r,p}\subset B_{d_{g_f}}(p,2r)
  \quad\text{and}\quad r^2\leq\int_{U_{r,p}}\dd\vol_{g_f}\leq5r^2.
  \end{equation}
  The universal constant $r^\ast\in(0,1)$ is fixed at a later stage in the proof. To keep the notation more readable, we abbreviate
  \begin{equation}    
    \rho^\ast=\frac{r^\ast}{2\sqrt2},
    \quad s^\ast=\frac{\rho^\ast}{4}=\frac{r^\ast}{8\sqrt2},
    \quad h^\ast=\frac{s^\ast}{4}.
  \end{equation}
  Choose successively $q_k\in\partial\Sigma\setminus\bigcup_{j<k}U_{r^\ast,q_j}$ until $\partial\Sigma\subset\bigcup_{j=1}^{N_2}U_{r^\ast,q_j}$. Using \eqref{eq:1-pf-lem:greedy-selection-of-Urps}, this process actually terminates at $N_2\in\N$ since $\Sigma$ and $\partial\Sigma$ are compact. Let $K\vcentcolon=\Sigma\setminus\bigcup_{j=1}^{N_2}U_{4r^\ast,q_j}$. For $p\in K$, choose a nearest boundary point $q$ and an
  index $j$ with $q\in U_{r^\ast,q_j}$. By \eqref{eq:1-pf-lem:greedy-selection-of-Urps}, we have
  \begin{equation}\label{eq:3-pf-lem:greedy-selection-of-Urps}
    d_{g_f}(p,\partial\Sigma) = d_{g_f}(p,q)\geq d(p,q_j)-d(q_j,q)>4r^\ast-2r^\ast = 2r^\ast = 4\sqrt2\rho^\ast\quad\text{for all $p\in K$}.
  \end{equation}
  In particular, \Cref{prop:ex-of-graph-charts} applies at every $p\in K$ with any radius $0<r\leq 2\rho^\ast$. Now choose successively $p_k\in K\setminus\bigcup_{j<k}U_{s^\ast,p_j}$. 
  Again, using \eqref{eq:1-pf-lem:greedy-selection-of-Urps} and compactness of $\Sigma$, this process terminates after $N_1\in\N$ steps, and using $4r^\ast=8\sqrt2\rho^\ast$, we get the first inequality in \eqref{eq:1-lem:greedy-selection-of-Urps}. For the second inequality, it only remains to show that $Z^\ast$ can be chosen universally. To this end, first note that the above construction and \eqref{eq:1-pf-lem:greedy-selection-of-Urps} together yield that
  \begin{equation}
    d_{g_f}(q_k,q_j)>s^\ast\text{ for all $k\neq j$}\quad\text{and}\quad d_{g_f}(p_k,p_j)>s^\ast\text{ for all $k\neq j$}.
  \end{equation}
  Moreover, for $p_j\notin U_{4r^\ast,q_k}$, \eqref{eq:1-pf-lem:greedy-selection-of-Urps} yields that $d_{g_f}(p_j,q_k)\geq4r^\ast>s^\ast$. Combined with \eqref{eq:1-pf-lem:greedy-selection-of-Urps}, we get that 
  \begin{equation}\label{eq:2-pf-lem:greedy-selection-of-Urps}
    \text{all sets in the family $\{U_{h^\ast,q_j},U_{h^\ast,p_k}\}_{j,k}$ are pairwise disjoint.}
  \end{equation}  
  Now suppose that for $J_{\mathrm{int}}\subset\{1,\dots,N_1\}$, $J_{\mathrm{bdr}}\subset\{1,\dots,N_2\}$, and $x\in\Sigma$, we have
  \begin{equation}
    x\in U_{2\rho^\ast,p_j}\text{ for all $j\in J_{\mathrm{int}}$}\quad\text{and}\quad x\in U_{32\rho^\ast,q_j}\text{ for all $j\in J_{\mathrm{bdr}}$}.
  \end{equation}
  Suppose that $r^\ast$ is sufficiently small such that $32\rho^\ast<1$. By \eqref{eq:3-pf-lem:greedy-selection-of-Urps} and \eqref{eq:1-pf-lem:greedy-selection-of-Urps} and the triangle inequality, there is a universal constant $C>0$ such that, for any $c\in\{p_j,q_k\}_{j\in J_{\mathrm{int}},k\in J_{\mathrm{bdr}}}$, we get
  \begin{equation}
    \bigcup_{j\in J_{\mathrm{int}}} U_{2\rho^\ast,p_j} \,\cup\,\bigcup_{j\in J_{\mathrm{bdr}}}U_{32\rho^\ast,q_j} \subset B_{d_{g_f}}(c,Cr^\ast).
  \end{equation}
  Further reduce $r^\ast$ so that $Cr^\ast<1$, and fix any of these centers $c$. If $d_{g_f}(c,\partial\Sigma)\geq2\sqrt2Cr^\ast$, we have $B_{d_{g_f}}(c,Cr^\ast)\subset U_{Cr^\ast,c}$ by \eqref{eq:1-pf-lem:greedy-selection-of-Urps}. Otherwise, choose a nearest boundary point $q$ and a suitable $q_j$ with $q\in U_{r^\ast,q_j}$. Then \eqref{eq:1-pf-lem:greedy-selection-of-Urps} yields that $d_{g_f}(c,q_j)<(2\sqrt2 C +2)r^\ast$ and therefore, writing $C'=C+(2\sqrt 2C+2)$, and finally fixing $r^\ast$ sufficiently small such that also $C'r^\ast < 1$, we get
  \begin{equation}
    \bigcup_{j\in J_{\mathrm{int}}} U_{2\rho^\ast,p_j} \,\cup\,\bigcup_{j\in J_{\mathrm{bdr}}}U_{32\rho^\ast,q_j} \subset B_{d_{g_f}}\bigl(q_j,C' r^\ast\bigr)
    \subset U_{C' r^\ast,q_j}.
  \end{equation}
  Using \eqref{eq:2-pf-lem:greedy-selection-of-Urps}, for $m=\# J_{\mathrm{int}}+\# J_{\mathrm{bdr}}$, we get 
  \begin{equation}
    m(h^\ast)^2 \leq \int_{\bigcup_{j\in J_{\mathrm{int}}} U_{2\rho^\ast,p_j} \,\cup\,\bigcup_{j\in J_{\mathrm{bdr}}}U_{32\rho^\ast,q_j}}\dd\vol_{g_f} \leq 5(C'r^\ast)^2.
  \end{equation}
  This shows that we can take $Z^\ast$ universally in \eqref{eq:1-lem:greedy-selection-of-Urps}. 
\end{proof}

\begin{lemma}[Interior isothermal charts]\label{lem:interior-isothermal-charts}
  Work in the setting of \Cref{lem:greedy-selection-of-Urps}, and fix its universal constant
  $\rho^\ast$. After reducing $\varepsilon_{\ref{prop:ex-good-atlas}}$ universally, there exist orientation-preserving $C^{1+\alpha}$-charts
  $\alpha_j\colon Q_{\rho^\ast}\to\Omega_j$, $j\in \mathbf I_{\mathrm{int}}$,
  with $\alpha_j(0)=p_j$ and inverses $\varphi_j$, such that
  \begin{equation}\label{eq:1-lem:interior-isothermal-charts}
    U_{\rho^\ast/4,p_j}
    \subset\alpha_j(Q_{\rho^\ast/2})
    \subset\Omega_j
    \subset U_{2\rho^\ast,p_j}.
  \end{equation}
  Moreover, there exist $\lambda_j\in C^\alpha(Q_{\rho^\ast})$,
  $\omega_j\in C^{1+\alpha}(Q_{\rho^\ast},\R^2)$ and
  $h_j\in C^{1+\alpha}(Q_{\rho^\ast},\R)$ such that
  \begin{equation}\label{eq:2-lem:interior-isothermal-charts}
    \alpha_j^*g_f=e^{2\lambda_j}g_{\R^2},
    \qquad f\circ\alpha_j=A_{p_j}(\omega_j,h_j),
  \end{equation}
  and
  \begin{equation}\label{eq:3-lem:interior-isothermal-charts}
    \Vert \lambda_j\Vert _{C^\alpha(Q_{\rho^\ast})}
    +\Vert \omega_j-\operatorname{id}\Vert _{C^{1+\alpha}(Q_{\rho^\ast})}
    \leq C\varepsilon^2,
    \qquad
    \Vert h_j\Vert _{C^{1+\alpha}(Q_{\rho^\ast})}\leq C\varepsilon.
  \end{equation}
  Here, $C$ is universal. Moreover, the transition functions $T_{jk}\colon\varphi_k(\Omega_j\cap\Omega_k)\to \varphi_j(\Omega_j\cap\Omega_k)$ are holomorphic, where $\varphi_k(\Omega_j\cap\Omega_k)$ is open in $\C$, and for a universal constant $C_1>0$, we have 
  \begin{equation}\label{eq:4.1-lem:interior-isothermal-charts}
    \Vert T_{jk}\Vert_{C^0(\varphi_k(\Omega_j\cap\Omega_k))} + \Vert T_{jk}'\Vert_{C^0(\varphi_k(\Omega_j\cap\Omega_k))} \leq C_1.
  \end{equation}
  Furthermore, for every $l\in\N$ with $l\geq 2$, there is a constant
  $C_l$, universal up to dependence on $l$, such that, writing $z=x+iy$ for the complex coordinate, 
  \begin{equation}\label{eq:4-lem:interior-isothermal-charts}
    \Vert \partial_z^l T_{jk}\Vert_{C^0(\varphi_k(\Omega_j\cap\Omega_k))}\leq C_l \varepsilon
    \quad\text{for all $j,k\in \mathbf I_{\mathrm{int}}$ with
    $\Omega_j\cap\Omega_k\neq\emptyset$.}
  \end{equation}
\end{lemma}
\begin{proof}
  By \Cref{lem:greedy-selection-of-Urps}, for $j\in \mathbf I_{\mathrm{int}}$, \Cref{prop:ex-of-graph-charts} applies at radius $r=2\rho^\ast$ and yields
  \begin{equation}\label{eq:1-pf:interior-isothermal-charts}
    \Psi_j\vcentcolon=
    (P_{p_j}|_{U_{2\rho^\ast,p_j}})^{-1}\colon Q_{2\rho^\ast}\to U_{2\rho^\ast,p_j} \quad\text{and}\quad f\circ\Psi_j(x)=A_{p_j}(x,u_{p_j}(x)).
  \end{equation}
  Let $g_j=\Psi_j^*g_f$. Its coefficients satisfy
  $(g_j)_{ab}=\delta_{ab}+\partial_a u_{p_j}\partial_b u_{p_j}$. Identify $\R^2$ with $\C$ by $z=x^{(1)}+ix^{(2)}$, and use $\partial_z=\frac12(\partial_1-i\partial_2)$ and $\partial_{\bar z}=\frac12(\partial_1+i\partial_2)$. Set
  \begin{equation}\label{eq:3-pf:interior-isothermal-charts}
    \mu_j\vcentcolon=
    \frac{(g_j)_{11}-(g_j)_{22}+2i(g_j)_{12}}
    {\tr g_j+2\sqrt{\det g_j}} \quad\text{and}\quad a_j\vcentcolon=\frac{\tr g_j+2\sqrt{\det g_j}}4.
  \end{equation}
  A direct computation yields $g_j=a_j\vert\mathrm d z+\mu_j\dd\bar z\vert^2$. Since $g_j-I$ is quadratic in $Du_{p_j}$, \eqref{eq:3-prop:ex-of-graph-charts} yields
  \begin{equation}\label{eq:5-pf:interior-isothermal-charts}
    \Vert \mu_j\Vert _{C^\alpha(Q_{2\rho^\ast})}
    +\Vert a_j-1\Vert _{C^\alpha(Q_{2\rho^\ast})}
    \leq C\varepsilon^2.
  \end{equation}
  Extend $\mu_j$ to a function $\nu_j\in C^\alpha_c(\C)$ such that
  \begin{equation}
    \nu_j=\mu_j\quad\text{on $Q_{2\rho^\ast}$}, \qquad \operatorname{supp}\nu_j\subset Q_{3\rho^\ast}, \quad\text{and}\quad \Vert \nu_j\Vert _{C^\alpha(\C)} \leq C\varepsilon^2.
  \end{equation}
  Suppose that $\varepsilon_{\ref{prop:ex-good-atlas}}$ is sufficiently small such that $\Vert \nu_j\Vert _{L^\infty(\C)}\leq 1/2$. We next argue that there exists $\Phi_j\in C^{1+\alpha}(\C)$ uniformly close to the identity solving the associated Beltrami equation, i.e.\
  \begin{equation}\label{eq:6-pf:interior-isothermal-charts}
    \Phi_j(0)=0,\qquad \partial_{\bar z}\Phi_j=\nu_j\partial_z\Phi_j, \quad\text{and}\quad \Vert \Phi_j-\operatorname{id}\Vert_{C^{1+\alpha}(\C)}
    \leq C\varepsilon^2.
  \end{equation}
  To this end, by \cite[Theorems~14.0.4 and~14.2.2]{astalaiwaniecmartin2009} applied with $p=2$, there exists a solution $w_j$, unique up to an additive constant, satisfying
  \begin{equation}\label{eq:7-pf:interior-isothermal-charts}
    \partial_{\bar z}w_j-\nu_j\partial_z w_j=\nu_j
    \quad\text{on $\C$} \quad\text{and}\quad \Vert Dw_j\Vert _{L^2(\C)}
    \leq C\Vert \nu_j\Vert _{L^2(\C)}
    \leq C\varepsilon^2.
  \end{equation}
  Since $\operatorname{supp}\nu_j\subset Q_{3\rho^\ast}$, we also have $\Vert \nu_j\Vert_{L^1(\C)}\leq C\varepsilon^2$. Choosing the additive constant as in \cite[Lemma~15.5.1, with $s=1$]{astalaiwaniecmartin2009}, we obtain
  \begin{equation}
    \Vert w_j\Vert _{L^\infty(\C)}
    \leq C
    \Vert \nu_j\Vert _{L^1(\C)}^{\frac{1+\alpha}{2+\alpha}}
    [\nu_j]_{C^\alpha(\C)}^{\frac{1}{2+\alpha}}
    \leq C\varepsilon^2.
  \end{equation}
  Applying \cite[Theorem~15.0.6, with $s=2$]{astalaiwaniecmartin2009} to \eqref{eq:7-pf:interior-isothermal-charts} yields
  \begin{equation}
    [Dw_j]_{C^\alpha(\C)}
    \leq C\Bigl(
      [\nu_j]_{C^\alpha(\C)}
      +[\nu_j]_{C^\alpha(\C)}^{1+1/\alpha}
       \Vert \nu_j\Vert _{L^2(\C)}
    \Bigr)
    \leq C\varepsilon^2.
  \end{equation}
  Together with \cite[Lemma~15.4.1]{astalaiwaniecmartin2009}, we get
  \begin{equation}
    \Vert Dw_j\Vert _{L^\infty(\C)}
    \leq C\bigl(
      \Vert Dw_j\Vert _{L^2(\C)}+[Dw_j]_{C^\alpha(\C)}
    \bigr)
    \leq C\varepsilon^2.
  \end{equation}
  Altogether, taking $\Phi_j=\mathrm{id}+w_j-w_j(0)$ yields \eqref{eq:6-pf:interior-isothermal-charts}. Further reducing $\varepsilon_{\ref{prop:ex-good-atlas}}$, \eqref{eq:6-pf:interior-isothermal-charts} and the inverse function theorem show that $\Phi_j$ is a $C^{1+\alpha}$-diffeomorphism such that
  \begin{equation}\label{eq:14-pf:interior-isothermal-charts}
    \Vert \Phi_j-\operatorname{id}\Vert _{C^0(\C)}
    +\Vert \Phi_j^{-1}-\operatorname{id}\Vert _{C^0(\C)}\leq\frac{\rho^\ast}8,
  \end{equation}
  and thus
  \begin{equation}\label{eq:15-pf:interior-isothermal-charts}
    \Phi_j(Q_{\rho^\ast/4})\subset Q_{3\rho^\ast/8}\subset Q_{\rho^\ast/2} \quad\text{and}\quad \Phi_j^{-1}(Q_{3\rho^\ast/2})\subset Q_{13\rho^\ast/8}\Subset Q_{2\rho^\ast}.
  \end{equation}
  Define the chart $\alpha_j \vcentcolon= \Psi_j\circ\Phi_j^{-1}$ and let $\Omega_j\vcentcolon=\alpha_j(Q_{\rho^\ast})$. Its inverse on $\Omega_j$ is $\varphi_j=\Phi_j\circ P_{p_j}|_{\Omega_j}$. Equations~\eqref{eq:14-pf:interior-isothermal-charts} and \eqref{eq:15-pf:interior-isothermal-charts}, together with $\Psi_j(Q_s)=U_{s,p_j}$, imply \eqref{eq:1-lem:interior-isothermal-charts}, where the upper inclusion uses $\Phi_j^{-1}(Q_\rho^\ast)\subset Q_{9\rho^\ast/8}\subset Q_{3\rho^\ast/2}$. By \eqref{eq:6-pf:interior-isothermal-charts}, we have
  \begin{equation}\label{eq:17-pf:interior-isothermal-charts}
    \Phi_j^*g_{\R^2} = |\partial_z\Phi_j|^2|\dd z+\mu_j\dd\bar z|^2 \quad\text{and}\quad g_j=\frac{a_j}{|\partial_z\Phi_j|^2}\Phi_j^*g_{\R^2}
    \quad\text{on $Q_{2\rho^\ast}$.}
  \end{equation}
  Setting $\lambda_j=(\frac12\log a_j-\log|\partial_z\Phi_j|) \circ\Phi_j^{-1}$, $\omega_j=\Phi_j^{-1}$, and $h_j=u_{p_j}\circ\omega_j$, we get \eqref{eq:2-lem:interior-isothermal-charts} and \eqref{eq:3-lem:interior-isothermal-charts}.

  It remains to prove \eqref{eq:4-lem:interior-isothermal-charts}.
  Write $\widehat\Omega_j=\alpha_j(Q_{3\rho^\ast/2})$.
  By \eqref{eq:15-pf:interior-isothermal-charts}, the chart $\alpha_j$
  is defined on a neighborhood of $\overline{Q_{3\rho^\ast/2}}$.
  Reducing $\varepsilon_{\ref{prop:ex-good-atlas}}$ if necessary, we have
  $1/2\leq e^{\lambda_j}\leq2$ there. Every path from a point of $\Omega_j$ to the complement of
  $\widehat\Omega_j$ must first cross the image of $\partial Q_{3\rho^\ast/2}$. Its portion inside $\widehat\Omega_j$ has length at least half the Euclidean distance from $Q_\rho^\ast$ to $\partial Q_{3\rho^\ast/2}$. Hence
  \begin{equation}\label{eq:21-pf:interior-isothermal-charts}
    d_{g_f}(\Omega_j,\Sigma\setminus\widehat\Omega_j) \geq\frac{\rho^\ast}4.
  \end{equation}
  If $x\in\Omega_j\cap\Omega_k$ and $w=\varphi_k(x)$, then $B_{\rho^\ast/16}(w)\subset Q_{3\rho^\ast/2}$, and the image under $\alpha_k$ of a line segment in this disk starting at $w$ has length less than $\rho^\ast/8$. By \eqref{eq:21-pf:interior-isothermal-charts}, it stays in $\widehat\Omega_j$. Consequently, for $r^\ast=\rho^\ast/16$,
  \begin{equation}\label{eq:18-pf:interior-isothermal-charts}
    T_{jk}\vcentcolon=\alpha_j^{-1}\circ\alpha_k
    \colon B_{r^\ast}(w)\to Q_{3\rho^\ast/2}
  \end{equation}
  is well-defined. Since we have 
  \begin{equation}\label{eq:22-pf:interior-isothermal-charts}
    T_{jk}^*g_{\R^2} = e^{2(\lambda_k-\lambda_j\circ T_{jk})}g_{\R^2}
  \end{equation}
  we get for the complex derivative $T_{jk}'$, on $B_{r^\ast}(w)$ that
  \begin{equation}\label{eq:23-pf:interior-isothermal-charts}
    \log|T_{jk}'|=\lambda_k-\lambda_j\circ T_{jk}, \quad\text{so using \eqref{eq:3-lem:interior-isothermal-charts},}\quad |\log|T_{jk}'||\leq C\varepsilon.
  \end{equation}
  On $B_{r^\ast}(w)$, choose a holomorphic branch of the logarithm $q_{jk}=\log T_{jk}'$. By \eqref{eq:23-pf:interior-isothermal-charts}, we have $|\operatorname{Re}q_{jk}|\leq C\varepsilon$. Since $\operatorname{Re}q_{jk}$ is harmonic, elliptic interior estimates yield for all $l\geq 0$ that
  \begin{equation}
    \Vert \operatorname{Re}q_{jk} \Vert_{C^l(B_{r^\ast/2}(w))} \leq C_l \varepsilon
  \end{equation}
  for constants $C_l$ which are universal up to $l$. As $q_{jk}$ is holomorphic, this bounds all derivatives of $\operatorname{Im}q_{jk}$ by Cauchy--Riemann, such that $\Vert q_{jk}'\Vert_{C^{l-1}(B_{r^\ast/2}(w))}\leq C_l\varepsilon$ for all $l\in\N$. Therefore, using $T_{jk}'=e^{q_{jk}}$, we obtain that, for all $l\geq 0$,
  \begin{equation}\label{eq:24-pf:interior-isothermal-charts}
    \Vert T_{jk}'-e^{i \operatorname{Im}q_{jk}(w)}\Vert_{C^{l}(B_{r^\ast/2}(w))}\leq C_l\varepsilon,
  \end{equation}
  which yields \eqref{eq:4-lem:interior-isothermal-charts}.
\end{proof}

\begin{lemma}[Boundary isothermal charts]\label{lem:boundary-isothermal-charts}
  Work in the setting of \Cref{lem:greedy-selection-of-Urps}, and fix its universal constant $\rho^\ast$. Choosing $\varepsilon_{\ref{prop:ex-good-atlas}}$ smaller if necessary, there exist $C^{1+\alpha}$-charts $\alpha_j\colon Q_{16\rho^\ast}^+\to\Omega_j$, $j\in \mathbf I_{\mathrm{bdr}}$, with $\alpha_j(0)=p_j$ and inverses $\varphi_j$, such that
  \begin{equation}\label{eq:1-lem:boundary-isothermal-charts}
    U_{8\sqrt2\rho^\ast,p_j}\subset\alpha_j(Q_{12\rho^\ast}^+)\subset\Omega_j\subset U_{32\rho^\ast,p_j}.
  \end{equation}
  Moreover, there exist $\lambda_j\in C^\alpha(Q_{16\rho^\ast}^+)$, $\omega_j\in C^{1+\alpha}(Q_{16\rho^\ast}^+,\R^2)$ and $h_j\in C^{1+\alpha}(Q_{16\rho^\ast}^+,\R)$ such that
  \begin{equation}\label{eq:2-lem:boundary-isothermal-charts}
    \alpha_j^*g_f=e^{2\lambda_j}g_{\R^2}\quad\text{and}\quad f\circ\alpha_j=A_{p_j}(\omega_j,h_j),
  \end{equation}
  and
  \begin{equation}\label{eq:3-lem:boundary-isothermal-charts}
    \Vert \lambda_j\Vert _{C^\alpha(Q_{16\rho^\ast}^+)}+\Vert \omega_j-\operatorname{id}\Vert _{C^{1+\alpha}(Q_{16\rho^\ast}^+)}+\Vert h_j\Vert _{C^{1+\alpha}(Q_{16\rho^\ast}^+)}\leq C\varepsilon.
  \end{equation}
  Here, $C$ is universal. Together with the charts of \Cref{lem:interior-isothermal-charts}, the transition functions $T_{jk}\colon\varphi_k(\Omega_j\cap\Omega_k)\to \varphi_j(\Omega_j\cap\Omega_k)$ are holomorphic, where $\varphi_k(\Omega_j\cap\Omega_k)$ open in $\C$ if $j$ or $k$ lies in $\mathbf I_{\mathrm{int}}$ resp.\ relatively open in the closed upper half-plane for $j,k\in \mathbf I_{\mathrm{bdr}}$, and we have the analog estimates to \eqref{eq:4.1-lem:interior-isothermal-charts} and \eqref{eq:4-lem:interior-isothermal-charts}.

  Finally, there exists a universal constant $r_{\ref{lem:boundary-isothermal-charts}}>0$ such that, if both $j,k\in \mathbf I_{\mathrm{bdr}}$, we have
  \begin{equation}\label{eq:4-lem:boundary-isothermal-charts}
    \sup_{\substack{z\in\varphi_k(\Omega_j\cap\Omega_k)\\0\leq\operatorname{Im}z\leq r_{\ref{lem:boundary-isothermal-charts}}}}\big|\partial_z^l(T_{jk}'(z)-1)\big|\leq C_l\varepsilon,
  \end{equation}
  for every $l\geq 0$ with a constant $C_l$ which is universal up to dependence on $l$.
\end{lemma}
\begin{proof}
  Fix $j\in \mathbf I_{\mathrm{bdr}}$. By \Cref{prop:ex-of-graph-charts}, we have
  \begin{equation}
    \Psi_j\vcentcolon=(P_{p_j}|_{U_{32\rho^\ast,p_j}})^{-1}\colon D_{32\rho^\ast,p_j}\to U_{32\rho^\ast,p_j}\quad\text{and}\quad f\circ\Psi_j(x)=A_{p_j}(x,u_{p_j}(x)).
  \end{equation}
  Write $D_{32\rho^\ast,p_j}=\{(x,y)\in Q_{32\rho^\ast}:y\geq\psi_j(x)\}$, where $\psi_j(0)=\psi_j'(0)=0$, and set
  \begin{equation}
    v_j(x)\vcentcolon=-\pi(DA_{p_j})^{-1}\partial_\eta f\bigl(\Psi_j(x,\psi_j(x))\bigr).
  \end{equation}
  The estimates \eqref{eq:3-prop:ex-of-graph-charts}, \eqref{eq:3-pf-prop:ex-of-graph-charts} on $u_{p_j}$ resp.\ $\psi_j$, and \eqref{eq:1-prop:ex-good-atlas} yield
  \begin{equation}
    \Vert \psi_j\Vert _{C^{1+\alpha}((-32\rho^\ast,32\rho^\ast))}+\Vert v_j-e_2\Vert _{C^\alpha((-32\rho^\ast,32\rho^\ast))}\leq C\varepsilon.
  \end{equation}
  Extend $v_j-e_2$ to $\R$ with the same bound up to a universal constant, and use the operator from \Cref{lem:whitney-extension} to define
  \begin{equation}\label{eq:1-pf:boundary-isothermal-charts}
    \phi_j(x,y)\vcentcolon=(x,\psi_j(x))+ye_2+\mathcal E(v_j-e_2)(x,y)\quad\text{for $(x,y)\in Q_{24\rho^\ast}^+$}.
  \end{equation}
  Then \eqref{eq:whitney-extension-traces} and \eqref{eq:whitney-extension-estimates} yield that
  \begin{equation}\label{eq:2-pf:boundary-isothermal-charts}
    \Vert \phi_j-\operatorname{id}\Vert_{C^{1+\alpha}(Q_{24\rho^\ast}^+)}\leq C\varepsilon,\quad \phi_j(x,0)=(x,\psi_j(x))\quad\text{and}\quad \partial_y\phi_j(x,0)=v_j(x).
  \end{equation}
  Thus, for $\varepsilon_{\ref{prop:ex-good-atlas}}$ sufficiently small, $\phi_j$ is a diffeomorphism onto a subset of $D_{32\rho^\ast,p_j}$ containing $D_{20\rho^\ast,p_j}$. Set $g_j=(\Psi_j\circ\phi_j)^*g_f$. By the definition of $v_j$, we have
  \begin{equation}
    \partial_y(f\circ\Psi_j\circ\phi_j)(x,0)=-\partial_\eta f\bigl(\Psi_j(x,\psi_j(x))\bigr).
  \end{equation}
  Since $\partial_x(f\circ\Psi_j\circ\phi_j)(x,0)$ is tangent to $\Gamma$, it follows that $(g_j)_{12}(x,0)=0$. Define $\mu_j$ and $a_j$ from $g_j$ as in \eqref{eq:3-pf:interior-isothermal-charts}. Then
  \begin{equation}
    \Vert \mu_j\Vert _{C^\alpha(Q_{24\rho^\ast}^+)}+\Vert a_j-1\Vert _{C^\alpha(Q_{24\rho^\ast}^+)}\leq C\varepsilon\quad\text{and}\quad \mu_j(x,0)\in\R.
  \end{equation}
  Consequently, reflection by $\mu_j(\bar z)=\overline{\mu_j(z)}$ defines a $C^\alpha$-coefficient on $Q_{24\rho^\ast}$. Extend it symmetrically to $\nu_j\in C_c^\alpha(\C)$, with support in $Q_{32\rho^\ast}$ and $\Vert \nu_j\Vert _{C^\alpha(\C)}\leq C\varepsilon$. From here, we can continue as in the proof of \Cref{lem:interior-isothermal-charts} to obtain a $C^{1+\alpha}$-diffeomorphism $\Phi_j$ satisfying
  \begin{equation}\label{eq:3-pf:boundary-isothermal-charts}
    \Phi_j(0)=0,\quad \partial_{\bar z}\Phi_j=\nu_j\partial_z\Phi_j\quad\text{and}\quad \Vert \Phi_j-\operatorname{id}\Vert _{C^{1+\alpha}(\C)}+\Vert \Phi_j^{-1}-\operatorname{id}\Vert _{C^{1+\alpha}(\C)}\leq C\varepsilon.
  \end{equation}
  The uniqueness used there implies that $\Phi_j(\bar z)=\overline{\Phi_j(z)}$, so $\Phi_j$ preserves the upper half-plane. Set $\alpha_j=\Psi_j\circ\phi_j\circ\Phi_j^{-1}$, initially on $Q_{20\rho^\ast}^+$, and let $\Omega_j=\alpha_j(Q_{16\rho^\ast}^+)$. Its inverse on $\Omega_j$ is $\varphi_j=\Phi_j\circ\phi_j^{-1}\circ P_{p_j}$. The inclusions and estimates follow as in the proof of \Cref{lem:interior-isothermal-charts}, taking
  \begin{equation}\label{eq:def-lambda_j}
    \omega_j=\phi_j\circ\Phi_j^{-1},\quad h_j=u_{p_j}\circ\omega_j\quad\text{and}\quad \lambda_j=\Big(\frac12\log a_j-\log|\partial_z\Phi_j|\Big)\circ\Phi_j^{-1}.
  \end{equation}
  Here, we reduce $\varepsilon_{\ref{prop:ex-good-atlas}}$ so that the sum of the $C^0$-errors of $\phi_j-\operatorname{id}$ and $\Phi_j-\operatorname{id}$ in \eqref{eq:2-pf:boundary-isothermal-charts} resp.\ \eqref{eq:3-pf:boundary-isothermal-charts} is at most $\rho^\ast/4$.
  Next we proceed as in the proof of \Cref{lem:interior-isothermal-charts} starting around \eqref{eq:21-pf:interior-isothermal-charts}. That is, let $j,k\in \mathbf I_{\mathrm{int}}\cup \mathbf I_{\mathrm{bdr}}$, and $w\in\varphi_k(\Omega_j\cap\Omega_k)$. In the case where $j\in \mathbf I_{\mathrm{bdr}}$, arguing analogously as for \eqref{eq:21-pf:interior-isothermal-charts}, there exists a universal $0<r^\ast\leq \rho^\ast/16$ such that, writing $\widehat\Omega_j=\alpha_j(Q_{18\rho^\ast}^+)$, we have $d_{g_f}(\Omega_j,\Sigma\setminus\widehat\Omega_j)\geq 2r^\ast$. 
  
  If $k\in \mathbf I_{\mathrm{int}}$, then $d_{g_f}(\Omega_k,\partial\Sigma)\geq 4r^\ast$ by \eqref{eq:21-pf:interior-isothermal-charts}. Therefore, we get that $T_{jk}=\alpha_j^{-1}\circ\alpha_k\colon B_{r^\ast}(w)\to Q_{18\rho^\ast}$ is well-defined, holomorphic, and we can proceed analogously to the proof of \Cref{lem:interior-isothermal-charts} after \eqref{eq:18-pf:interior-isothermal-charts}. This also covers the case where the roles of $j$ and $k$ are reversed. 

  Finally, suppose that both $j,k\in\mathbf I_{\mathrm{bdr}}$. If $\operatorname{Im}w>r^\ast/2$, we again get a holomorphic transition on $B_{r^\ast/2}(w)$ and proceed as in the interior case. If $0\leq\operatorname{Im}w\leq r^\ast/2$, set $w_0=\operatorname{Re}w$. After reducing $r^\ast$ universally, the margins of the enlarged charts and the same path-length argument as for \eqref{eq:21-pf:interior-isothermal-charts} show that
  \begin{equation}
    T_{jk}=\alpha_j^{-1}\circ\alpha_k\colon B_{2r^\ast}^+(w_0)\to Q_{18\rho^\ast}^+
  \end{equation}
  is well-defined and holomorphic, where $B_{2r^\ast}^+(w_0)=B_{2r^\ast}(w_0)\cap\{\operatorname{Im}z\geq0\}$. It maps the real interval into $\R$ and preserves the upper half-plane, so $T_{jk}'(w_0)>0$. By Schwarz reflection, it extends holomorphically to $B_{2r^\ast}(w_0)$, preserving the estimate for $\log|T_{jk}'|$ in \eqref{eq:23-pf:interior-isothermal-charts}. Choose the branch $q_{jk}=\log T_{jk}'$ with $q_{jk}(w_0)\in\R$. The harmonic estimates and the Cauchy--Riemann equations used above yield that $\Vert q_{jk}\Vert_{C^l(B_{r^\ast}(w_0))}\leq C_l\varepsilon$, and therefore $\Vert T_{jk}'-1\Vert_{C^l(B_{r^\ast}(w_0))}\leq C_l\varepsilon$. Since $w\in B_{r^\ast}(w_0)$, this proves \eqref{eq:4-lem:boundary-isothermal-charts} with $r_{\ref{lem:boundary-isothermal-charts}}=r^\ast/2$, and also yields the higher derivative estimates at the boundary.
\end{proof}

\subsection{Smooth approximation of clamped $C^{1+\alpha}$-immersions with estimates}

In this subsection, we prove the following approximation result.

\begin{proposition}\label{prop:smooth-clamped-immersion}
  There exist universal constants $C_{\ref{prop:smooth-clamped-immersion}}>0$ and $0<\varepsilon_{\ref{prop:smooth-clamped-immersion}}\leq\varepsilon_{\ref{prop:ex-good-atlas}}$ with the following property. Let $\Gamma\subset\R^3$ be a smooth compact not necessarily connected one-dimensional submanifold and assume that
  \begin{equation}\label{eq:1-prop:smooth-clamped-immersion}
    \mathfrak r_\Gamma>16 \quad\text{and}\quad 
    [\bm\tau_{\Gamma}]_{C^\alpha(\Gamma)}\leq\varepsilon\leq\varepsilon_{\ref{prop:smooth-clamped-immersion}}.
  \end{equation}
  Moreover, suppose that $\mathbf N\colon\Gamma\to\R^3$ is a smooth unit normal field along $\Gamma$ and that $f\in C^{1+\alpha}(\Sigma,\R^3)$ is an immersion with $f(\partial\Sigma)=\Gamma$ and $\partial_{\eta_f}f=\mathbf N\circ f$ on $\partial\Sigma$. Moreover, assume that
  \begin{equation}\label{eq:2-prop:smooth-clamped-immersion}
    [\mathbf n_f]_{\alpha,1}\leq\varepsilon\leq\varepsilon_{\ref{prop:smooth-clamped-immersion}},
  \end{equation}
  and fix the smooth structure induced by the isothermal atlas in \Cref{prop:ex-good-atlas}. There exist a $C^{1+\alpha}$-diffeomorphism $\Xi\colon\Sigma\to\Sigma$ and a smooth immersion $G\colon\Sigma\to\R^3$ such that
  \begin{enumerate}
    \item[(a)] $f\circ\Xi=G$ and $d(f\circ\Xi-G)=0$ on $\partial\Sigma$,
    \item[(b)] $\Vert G\circ\alpha_j-A_j(\cdot,0)\Vert _{C^{1+\alpha}(D_j)}\leq C_{\ref{prop:smooth-clamped-immersion}}\varepsilon$ for all $j\in \mathbf I_{\mathrm{int}}\cup \mathbf I_{\mathrm{bdr}}$ with the Euclidean isometries $A_j$ from \Cref{prop:ex-good-atlas}. Moreover, $\Xi=\operatorname{id}$ on $\Omega_j$ for all $j\in \mathbf I_{\mathrm{int}}$, while $\Xi(\alpha_j(Q_{12\rho^\ast}^+))\subset\alpha_j(Q_{16\rho^\ast}^+)$ and
    \begin{equation}\label{eq:3-prop:smooth-clamped-immersion}
      \Vert \varphi_j\circ\Xi\circ\alpha_j-\operatorname{id}\Vert _{C^{1+\alpha}(Q_{12\rho^\ast}^+)}\leq C_{\ref{prop:smooth-clamped-immersion}}\varepsilon \quad\text{for all $j\in \mathbf I_{\mathrm{bdr}}$}.
    \end{equation}
    In particular, we have $\Vert f\circ \Xi-G\Vert_{C^{1+\alpha}(\Sigma)}\leq C_{\ref{prop:smooth-clamped-immersion}}\varepsilon$ using \Cref{lem:unif-equivalent-norms} with respect to the atlas from \Cref{prop:ex-good-atlas} for $f$.
    \item[(c)] Let $k\geq2$ and additionally suppose that
    \begin{equation}\label{eq:4-prop:smooth-clamped-immersion}
      \Vert \bm\kappa_\Gamma\Vert _{C^{k-2}(\Gamma)}+\Vert  D_s\mathbf N\Vert _{C^{k-2}(\Gamma)}\leq M_k\varepsilon.
    \end{equation}
    Then, the same $G$ satisfies
    \begin{equation}\label{eq:5-prop:smooth-clamped-immersion}
      \Vert G\circ\alpha_j-A_j(\cdot,0)\Vert _{C^k(D_j)}\leq C_k(1+M_k)\varepsilon \quad\text{for all $j\in \mathbf I_{\mathrm{int}}\cup \mathbf I_{\mathrm{bdr}}$},
    \end{equation}
    and in particular,
    \begin{equation}
      \max_{0\leq l\leq k-2} \Vert \nabla_{g_G}^l \mathbf A_G\Vert_{C^0(\Sigma)}\leq C_{k,M_k}\varepsilon.
    \end{equation}    
  \end{enumerate}
  Here, $C_k$ depends at most on $k$ and $\alpha$. In \eqref{eq:4-prop:smooth-clamped-immersion}, $\bm\kappa_\Gamma=D_s\bm \tau_\Gamma$ denotes the curvature vector, and the $C^m(\Gamma)$-norms of vector fields use ambient arclength derivatives $D_s^\ell$, $0\leq\ell\leq m$. 
\end{proposition}

We split its proof in two parts. First, we construct a family of diffeomorphisms $(\Xi_u)_{0\leq u\leq 1}$ on $\partial\Sigma$ which, using a collar-neighborhood, is later extended to $\Xi$ with $\Xi\vert_{\partial\Sigma}=\Xi_1$. Moreover, we abbreviate $I_r=(-r,r)$ for $r>0$. 

\begin{lemma}\label{lem:boundary-reparametrisation}
  Work in the setting of \Cref{prop:ex-good-atlas,lem:interior-isothermal-charts,lem:boundary-isothermal-charts}, and fix the smooth structure induced by the isothermal atlas. After reducing $\varepsilon_{\ref{prop:ex-good-atlas}}$ universally, there exists a family of $C^{1+\alpha}$-diffeomorphisms $\Xi_u\colon\partial\Sigma\to\partial\Sigma$ for $0\leq u\leq 1$ such that, for all $j\in \mathbf I_{\mathrm{bdr}}$,
  \begin{align}
    \Xi_u(\alpha_j(I_{14\rho^\ast}\times\{0\}))\subset\alpha_j(I_{16\rho^\ast}\times\{0\}) &\quad\text{and}\\ 
    \Vert (x,u)\mapsto \varphi_j\circ\Xi_u\circ\alpha_j(x,0)-(x,0)\Vert _{C^{1+\alpha}(I_{14\rho^\ast}\times [0,1])}&\leq C\varepsilon.\label{eq:1-lem:boundary-reparametrisation}
  \end{align}
  Moreover, $\Xi_0=\mathrm{id}_{\partial\Sigma}$ and there exist smooth $\gamma_j\colon\R\to\Gamma$ with $\vert\gamma_j'\vert\equiv 1$, and smooth increasing $\xi_j\colon I_{14\rho^*}\to\R$ such that
  \begin{equation}\label{eq:2-lem:boundary-reparametrisation}
    f\circ\Xi_1\circ\alpha_j(x,0)=\gamma_j(\xi_j(x))\text{ for all $x\in I_{14\rho^*}$} \quad\text{and}\quad \Vert \xi_j-\operatorname{id}\Vert _{C^m(I_{14\rho^*})}\leq C_m\varepsilon \quad\text{for all $m\in\N_0$}.
  \end{equation}
  Here, $C$ is universal and $C_m$ depends at most on $m$.
\end{lemma}

We briefly comment on the underlying idea. In general, $f$ is only $C^{1+\alpha}$ on the boundary. But assuming that its image $\Gamma=f(\partial\Sigma)$ is smooth, each arclength parametrization $\gamma_j$ is necessarily smooth. In the simplest case, namely if there was a $C^{1+\alpha}$ boundary-diffeomorphism $\Xi_1\colon\partial\Sigma\to\partial\Sigma$ such that $\gamma_j=f\circ\Xi_1\circ\alpha_j(\cdot,0)$ is parametrized by \textbf{arclength} on $I_{15\rho^\ast}$ for all charts $j\in \mathbf I_{\mathrm{bdr}}$ (this corresponds to taking $\xi_j=\operatorname{id}$ in \eqref{eq:2-lem:boundary-reparametrisation}), then smoothness of each $\gamma_j$ yields that $f\circ\Xi_1\big\vert_{\partial\Sigma}$ is smooth in the isothermal atlas. A simple calculation shows that \eqref{eq:2-lem:boundary-reparametrisation} implies that the $\xi_j$ must necessarily satisfy $(\xi_j-\xi_k\circ\tau_{kj})'\equiv 0$ for the boundary transition functions $\tau_{kj}=\varphi_k^{(1)}\circ\alpha_j(\cdot,0)$, so taking $\xi_j=\operatorname{id}$ is only possible if the transitions $\tau_{kj}$ are translations on connected boundary overlaps. While this is not true in general, \eqref{eq:4-lem:boundary-isothermal-charts} yields that $\tau_{kj}'-1$ is uniformly small in $C^m$, which allows proving the result with the small corrections $\xi_j-\operatorname{id}$ in \eqref{eq:2-lem:boundary-reparametrisation}. 

\begin{proof}[Proof of \Cref{lem:boundary-reparametrisation}]
  Let $j,k\in\mathbf I_{\mathrm{bdr}}$ such that $\partial\Sigma\cap\Omega_j\cap\Omega_k\neq\emptyset$. Then let $\tau_{kj}=\varphi_k^{(1)}\circ\alpha_j(\cdot,0)$ be the restriction of $T_{kj}$ from \eqref{eq:4-lem:boundary-isothermal-charts} to the boundary charts. In particular, we get that 
  \begin{equation}\label{eq:12-pf:boundary-reparametrisation}
    \Vert \tau_{kj}'-1\Vert_{C^{m-1}(\varphi_j(\partial\Sigma\cap\Omega_j\cap\Omega_k))}\leq C_m\varepsilon \quad\text{for every $m\geq1$}.
  \end{equation}
  Next, we construct the \emph{corrections} $\xi_j$ in \eqref{eq:2-lem:boundary-reparametrisation}. To this end, we first fix $C^{1+\alpha}$-arclength reparametrizations $s_j$ in each coordinate chart such that $\gamma_j(s)\vcentcolon=f(\alpha_j(s_j^{-1}(s),0))$ is parametrized by arc-length, and, by adding small correction terms accounting for the error in $\tau_{kj}'-1$, construct $\xi_j$ such that
  \begin{equation}\label{eq:1.5-pf:boundary-reparametrization}
    \xi_j\text{ is smooth, and} \qquad \Vert \xi_j-\operatorname{id}\Vert _{C^m(I_{15\rho^\ast})}\leq C_m\varepsilon \quad\text{for every $m\in\N_0$}.
  \end{equation}
  Setting $s_j(x)=\int_0^x e^{\lambda_j(h,0)}\dd h$ the (oriented) arclength of the boundary arc between $p_j=\alpha_j(0,0)$ and $\alpha_j(x,0)$, by \Cref{prop:ex-good-atlas}, we have
  \begin{equation}\label{eq:2-pf:boundary-reparametrisation}
    \Vert e^{\lambda_j(\cdot,0)}-1\Vert _{C^\alpha(I_{16\rho^\ast})}+\Vert s_j-\operatorname{id}\Vert _{C^{1+\alpha}(I_{16\rho^\ast})}\leq C\varepsilon.
  \end{equation}
  Thus, defining $d_k(\alpha_k(x,0))=x-s_k(x)$ for $k\in \mathbf I_{\mathrm{bdr}}$ and $d=\sum_{k\in \mathbf I_{\mathrm{bdr}}}\pi_k^2 d_k$ with the localization system from \Cref{rem:local-sys}, \eqref{eq:12-pf:boundary-reparametrisation}, \eqref{eq:2-pf:boundary-reparametrisation} and \eqref{eq:unif-bounds-pij} yield
  \begin{equation}\label{eq:5-pf:boundary-reparametrisation}
    \Vert d\circ\alpha_j(\cdot,0)\Vert _{C^{1+\alpha}(I_{16\rho^\ast})}\leq C\varepsilon.
  \end{equation}
  We define
  \begin{equation}\label{eq:6-pf:boundary-reparametrisation}
    \xi_j(x)\vcentcolon=s_j(x)+d(\alpha_j(x,0)).
  \end{equation}
  As in \eqref{eq:23-pf:interior-isothermal-charts}, we have $\tau_{kj}'=e^{\lambda_j-\lambda_k\circ\tau_{kj}}$ and thus,
  \begin{equation}\label{eq:7-pf:boundary-reparametrisation}
    (s_j-s_k\circ\tau_{kj})'\equiv 0\quad\text{so that $s_j-s_k\circ\tau_{kj}$ is smooth}.
  \end{equation}
  Using \eqref{eq:6-pf:boundary-reparametrisation}, we compute
  \begin{equation}
    \xi_j(x)-x=\sum_k\pi_k^2\circ\alpha_j(\cdot,0)\big([s_j-s_k\circ\tau_{kj}]+[\tau_{kj}-\operatorname{id}]\big)\Big\vert_x = \sum_k\pi_k^2\big(d_k-d_j\big)\Big\vert_{\alpha_j(x,0)}.
  \end{equation}
  The first equality and \eqref{eq:7-pf:boundary-reparametrisation} show that $\xi_j$ is smooth. For $\xi_j-\operatorname{id}$, the second equality and \eqref{eq:2-pf:boundary-reparametrisation} together yield that the $C^0(I_{16\rho^\ast})$-norm is at most $C\varepsilon$, and the higher-order norms are estimated by the first equality and \eqref{eq:7-pf:boundary-reparametrisation} together with \eqref{eq:12-pf:boundary-reparametrisation}. Altogether, this yields \eqref{eq:1.5-pf:boundary-reparametrization}.

  Using \eqref{eq:2-pf:boundary-reparametrisation} and \eqref{eq:5-pf:boundary-reparametrisation}, after universally shrinking $\varepsilon_{\ref{prop:ex-good-atlas}}$, we have that $s_j\colon I_{16\rho^\ast}\to\R$ is a $C^{1+\alpha}$-diffeomorphism onto its image, $s_j(x)+ud(\alpha_j(x,0))\in I_{15\rho^\ast}$ for all $x\in I_{14\rho^\ast}$ and $0\leq u\leq 1$, and $I_{15\rho^\ast}\subset s_j(I_{16\rho^\ast})$.  Moreover, \eqref{eq:2-pf:boundary-reparametrisation} yields that 
  \begin{equation}\label{eq:9.5-pf:boundary-reparametrisation}
    \Vert (x,u)\mapsto (s_j^{-1}\big(s_j(x)+u\,d(\alpha_j(x,0))\big)-x)\Vert_{C^{1+\alpha}(I_{14\rho^\ast}\times[0,1])}\leq C\varepsilon.
  \end{equation}
  Combined with \eqref{eq:7-pf:boundary-reparametrisation}, setting
  \begin{equation}\label{eq:9-pf:boundary-reparametrisation}
    \Xi_u(\alpha_j(x,0))=\alpha_j\Big(s_j^{-1}\big(s_j(x)+u\,d(\alpha_j(x,0))\big),0\Big)\quad\text{for all $j\in \mathbf I_{\mathrm{bdr}}$ and $x\in I_{14\rho^\ast}$},
  \end{equation}
  we get a $C^{1+\alpha}$-diffeomorphism $\Xi_u\colon\partial\Sigma\to\partial\Sigma$ for any $0\leq u\leq 1$, after possibly reducing $\varepsilon_{\ref{prop:ex-good-atlas}}$, and \eqref{eq:1-lem:boundary-reparametrisation} follows from \eqref{eq:9.5-pf:boundary-reparametrisation}. Moreover, $\Xi_0=\operatorname{id}$ and $\Xi_1$ satisfies that $f\circ\Xi_1\circ\alpha_j(x,0)=\gamma_j(\xi_j(x))$ with the arclength parametrizations $\gamma_j$ satisfying $\gamma_j\circ s_j=f\circ\alpha_j(\cdot,0)$. 
\end{proof}

\begin{proof}[Proof of \Cref{prop:smooth-clamped-immersion}]
  After reducing $\varepsilon_{\ref{prop:smooth-clamped-immersion}}$ universally so that $2\varepsilon_{\ref{prop:smooth-clamped-immersion}}\leq\varepsilon_{\ref{prop:ex-good-atlas}}$, apply \Cref{prop:ex-good-atlas} with $2\varepsilon$ in place of $\varepsilon$, and apply \Cref{lem:boundary-reparametrisation} to this atlas.
  Using the localization system from \Cref{rem:local-sys}, the idea is to define
  \begin{equation}\label{eq:14-pf:smooth-clamped-immersion}
    G\vcentcolon=\sum_{j\in \mathbf I_{\mathrm{int}}\cup \mathbf I_{\mathrm{bdr}}}\pi_j^2 \cdot G_j\circ\varphi_j
  \end{equation}
  where the $G_j$ are local approximations defined on neighborhoods of $\varphi_j(\operatorname{supp}\pi_j)$. For $j\in \mathbf I_{\mathrm{int}}$, we will take $G_j$ to parametrize a rectangle in the shifted tangent space $f(p_j)+df_{p_j}(T_{p_j}\Sigma)$. For $j\in\mathbf I_{\mathrm{bdr}}$, we use a construction similar to the proof of \Cref{lem:boundary-isothermal-charts}, namely extending the boundary parametrizations $\gamma_j\circ\xi_j$ from \Cref{lem:boundary-reparametrisation} to a surface using the outer conormal $\mathbf N$.

  More precisely, with $A_j$ from \Cref{prop:ex-good-atlas}, let $P_j(x,y)\vcentcolon=A_j(x,y,0)=c_j+a_jx+b_jy$ where $c_j=f(p_j)$, and $a_j,b_j$ are an orthonormal basis of the tangent space $df_{p_j}(T_{p_j}\Sigma)\subset\R^3$. For $j\in \mathbf I_{\mathrm{int}}$, set $G_j=P_j$ on $D_j$. For the definition of $G_j$ for $j\in\mathbf I_{\mathrm{bdr}}$, we use the operator $\mathcal E$ from \Cref{lem:whitney-extension}.
  For $j\in \mathbf I_{\mathrm{bdr}}$, using $\Xi_1$, $\gamma_j$ and $\xi_j$ from \Cref{lem:boundary-reparametrisation}, abbreviate
  \begin{equation}\label{eq:4-pf:smooth-clamped-immersion}
    \beta_j(x)\vcentcolon=f\circ\Xi_1\circ\alpha_j(x,0)=\gamma_j\circ\xi_j(x) \quad\text{and}\quad \nu_j(x)\vcentcolon=\mathbf N(\beta_j(x)).
  \end{equation}
  In particular, both are smooth, $|\nu_j|=1$ and $\nu_j\cdot\beta_j'=0$. Moreover, using \eqref{eq:almost-graph-charts} and \eqref{eq:1-lem:boundary-reparametrisation},
  \begin{equation}\label{eq:5-pf:smooth-clamped-immersion}
    \Vert \beta_j-P_j(\cdot,0)\Vert _{C^{1+\alpha}(I_{14\rho^\ast})}+\Vert \nu_j+b_j\Vert _{C^\alpha(I_{14\rho^\ast})}\leq C\varepsilon.
  \end{equation}
  When additionally assuming \eqref{eq:4-prop:smooth-clamped-immersion} for a fixed $k\geq 2$, using $(\gamma_j'')^{(l)}=\partial_s^l\bm\kappa_\Gamma\circ\gamma_j$ for $l\geq 0$ and \eqref{eq:2-lem:boundary-reparametrisation} and that $P_j(\cdot,0)''\equiv 0$, \eqref{eq:4-pf:smooth-clamped-immersion} yields 
  \begin{equation}\label{eq:6-pf:smooth-clamped-immersion}
    \Vert \beta_j-P_j(\cdot,0)\Vert _{C^{k}(I_{14\rho^\ast})}+\Vert \nu_j+b_j\Vert _{C^{k-1}(I_{14\rho^\ast})}\leq C_{k}(1+M_k)\varepsilon.
  \end{equation}
  Then, multiply $\nu_j+b_j$ by a fixed smooth cut-off equal to one on $I_{13\rho^\ast}$ and supported in $I_{14\rho^\ast}$, and extend the product by zero to $u_j\colon\R\to\R^3$. Define $G_j(x,y)\vcentcolon=P_j(x,y)+\beta_j(x)-P_j(x,0)-\mathcal Eu_j(x,y)$ on $Q_{13\rho^\ast}^+$. By \eqref{eq:whitney-extension-traces}, we have
  \begin{equation}\label{eq:12-pf:smooth-clamped-immersion}
    G_j(x,0)=\beta_j(x),\quad \partial_xG_j(x,0)=\beta_j'(x)\quad\text{and}\quad \partial_yG_j(x,0)=-\nu_j(x),
  \end{equation}
  so that \eqref{eq:whitney-extension-estimates} and \eqref{eq:5-pf:smooth-clamped-immersion} yield
  \begin{equation}\label{eq:13-pf:smooth-clamped-immersion}
    \Vert G_j-P_j\Vert _{C^{1+\alpha}(Q_{13\rho^\ast}^+)}
    \leq C\varepsilon.
  \end{equation}
  And additionally assuming \eqref{eq:4-prop:smooth-clamped-immersion} for a fixed $k\geq 2$, \eqref{eq:whitney-extension-estimates} and \eqref{eq:6-pf:smooth-clamped-immersion} imply that
  \begin{equation}\label{eq:13.1-pf:smooth-clamped-immersion}
    \Vert G_j-P_j\Vert _{C^{k}(Q_{13\rho^\ast}^+)}\leq C_k(1+M_k)\varepsilon.
  \end{equation}
  Then $G$ is defined by \eqref{eq:14-pf:smooth-clamped-immersion}, and we have
  \begin{align}
    G\circ\alpha_j-P_j = \sum_{k\in\mathbf I_{\mathrm{int}}\cup\mathbf I_{\mathrm{bdr}}} \pi_k^2\circ\alpha_j \big((G_k-P_k)\circ T_{kj} + P_k\circ T_{kj}-P_j \big). 
  \end{align}
  By \eqref{eq:almost-graph-charts}, we obtain that $\Vert P_k\circ T_{kj}-P_j\Vert_{C^{1}(\varphi_j(\Omega_j\cap\Omega_k))}\leq C\varepsilon$. Thus, additionally using that $P_j$ and $P_k$ are affine, and \eqref{eq:4-lem:interior-isothermal-charts} and \eqref{eq:4-lem:boundary-isothermal-charts}, we get that
  \begin{equation}\label{eq:13.2-pf:smooth-clamped-immersion}
    \Vert P_k\circ T_{kj}-P_j\Vert_{C^{l}(\varphi_j(\Omega_j\cap\Omega_k))}\leq C_l\varepsilon \quad\text{for all $l\geq 0$ and constants $C_l>0$ depending at most on $l$}.
  \end{equation}
  Combined with \eqref{eq:13-pf:smooth-clamped-immersion}, \eqref{eq:unif-bounds-pij} and \eqref{eq:2-prop:ex-good-atlas}, we get
  \begin{equation}\label{eq:13.3-pf:smooth-clamped-immersion}
    \Vert G\circ\alpha_j-P_j\Vert_{C^{1+\alpha}(D_j)}\leq C\varepsilon\quad\text{for all $j\in \mathbf I_{\mathrm{int}}\cup \mathbf I_{\mathrm{bdr}}$},
  \end{equation}
  and when again additionally assuming \eqref{eq:4-prop:smooth-clamped-immersion} for a fixed $k\geq 2$, \eqref{eq:13.1-pf:smooth-clamped-immersion} yields that
  \begin{equation}\label{eq:13.4-pf:smooth-clamped-immersion}
    \Vert G\circ\alpha_j-P_j\Vert_{C^{k}(D_j)}\leq C_k(1+M_k)\varepsilon\quad\text{for all $j\in \mathbf I_{\mathrm{int}}\cup \mathbf I_{\mathrm{bdr}}$}.
  \end{equation}
  After universally reducing $\varepsilon_{\ref{prop:smooth-clamped-immersion}}$ if necessary, \eqref{eq:13.3-pf:smooth-clamped-immersion} shows that $G$ is an immersion and proves the estimate for $G$ in (b). The coordinate formulae for the metric, its inverse and the second fundamental form, together with \eqref{eq:13.3-pf:smooth-clamped-immersion} and \eqref{eq:13.4-pf:smooth-clamped-immersion}, yield the curvature estimate in (c).

  We next compute the zeroth and first order boundary conditions of $G$. As in \eqref{eq:22-pf:interior-isothermal-charts}, there are $a_{jk}>0$ with
  \begin{equation}\label{eq:16-pf:smooth-clamped-immersion}
    DT_{jk}(x,0)e_2=a_{jk}(x)e_2.
  \end{equation}
  At the boundary, all local maps in \eqref{eq:14-pf:smooth-clamped-immersion} have the same value, since $\beta_j\circ\tau_{jk}=\beta_k$. Thus, the terms involving derivatives of the weights cancel by $d(\sum_j\pi_j^2)=0$. Using \eqref{eq:12-pf:smooth-clamped-immersion} and \eqref{eq:16-pf:smooth-clamped-immersion}, we get
  \begin{equation}\label{eq:17-pf:smooth-clamped-immersion}
    \partial_1(G\circ\alpha_k)(x,0)=\beta_k'(x) \quad\text{and}\quad \partial_2(G\circ\alpha_k)(x,0)=-\sum_j\pi_j^2(\alpha_k(x,0))a_{jk}(x) \cdot \nu_k(x).
  \end{equation}
  Altogether, we get
  \begin{equation}\label{eq:18-pf:smooth-clamped-immersion}
    G\circ\alpha_k(x,0)=\beta_k(x) \quad\text{and}\quad \partial_{\eta_G}G\circ\alpha_k(x,0)=\nu_k(x)=\mathbf N(\beta_k(x)).
  \end{equation}
  It remains to construct $\Xi$ and prove (a) and the estimate for $\Xi$ in (b). To this end, define the smooth vector field
  \begin{equation}\label{eq:collar-field-smooth-clamped-immersion}
    V(p)\vcentcolon=\sum_{j\in\mathbf I_{\mathrm{bdr}}}\pi_j^2(p)(d\alpha_j)_{\varphi_j(p)}.e_2.
  \end{equation}
  By the construction in \Cref{lem:interior-isothermal-charts}, $d_{g_f}(\Omega_j,\partial\Sigma)\geq \frac{\rho^\ast}{4}$. Therefore, for any $p\in\Sigma$ such that $d_{g_f}(p,\partial\Sigma)$ is at most a small universal constant, we have $\sum_{j\in\mathbf I_{\mathrm{bdr}}}\pi_j^2(p)=1$. Using \eqref{eq:4-lem:boundary-isothermal-charts} and \eqref{eq:unif-bounds-pij}, we get, for a universal $r_0>0$,
  \begin{equation}\label{eq:collar-field-estimate-smooth-clamped-immersion}
    \Vert d\varphi_j.(V\circ\alpha_j)-e_2\Vert_{C^l(I_{15\rho^\ast}\times[0,r_0])}\leq C_l\varepsilon\quad\text{for all $j\in\mathbf I_{\mathrm{bdr}}$ and $l\in\N_0$}.
  \end{equation}
  Thus, the flow of $V$ starting at $\partial\Sigma$ defines a smooth collar $\mathbf C\colon\partial\Sigma\times[0,2a^\ast)\to\Sigma$ for a sufficiently small universal $a^\ast>0$. Differentiating the flow equation and using \eqref{eq:collar-field-estimate-smooth-clamped-immersion} yields
  \begin{gather}
    \mathbf C(p,0)=p,\quad \partial_y\mathbf C(p,0)=V(p)\quad\text{and}\quad \mathbf C(\alpha_j(I_{14\rho^\ast}\times\{0\})\times[0,2a^\ast))\subset\Omega_j,\label{eq:19-pf:smooth-clamped-immersion}\\
    \Vert\varphi_j\circ\mathbf C(\alpha_j(x,0),y)-(x,y)\Vert_{C^l(I_{14\rho^\ast}\times[0,2a^\ast))}\leq C_l\varepsilon\quad\text{for all $l\in\N_0$}.\label{eq:20-pf:smooth-clamped-immersion}
  \end{gather}
  The inverse collar maps satisfy the corresponding estimates. By choosing $a^\ast$ smaller if necessary, we get $\mathbf C(\partial\Sigma\times[0,a^\ast))\cap\Omega_j = \emptyset$ for all $j\in\mathbf I_{\mathrm{int}}$.

  Moreover, \eqref{eq:16-pf:smooth-clamped-immersion} shows that $V(\alpha_j(x,0))$ is a positive multiple of $d\alpha_j.e_2$. Using \eqref{eq:3-prop:ex-good-atlas} and \eqref{eq:17-pf:smooth-clamped-immersion}, $V(\alpha_j(x,0))$ is thus normal to $T_{\alpha_j(x,0)}\partial\Sigma$ with respect to both $g_f(\alpha_j(x,0))$ and $g_G(\alpha_j(x,0))$. Set
  \begin{equation}\label{eq:boundary-speeds-smooth-clamped-immersion}
    m_f(p)\vcentcolon=|df_p.V(p)|\quad\text{and}\quad m_G(p)\vcentcolon=|dG_p.V(p)|\quad\text{for $p\in\partial\Sigma$}.
  \end{equation}
  These functions are positive and satisfy
  \begin{gather}
    df_p.V(p)=-m_f(p)\mathbf N(f(p))\quad\text{and}\quad dG_p.V(p)=-m_G(p)\mathbf N(G(p)),\label{eq:normal-collar-derivatives-smooth-clamped-immersion}\\
    \Vert m_f-1\Vert_{C^\alpha(\partial\Sigma)}+\Vert m_G-1\Vert_{C^\alpha(\partial\Sigma)}\leq C\varepsilon.\label{eq:boundary-speeds-estimate-smooth-clamped-immersion}
  \end{gather}
  The estimate follows from \eqref{eq:almost-graph-charts}, \eqref{eq:13.3-pf:smooth-clamped-immersion} and \eqref{eq:collar-field-estimate-smooth-clamped-immersion}. Define
  \begin{equation}\label{eq:normal-speed-ratio-smooth-clamped-immersion}
    \vartheta(p)\vcentcolon=\frac{m_G(p)}{m_f(\Xi_1(p))}\quad\text{and}\quad h(p,y)\vcentcolon=y+\mathcal E_{\partial\Sigma}(\vartheta-1)(p,y).
  \end{equation}
  By \eqref{eq:1-lem:boundary-reparametrisation}, \eqref{eq:boundary-speeds-estimate-smooth-clamped-immersion} and \Cref{lem:whitney-extension}, we have
  \begin{equation}\label{eq:height-correction-smooth-clamped-immersion}
    h(p,0)=0,\quad \partial_yh(p,0)=\vartheta(p)\quad\text{and}\quad \Vert\vartheta-1\Vert_{C^\alpha(\partial\Sigma)}+\Vert h-y\Vert_{C^{1+\alpha}(\partial\Sigma\times[0,a^\ast])}\leq C\varepsilon.
  \end{equation}
  Fix a smooth $\chi\colon[0,a^\ast]\to[0,1]$ equal to one on $[0,a^\ast/4]$ and zero on $[a^\ast/2,a^\ast]$. On the collar, define
  \begin{equation}\label{eq:def-Xi-smooth-clamped-immersion}
    \Xi(\mathbf C(p,y))\vcentcolon=\mathbf C\big(\Xi_{\chi(y)}(p),\chi(y)h(p,y)+(1-\chi(y))y\big)\quad\text{for $0\leq y<a^\ast$},
  \end{equation}
  and set $\Xi=\operatorname{id}$ outside the collar. For sufficiently small $\varepsilon_{\ref{prop:smooth-clamped-immersion}}$, the second argument on the right hand side is strictly increasing in $y$, is zero at $y=0$ and equals $y$ for $y\geq a^\ast/2$. Thus, it takes values in $[0,a^\ast)$, and the definition agrees with the identity near the outer edge of the collar. Equations~\eqref{eq:1-lem:boundary-reparametrisation}, \eqref{eq:20-pf:smooth-clamped-immersion} and \eqref{eq:height-correction-smooth-clamped-immersion} show that $\Xi$ is $C^{1+\alpha}$ and satisfies \eqref{eq:3-prop:smooth-clamped-immersion}. Moreover, $\Xi$ is a $C^{1+\alpha}$-diffeomorphism.

  By construction, $\Xi\vert_{\partial\Sigma}=\Xi_1$ and hence $f\circ\Xi=G$ on $\partial\Sigma$. Thus, $d(f\circ \Xi-G)\big\vert_{T_p\partial\Sigma}=0$ on $\partial\Sigma$. Since $\chi$ is constant near zero, differentiating \eqref{eq:def-Xi-smooth-clamped-immersion} at $y=0$ and using \eqref{eq:normal-collar-derivatives-smooth-clamped-immersion} and \eqref{eq:normal-speed-ratio-smooth-clamped-immersion} yields that
  \begin{equation}\label{eq:matching-first-derivative-smooth-clamped-immersion}
    d(f\circ\Xi)_p.V(p)=\vartheta(p)df_{\Xi_1(p)}.V(\Xi_1(p))=-m_G(p)\mathbf N(G(p))=dG_p.V(p).
  \end{equation}
  Since $V(p)\notin T_p\partial\Sigma$ for $p\in\partial\Sigma$, this proves $d(f\circ\Xi-G)=0$ on $\partial\Sigma$, and thus (a).
\end{proof}

\section{Life-span and parabolic smoothing estimates for the Willmore flow}\label{sec:lifespan-estimate}

For an introduction of the parabolic function spaces, fundamental estimates, and the linear parabolic theory in time-weighted parabolic $L^p$-spaces, see \Cref{app:max-reg-unif-surfaces}. Let $\alpha\in(0,1)$ be fixed, and choose
\begin{align}
  \quad p > \frac{10}{3\alpha}\quad\text{and}\quad \alpha'&\in \Big(0,\alpha-\frac{2}{p}\Big).\label{eq:choice-of-parameters}
\end{align}
The reason for this choice becomes apparent later in \Cref{rem:choice-of-parameters:resulting-estimates}. The time-weight parameter $\mu\in (0,1)$ is also fixed in \Cref{rem:choice-of-parameters:resulting-estimates} below. In particular, it will be chosen such that
\begin{equation}\label{eq:1-rem:choice-of-parameters:resulting-estimates}
  C^{1+\alpha}(\Sigma)\hookrightarrow W^{4(\mu-\frac1p),p}(\Sigma)\hookrightarrow C^{1+\alpha'}(\Sigma).
\end{equation}
In this section, we prove the following main result. We use the notation of \Cref{app:graph-rep}.

\begin{theorem}\label{thm:lifespan-and-smoothing}
  Fix $K\geq 1$. There exist constants $\tau_{\ref{thm:lifespan-and-smoothing}}\in(0,1)$, and $0<\varepsilon_{\ref{thm:lifespan-and-smoothing}}\leq\varepsilon_{\ref{prop:smooth-clamped-immersion}}$, only depending on $K$ and the fixed parameters with the following property. Suppose that $\Gamma\subset\R^3$ is a smooth, not necessarily connected one-dimensional submanifold, and $\mathbf N$ a smooth unit normal field along $\Gamma$, and let
  \begin{equation}
    \mathfrak r_\Gamma>16\quad\text{and}\quad [\bm\tau_\Gamma]_{C^\alpha(\Gamma)}\leq\varepsilon\leq\varepsilon_{\ref{thm:lifespan-and-smoothing}},
  \end{equation}
  using the notation in \Cref{def:rGamma}. Further, assume that, for some $k\geq 5$,
  \begin{equation}\label{eq:boundary-data-small}
    \Vert\bm\kappa_\Gamma\Vert_{C^{k-2}(\Gamma)}+\Vert D_s\mathbf N\Vert_{C^{k-2}(\Gamma)}\leq K\varepsilon.
  \end{equation}
  Here, $\bm\kappa_\Gamma=D_s\bm\tau_\Gamma$, and all derivatives of the boundary data are ambient arclength derivatives. Let $f_0\in C^{1+\alpha}(\Sigma,\R^3)$ be an immersion with $f_0(\partial\Sigma)=\Gamma$ and $\partial_{\eta_{f_0}}f_0=\mathbf N\circ f_0$ on $\partial\Sigma$ such that
  \begin{equation}
    [\mathbf n_{f_0}]_{\alpha,1}\leq\varepsilon\leq\varepsilon_{\ref{thm:lifespan-and-smoothing}}.
  \end{equation}
  Choose the smooth immersion $G_{f_0}$ and diffeomorphism $\Xi_{f_0}$ from \Cref{prop:smooth-clamped-immersion}. Then there exists $T>0$ and a family of immersions
  \begin{equation}\label{eq:1-thm:lifespan-and-smoothing}
    f\in C^0([0,T],C^{1+\alpha'}(\Sigma,\R^3))\cap W^{1,p}_\mu(0,T;L^p(\Sigma))\cap L^p_\mu(0,T;W^{4,p}(\Sigma))
  \end{equation}
  solving
  \begin{equation}
    \begin{cases}
      \langle\partial_t f,\mathbf n_f\rangle=-\frac12\big(\Delta_{g_f}\mathbf H_f+|\mathbf A_f^0|^2\mathbf H_f\big)&\text{a.e.\ on $[0,T]\times\Sigma$}\\
      f(0,\cdot)=f_0\circ\Xi_{f_0}&\text{on $\Sigma$}\\
      f=f_0\circ\Xi_{f_0}\text{ and }\partial_{\eta_f}f=\mathbf N\circ f&\text{on $[0,T]\times\partial\Sigma$}.
    \end{cases}
  \end{equation}
  Moreover, $T\geq\tau_{\ref{thm:lifespan-and-smoothing}}$ and we have that
  \begin{align}
    &\sup_{0\leq t\leq\tau_{\ref{thm:lifespan-and-smoothing}}}\Vert f(t)-f_0\circ\Xi_{f_0}\Vert_{\bm C^{1+\alpha'}(\Sigma)}\\
    &\quad+\Vert f-G_{f_0}\Vert_{\bm W^{1,p}_\mu(0,\tau_{\ref{thm:lifespan-and-smoothing}};\bm L^p(\Sigma))\cap\bm L^p_\mu(0,\tau_{\ref{thm:lifespan-and-smoothing}};\bm W^{4,p}(\Sigma))}\leq C \varepsilon\label{eq:2-thm:lifespan-and-smoothing}
  \end{align}
  for a universal constant $C>0$. Regularity of $f$ in \eqref{eq:1-thm:lifespan-and-smoothing} and all norms in \eqref{eq:2-thm:lifespan-and-smoothing} are understood with respect to the isothermal atlas from \Cref{prop:ex-good-atlas} with respect to $f_0$. 
  Moreover, $f(t)\in C^{k-1+\alpha'}(\Sigma)$ for all $t\in(0,\tau_{\ref{thm:lifespan-and-smoothing}}]$ and there exists a constant $C_k>0$, only depending on $k$, $K$ and the fixed parameters such that
  \begin{equation}\label{eq:curvature-smoothing}
    \max_{\frac12\tau_{\ref{thm:lifespan-and-smoothing}}\leq t\leq\tau_{\ref{thm:lifespan-and-smoothing}}}\sum_{\ell=0}^{k-3}\big\Vert\nabla_{g_{f(t)}}^\ell\mathbf A_{f(t)}\big\Vert_{C^0(\Sigma)}\leq C_k\varepsilon.
  \end{equation}
\end{theorem}

Throughout this section, constants are called universal if they depend at most on $K$ and the fixed parameters. We sometimes only write $G=G_{f_0}$ and $\Xi=\Xi_{f_0}$.

\begin{remark}[Choice of parameters]\label{rem:choice-of-parameters:resulting-estimates}
  We motivate \eqref{eq:choice-of-parameters} by the resulting Gagliardo--Nirenberg interpolation estimates. First, \eqref{eq:choice-of-parameters} allows choosing
  \begin{align}
    \max\Big\{\frac{1+\alpha'}{4}+\frac{3}{2p},\frac14+\frac{11}{6p}\Big\}&<\mu<\frac{1+\alpha}{4}+\frac1p \quad\text{and}\quad
  \frac34+\frac{1}{2p} <\beta < \mu-\frac4{3p}+\frac12.\label{eq:choice-of-parameters-2}
  \end{align}
  Then, we actually get \eqref{eq:1-rem:choice-of-parameters:resulting-estimates}
  where the first embedding is clear with \Cref{lem:unif-equivalent-norms}, and the second follows from \Cref{prop:embeddings-and-interpolation-results}\eqref{it:Morrey-prop:embeddings-and-interpolation-results}. Using \eqref{eq:choice-of-parameters} and \Cref{prop:embeddings-and-interpolation-results}\eqref{it:Sobolev-prop:embeddings-and-interpolation-results} and \eqref{it:gagliardonirenberg-prop:embeddings-and-interpolation-results}, we have for a universal constant $C>0$ only depending on the specific choice of parameters that
  \begin{align}
    \Vert u \Vert_{\bm W^{2,3p}(\Sigma)} &\leq C \Vert u\Vert_{\bm W^{4\beta,p}(\Sigma)}^{\theta_{2,4\beta}}\Vert u\Vert_{\bm W^{4(\mu-\frac1p),p}(\Sigma)}^{1-\theta_{2,4\beta}} \quad\text{with } \theta_{2,4\beta} = \frac{\frac12+\frac{4}{3p}-\mu}{\beta-\mu+\frac1p}\in\Big(0,\frac12\Big), \\
    \Vert u \Vert_{\bm W^{3,\frac{3p}{2}}(\Sigma)} &\leq C \Vert u\Vert_{\bm W^{4\beta,p}(\Sigma)},\quad \Vert u\Vert_{\bm C^3(\Sigma)}\leq C\Vert u\Vert_{\bm W^{4\beta,p}(\Sigma)}, \\
    \Vert u \Vert_{\bm W^{4\beta,p}(\Sigma)} &\leq C \Vert u\Vert_{\bm W^{4,p}(\Sigma)}^{\theta_{4\beta,4}}\Vert u\Vert_{\bm W^{4(\mu-\frac1p),p}(\Sigma)}^{1-\theta_{4\beta,4}} \quad\text{with } \theta_{4\beta,4} = \frac{\beta-\mu+\frac1p}{1-\mu+\frac1p}\in(0,1)    
  \end{align}
  where the exponents satisfy the relations
  \begin{align}
    3\theta_{2,4\beta} &(\beta-\mu+\frac1p) + (\mu-\frac1p) < (1+\theta_{2,4\beta})(\beta-\mu+\frac1p) + (\mu-\frac1p) <1,\\
    &\text{or equivalently} \quad 3 \theta_{2,4\beta}\theta_{4\beta,4} < (1+\theta_{2,4\beta})\theta_{4\beta,4} < 1.
  \end{align}
  These estimates later allow controlling the nonlinearities appearing in the DeTurck--Willmore flow PDE in terms of the norm of the parabolic function space $\mathbb E_{1,\mu}(J_T)$.
\end{remark}

In order to prove \Cref{thm:lifespan-and-smoothing}, we study the so-called \emph{DeTurck--Willmore flow} starting at $f_0$, with \emph{Dirichlet boundary conditions}. That is, we study families of immersions
\begin{equation}
  f\in C^0([0,T],C^{1+\alpha'}_{\mathrm{imm}}(\Sigma,\R^3) \quad\text{and}\quad f-G_{f_0}\in\mathbb E_{1,\mu}(J_T)
\end{equation}
solving the system of PDEs given by
\begin{equation}\label{eq:DeTurck-Willmore-flow}
  \begin{cases}
    \partial_t f = -\frac12\big( \Delta_{g_f}\mathbf H_f + \vert \mathbf A^0_f\vert^2\mathbf H_f \big)\bm n_f- (\Delta_{g_f}\mathbf W_{f,G})^k\partial_k f &\text{a.e.\ on $[0,T]\times\Sigma$}\\
    f(0,\cdot)=f_0\circ\Xi_{f_0}&\text{on $\Sigma$}\\
    f - G_{f_0} = 0\text{ and } d(f-G_{f_0})=0 &\text{on $[0,T]\times\partial\Sigma$}.   
  \end{cases}
\end{equation}
Here, the 1-contravariant tensor $\mathbf W_{f,{G_{f_0}}}$ is defined by
\begin{equation}
  \mathbf W_{f,G_{f_0}}^k = \frac12 g_f^{\mu\nu}\big((\bm \Gamma_{f})^{k}_{\mu\nu}-(\bm \Gamma_{{G_{f_0}}})^{k}_{\mu\nu}\big)
\end{equation}
where $\bm\Gamma_{f}$ and $\bm\Gamma_{{G_{f_0}}}$ are the Christoffel symbols of $g_f$ and $g_{G_{f_0}}$, respectively. The idea is to make the ansatz $f=G_{f_0}+u$, and also view \eqref{eq:DeTurck-Willmore-flow} as a quasilinear parabolic system of PDEs for $u$ with homogeneous boundary conditions $u(t,\cdot)=d u(t,\cdot)=0$ on $[0,T]\times\partial\Sigma$.

In order to prove \Cref{thm:lifespan-and-smoothing}, we follow the strategy in \cite{lecronepruesswilke2014}. To this end, we rewrite the right-hand-side of the PDE in \eqref{eq:DeTurck-Willmore-flow} as $-A(f).f+F(f)$ where each $A(f)$ is a fourth-order elliptic operator and $F(f)$ captures the remaining lower-order terms with suitable ``singular'' Lipschitz control. In particular, weighted maximal $L^p$-regularity theory is used for the linearized problem to account for the low initial regularity. We repeat the steps in the fixed-point argument using the localized parabolic norms to verify that all constants can be chosen universally, as claimed in \Cref{thm:lifespan-and-smoothing}. Finally, the higher-order estimates \eqref{eq:curvature-smoothing} follow from a bootstrapping argument using classical Schauder theory.

More precisely, for $f\in C^1_{\mathrm{imm}}(\Sigma,\R^3)$ and $v\in W^{4,p}(\Sigma,\R^3)$, we define in local coordinates 
\begin{equation}\label{eq:def-A(f).v}
  A(f).v \vcentcolon= \frac12 g_{f}^{ab}g_{f}^{cd}\bar\nabla_{abcd}^4v
\end{equation}
where $\bar\nabla$ denotes covariant differentiation in $(\Sigma,g_{G_{f_0}})$. For $f\in W^{4,p}_{\mathrm{imm}}(\Sigma,\R^3)$, $F(f)$ is determined by the identity
\begin{equation}\label{eq:def-F(f)}
  \nabla\W(f)+(\Delta_{g_f}\mathbf W_{f,G_{f_0}})^k\partial_k f =\vcentcolon A(f).f - F(f).
\end{equation}
In particular, the PDE in \eqref{eq:DeTurck-Willmore-flow} can be rewritten as
\begin{equation}
  \partial_t f + A(f).f = F(f)\quad\text{a.e.\ on $(0,T)\times\Sigma$}.
\end{equation}
We next state some structural properties for $A(f)$ and $F(f)$ which are needed to apply a proof strategy motivated by the abstract result \cite[Theorem~2.1]{lecronepruesswilke2014}. 

\begin{lemma}\label{lem:structure-for-A-and-F}
  After reducing $\varepsilon_{\ref{thm:lifespan-and-smoothing}}$ universally, there exists a universal $0<\varepsilon_{\ref{lem:structure-for-A-and-F}}\leq\varepsilon_{\ref{prop:smooth-clamped-immersion}}$ such that, for $G=G_{f_0}$ as in \Cref{thm:lifespan-and-smoothing}, the following holds. Let $f\in W^{4(\mu-1/p),p}(\Sigma,\R^3)$ satisfy
  \begin{equation}\label{eq:1-lem:structure-for-A-and-F}
    \Vert f-G\Vert_{\bm W^{4(\mu-1/p),p}(\Sigma)}
    \leq\varepsilon_{\ref{lem:structure-for-A-and-F}}.
  \end{equation}
  Write $f_j=f\circ\alpha_j$, $v_j=v\circ\alpha_j$ and $g_{h,j}=\alpha_j^*g_h$ for $h=f,G$.
  \begin{enumerate}[(a)]
    \item For $v\in W^{4,p}(\Sigma,\R^3)$, the local representation of $A(f)$ is
    \begin{equation}
      (A(f).v)\circ\alpha_j
      =\frac12 g_{f,j}^{ab}g_{f,j}^{cd}D^4_{abcd}v_j
      +\sum_{|\nu|\leq3}c_{j,\nu}(f)D^\nu v_j.
    \end{equation}
    Thus $A(f)$ has the form \eqref{eq:linear-operator-in-local-coords}, with $a_j^{ab}=2^{-1/2}g_{f,j}^{ab}$ and $\alpha'$ in place of $\alpha$. All coefficients have universal $C^{\alpha'}(D_j)$-bounds, uniformly in $f$ and $j$.

    \item Let $\bar{\bm\Gamma}_j$ denote the Christoffel symbols of $g_{G,j}$. There are smooth coefficient maps $\mathcal T_3$, $\mathcal T_{2,3}$, $\mathcal T_{2,2,2}$, $\mathcal T_{2,2}$, $\mathcal T_2$ and $T$, independent of $f$ and $j$, such that, if additionally $f\in W^{4\beta,p}_{\mathrm{imm}}(\Sigma,\R^3)$, then
    \begin{align}
      F(f)\circ\alpha_j={}&\mathcal T_3*D^3f_j +\mathcal T_{2,3}*D^2f_j*D^3f_j +\mathcal T_{2,2,2}*D^2f_j*D^2f_j*D^2f_j \\
        &+\mathcal T_{2,2}*D^2f_j*D^2f_j + \mathcal T_2*D^2f_j+T \quad\text{a.e.\ on $D_j$}.\label{eq:2-lem:structure-for-A-and-F}
    \end{align}
    Here, every coefficient is evaluated at $(\bar{\bm\Gamma}_j(x),D\bar{\bm\Gamma}_j(x),Df_j(x))$, and $*$ denotes tensor contraction. Moreover, $(x,P)\mapsto\mathcal T_\bullet(\bar{\bm\Gamma}_j(x),D\bar{\bm\Gamma}_j(x),P)$ and $(x,P)\mapsto T(\bar{\bm\Gamma}_j(x),D\bar{\bm\Gamma}_j(x),P)$, have uniform $C^1$-bounds on the relevant argument ranges, depending only on $K$ in \eqref{eq:boundary-data-small}.
  \end{enumerate}
\end{lemma}
\begin{proof}
  Together, \Cref{prop:smooth-clamped-immersion}(b),(c), and \eqref{eq:1-lem:structure-for-A-and-F} with the embedding \eqref{eq:1-rem:choice-of-parameters:resulting-estimates} 
  yield for a universal $C>0$ that
  \begin{equation}\label{eq:1-pf-lem:structure-for-A-and-F}
    \Vert Df_j-DA_j\Vert_{C^{\alpha'}(D_j)}\leq C(\varepsilon_{\ref{lem:structure-for-A-and-F}}+\varepsilon_{\ref{thm:lifespan-and-smoothing}}),\quad\text{and thus also}\quad \Vert g_{f,j}-(\delta_{ab})\Vert_{C^{\alpha'}(D_j)}\leq C(\varepsilon_{\ref{lem:structure-for-A-and-F}}+\varepsilon_{\ref{thm:lifespan-and-smoothing}})
  \end{equation}
  Choose $\varepsilon_{\ref{lem:structure-for-A-and-F}}$ and $\varepsilon_{\ref{thm:lifespan-and-smoothing}}$ sufficiently small such that $C(\varepsilon_{\ref{lem:structure-for-A-and-F}}+ \varepsilon_{\ref{thm:lifespan-and-smoothing}} )<1/2$, so that $g_{f,j}$ is uniformly positive definite and the $C^{\alpha'}$-bounds of the inverses $g_{f,j}^{ab}$ are universal. Then writing \eqref{eq:def-A(f).v} in local coordinates yields (a). Indeed, the lower-order coefficients $c_{j,\nu}$ depend smoothly on $Df_j$ and polynomially on $\bar{\bm\Gamma}_j$, $D\bar{\bm\Gamma}_j$ and $D^2\bar{\bm\Gamma}_j$. Their uniform bounds follow from \Cref{prop:smooth-clamped-immersion}(c) and \eqref{eq:boundary-data-small}.

  For (b), we refer to the computations in \cite[Appendix~A]{palmurellariviere2022} in local coordinates. When applied to our setting, the fourth-order terms reduce to $\frac12g_{f,j}^{ab}g_{f,j}^{cd}D^4_{abcd}f_j$. Moreover, all terms involving $D^2\bar{\bm\Gamma}_j$ cancel in $F(f)\circ\alpha_j=(A(f).f-\nabla\W(f)-(\Delta_{g_f}\mathbf W_{f,G})^k\partial_k f)\circ\alpha_j$.
  Furthermore, for the Christoffel symbols $\bm\Gamma_{f,j}$ of $g_{f,j}$, we have
  \begin{equation}
     (\bm\Gamma_{f,j})^k_{ab}
    =g_{f,j}^{k\ell}\langle D^2_{ab}f_j,D_\ell f_j\rangle.
  \end{equation}
  Since $g_{f,j}$ and $\mathbf n_{f\circ\alpha_j}=\mathbf n_f\circ\alpha_j$ depend smoothly only on $Df_j$, expanding the remaining terms in $F(f)\circ\alpha_j$ and grouping the factors $D^2f_j$ and $D^3f_j$ indeed yields \eqref{eq:2-lem:structure-for-A-and-F}. Finally, using \Cref{prop:smooth-clamped-immersion}(c) and \eqref{eq:boundary-data-small}, we get the $C^1$-bounds for the coefficients in (b).
\end{proof}

We now verify that this multilinear structure for the lower-order terms in $F(G+u)$ allows for a singular\footnote{Singular in the sense that the Lipschitz constant explicitly depends on $\Vert u(t)\Vert_{W^{4\beta,p}(\Sigma)}$ which is not controlled for $t\searrow 0$ by $\Vert u\Vert_{\mathbb E_{1,\mu}(J_T)}$.} Lipschitz estimate which corresponds to established regimes for well-posedness of quasilinear problems in literature, see \cite[Equation~(2.4)]{lecronepruesswilke2014} in our specific weighted $L^p$-setting, or \cite[Assumption~(H2)]{lecronesimonett2020} in the setting of so-called \emph{continuous maximal regularity}.
\begin{lemma}[Localized Lipschitz estimate for $F(f)$]\label{lem:lipschitz-estimate-for-F}
  There exists a universal constant $C_{\ref{lem:lipschitz-estimate-for-F}}>0$ such that, for $f_i=G+u_i\in W^{4\beta,p}_{\mathrm{imm}}(\Sigma,\R^3)$ with $\Vert u_i\Vert_{\bm W^{4(\mu-1/p),p}(\Sigma)}\leq\varepsilon_{\ref{lem:structure-for-A-and-F}}$ for $i=1,2$, we have
  \begin{align}
    \Vert &(F(f_1)-F(f_2))\circ\alpha_k\Vert_{L^p(D_k)} \\
    &\leq C_{\ref{lem:lipschitz-estimate-for-F}}\sum_{j=1}^{3}\Big(1+\Vert u_1\circ\alpha_k\Vert_{W^{4\beta,p}(D_k)}^{\rho_j}+\Vert u_2\circ\alpha_k\Vert_{W^{4\beta,p}(D_k)}^{\rho_j}\Big)\Vert(f_1-f_2)\circ\alpha_k\Vert_{W^{4\beta_j,p}(D_k)}\label{eq:1-lem:lipschitz-estimate-for-F}
  \end{align}
  for all $k\in\mathbf I_{\mathrm{bdr}}\cup\mathbf I_{\mathrm{int}}$, where $\beta$ is chosen as in \eqref{eq:choice-of-parameters-2} and
  \begin{equation}\label{eq:2-lem:lipschitz-estimate-for-F}
    \{(\rho_j,\beta_j)\}_{j=1}^3=\{(1/2,3/4),(1,1/2),(1+\theta_{2,4\beta},\mu-1/p)\}
  \end{equation}
  with $\theta_{2,4\beta}$ from \Cref{rem:choice-of-parameters:resulting-estimates} such that
  \begin{equation}\label{eq:3-lem:lipschitz-estimate-for-F}
    \rho_j(\beta-(\mu-1/p))+\beta_j<1\quad\text{for all $j=1,2,3$}.
  \end{equation}
\end{lemma}
\begin{proof}
  Fix $k$ and write $f_{i,k}=f_i\circ\alpha_k$ and $u_{i,k}=u_i\circ\alpha_k$. By \eqref{eq:1-rem:unlocalized-coordinate-norms}, \eqref{eq:5-prop:smooth-clamped-immersion}, \eqref{eq:boundary-data-small} and $4(\mu-1/p)>1+\frac2p$, we have
  \begin{align}
    \Vert u_{i,k}\Vert_{W^{4(\mu-1/p),p}(D_k)}+\Vert u_{i,k}\Vert_{C^1(D_k)}&\leq C,\quad \Vert D(G\circ\alpha_k)\Vert_{C^2(D_k)}\leq C,\\
    \Vert f_{1,k}-f_{2,k}\Vert_{C^1(D_k)}&\leq C\Vert f_{1,k}-f_{2,k}\Vert_{W^{4(\mu-1/p),p}(D_k)}.\label{eq:1-pf-lem:lipschitz-estimate-for-F}
  \end{align}
  Standard interpolation and embedding estimates on $D_k$ yield that
  \begin{align}
    \Vert u_{i,k}\Vert_{W^{2,3p}(D_k)}&\leq C\Vert u_{i,k}\Vert_{W^{4\beta,p}(D_k)}^{\theta_{2,4\beta}}\Vert u_{i,k}\Vert_{W^{4(\mu-1/p),p}(D_k)}^{1-\theta_{2,4\beta}}\leq C\Vert u_{i,k}\Vert_{W^{4\beta,p}(D_k)}^{\theta_{2,4\beta}},\\
    \Vert u_{i,k}\Vert_{W^{3,\frac{3p}{2}}(D_k)}+\Vert u_{i,k}\Vert_{C^3(D_k)}&\leq C\Vert u_{i,k}\Vert_{W^{4\beta,p}(D_k)},\\
    \Vert u_{i,k}\Vert_{C^2(D_k)}^2&\leq C\Vert u_{i,k}\Vert_{C^1(D_k)}\Vert u_{i,k}\Vert_{C^3(D_k)}\leq C\Vert u_{i,k}\Vert_{W^{4\beta,p}(D_k)}.\label{eq:2-pf-lem:lipschitz-estimate-for-F}
  \end{align}
  Combining \eqref{eq:1-pf-lem:lipschitz-estimate-for-F} and \eqref{eq:2-pf-lem:lipschitz-estimate-for-F}, we get that
  \begin{align}
    \Vert D^2f_{i,k}\Vert_{C^0(D_k)}&\leq C\big(1+\Vert u_{i,k}\Vert_{W^{4\beta,p}(D_k)}^{1/2}\big),\\
    \Vert D^3f_{i,k}\Vert_{C^0(D_k)}+\Vert D^3f_{i,k}\Vert_{L^{\frac{3p}{2}}(D_k)}&\leq C\big(1+\Vert u_{i,k}\Vert_{W^{4\beta,p}(D_k)}\big),\\
    \Vert D^2f_{i,k}\Vert_{L^{3p}(D_k)}&\leq C\big(1+\Vert u_{i,k}\Vert_{W^{4\beta,p}(D_k)}^{\theta_{2,4\beta}}\big).\label{eq:3-pf-lem:lipschitz-estimate-for-F}
  \end{align}

  As in \eqref{eq:2-lem:structure-for-A-and-F}, the arguments denoted by $\cdot$ in $\mathcal T_\bullet(\cdot,Df_{i,k})$ and $T(\cdot,Df_{i,k})$ are $(\bar{\bm\Gamma}_k,D\bar{\bm\Gamma}_k)$. Since $(1-t)f_2+tf_1$ satisfies \eqref{eq:1-lem:structure-for-A-and-F} for $0\leq t\leq1$, the uniform $C^1$-bounds in \Cref{lem:structure-for-A-and-F}(b) yield that
  \begin{align}
    \Vert\mathcal T_\bullet(\cdot,Df_{i,k})\Vert_{C^0(D_k)}+\Vert T(\cdot,Df_{i,k})\Vert_{C^0(D_k)}&\leq C,\\
    \Vert\mathcal T_\bullet(\cdot,Df_{1,k})-\mathcal T_\bullet(\cdot,Df_{2,k})\Vert_{C^0(D_k)}&\leq C\Vert f_{1,k}-f_{2,k}\Vert_{C^1(D_k)},\\
    \Vert T(\cdot,Df_{1,k})-T(\cdot,Df_{2,k})\Vert_{C^0(D_k)}&\leq C\Vert f_{1,k}-f_{2,k}\Vert_{C^1(D_k)}.\label{eq:4-pf-lem:lipschitz-estimate-for-F}
  \end{align}
  By Hölder's inequality and \eqref{eq:1-pf-lem:lipschitz-estimate-for-F}, \eqref{eq:3-pf-lem:lipschitz-estimate-for-F} and \eqref{eq:4-pf-lem:lipschitz-estimate-for-F}, we have for the second term in \eqref{eq:2-lem:structure-for-A-and-F}
  \begin{align}
    \Vert &\mathcal T_{2,3}(\cdot,Df_{1,k})*D^2f_{1,k}*D^3f_{1,k}-\mathcal T_{2,3}(\cdot,Df_{2,k})*D^2f_{2,k}*D^3f_{2,k}\Vert_{L^p(D_k)}\\
    &\leq C\Vert D^2f_{1,k}\Vert_{L^{3p}(D_k)}\Vert D^3f_{1,k}\Vert_{L^{\frac{3p}{2}}(D_k)}\Vert f_{1,k}-f_{2,k}\Vert_{C^1(D_k)}\\
    &\quad+C\Vert D^3f_{1,k}\Vert_{C^0(D_k)}\Vert f_{1,k}-f_{2,k}\Vert_{W^{2,p}(D_k)}+C\Vert D^2f_{2,k}\Vert_{C^0(D_k)}\Vert f_{1,k}-f_{2,k}\Vert_{W^{3,p}(D_k)}\\
    &\leq C\big(1+\Vert u_{1,k}\Vert_{W^{4\beta,p}(D_k)}^{1+\theta_{2,4\beta}}\big)\Vert f_{1,k}-f_{2,k}\Vert_{W^{4(\mu-1/p),p}(D_k)}\\
    &\quad+C\big(1+\Vert u_{1,k}\Vert_{W^{4\beta,p}(D_k)}\big)\Vert f_{1,k}-f_{2,k}\Vert_{W^{2,p}(D_k)}\\
    &\quad+C\big(1+\Vert u_{2,k}\Vert_{W^{4\beta,p}(D_k)}^{1/2}\big)\Vert f_{1,k}-f_{2,k}\Vert_{W^{3,p}(D_k)}.\label{eq:5-pf-lem:lipschitz-estimate-for-F}
  \end{align}
  Similarly, for the third term, we have
  \begin{align}
    \Vert &\mathcal T_{2,2,2}(\cdot,Df_{1,k})*D^2f_{1,k}*D^2f_{1,k}*D^2f_{1,k}-\mathcal T_{2,2,2}(\cdot,Df_{2,k})*D^2f_{2,k}*D^2f_{2,k}*D^2f_{2,k}\Vert_{L^p(D_k)}\\
    &\leq C\big(\Vert D^2f_{1,k}\Vert_{C^0(D_k)}^2+\Vert D^2f_{2,k}\Vert_{C^0(D_k)}^2\big)\Vert f_{1,k}-f_{2,k}\Vert_{W^{2,p}(D_k)}\\
    &\quad+C\big(\Vert D^2f_{1,k}\Vert_{L^{3p}(D_k)}^3+\Vert D^2f_{2,k}\Vert_{L^{3p}(D_k)}^3\big)\Vert f_{1,k}-f_{2,k}\Vert_{C^1(D_k)}\\
    &\leq C\big(1+\Vert u_{1,k}\Vert_{W^{4\beta,p}(D_k)}+\Vert u_{2,k}\Vert_{W^{4\beta,p}(D_k)}\big)\Vert f_{1,k}-f_{2,k}\Vert_{W^{2,p}(D_k)}\\
    &\quad+C\big(1+\Vert u_{1,k}\Vert_{W^{4\beta,p}(D_k)}^{3\theta_{2,4\beta}}+\Vert u_{2,k}\Vert_{W^{4\beta,p}(D_k)}^{3\theta_{2,4\beta}}\big)\Vert f_{1,k}-f_{2,k}\Vert_{W^{4(\mu-1/p),p}(D_k)}.\label{eq:6-pf-lem:lipschitz-estimate-for-F}
  \end{align}
  For the remaining terms, \eqref{eq:1-pf-lem:lipschitz-estimate-for-F}, \eqref{eq:3-pf-lem:lipschitz-estimate-for-F} and \eqref{eq:4-pf-lem:lipschitz-estimate-for-F} yield that
  \begin{align}
    &\Vert\mathcal T_3(\cdot,Df_{1,k})*D^3f_{1,k}-\mathcal T_3(\cdot,Df_{2,k})*D^3f_{2,k}\Vert_{L^p(D_k)}\\
    &\quad+\Vert\mathcal T_{2,2}(\cdot,Df_{1,k})*D^2f_{1,k}*D^2f_{1,k}-\mathcal T_{2,2}(\cdot,Df_{2,k})*D^2f_{2,k}*D^2f_{2,k}\Vert_{L^p(D_k)}\\
    &\quad+\Vert\mathcal T_2(\cdot,Df_{1,k})*D^2f_{1,k}-\mathcal T_2(\cdot,Df_{2,k})*D^2f_{2,k}\Vert_{L^p(D_k)}\\
    &\quad+\Vert T(\cdot,Df_{1,k})-T(\cdot,Df_{2,k})\Vert_{L^p(D_k)}\\
    &\leq C\Vert f_{1,k}-f_{2,k}\Vert_{W^{3,p}(D_k)}\\
    &\quad+C\big(1+\Vert u_{1,k}\Vert_{W^{4\beta,p}(D_k)}^{1/2}+\Vert u_{2,k}\Vert_{W^{4\beta,p}(D_k)}^{1/2}\big)\Vert f_{1,k}-f_{2,k}\Vert_{W^{2,p}(D_k)}\\
    &\quad+C\big(1+\Vert u_{1,k}\Vert_{W^{4\beta,p}(D_k)}+\Vert u_{2,k}\Vert_{W^{4\beta,p}(D_k)}\big)\Vert f_{1,k}-f_{2,k}\Vert_{W^{4(\mu-1/p),p}(D_k)}.\label{eq:7-pf-lem:lipschitz-estimate-for-F}
  \end{align}
  Thus \eqref{eq:2-lem:structure-for-A-and-F}, \eqref{eq:5-pf-lem:lipschitz-estimate-for-F} and \eqref{eq:7-pf-lem:lipschitz-estimate-for-F}, and $3\theta_{2,4\beta}<1+\theta_{2,4\beta}$ yield \eqref{eq:1-lem:lipschitz-estimate-for-F} with \eqref{eq:2-lem:lipschitz-estimate-for-F}. Finally, \Cref{rem:choice-of-parameters:resulting-estimates} yields \eqref{eq:3-lem:lipschitz-estimate-for-F}.
\end{proof}

\begin{proof}[Proof of \Cref{thm:lifespan-and-smoothing}]
  We write $G=G_{f_0}$ and $\widetilde f_0=f_0\circ\Xi_{f_0}$. Let
  $$V_\mu=\{f\in W^{4(\mu-1/p),p}(\Sigma,\R^3)\mid \Vert f-G\Vert_{\bm W^{4(\mu-1/p),p}(\Sigma)}<\varepsilon_{\ref{lem:structure-for-A-and-F}}\}\subset W^{4(\mu-1/p),p}_{\mathrm{imm}}(\Sigma,\R^3).$$
  We can choose $\varepsilon_{\ref{thm:lifespan-and-smoothing}}>0$ sufficiently small, only depending on universal constants, such that \Cref{prop:smooth-clamped-immersion} and \eqref{eq:1-rem:choice-of-parameters:resulting-estimates} yield $\widetilde f_0\in V_\mu$. Using \Cref{rem:unlocalized-coordinate-norms}, \eqref{eq:def-A(f).v}, \Cref{rem:choice-of-parameters:resulting-estimates}, and \Cref{prop:smooth-clamped-immersion}(c), for all $f_1,f_2\in V_\mu$ and $w\in W^{4,p}(\Sigma,\R^3)$, we have
  \begin{align}
    \Vert((A(f_1)-A(f_2)).w)\circ\alpha_k\Vert_{L^p(D_k)} &\leq C\Vert (f_1-f_2)\circ\alpha_k\Vert_{W^{4(\mu-1/p),p}(D_k)}\Vert w\circ\alpha_k\Vert_{W^{4,p}(D_k)},\\
    \Vert((A(f_1)-A(f_2)).G)\circ\alpha_k\Vert_{L^p(D_k)} &\leq C\Vert (f_1-f_2)\circ\alpha_k\Vert_{W^{4(\mu-1/p),p}(D_k)}\label{eq:2-proof-prop:lifespan-estimate}
  \end{align}
  for all $k\in\mathbf{I}_{\mathrm{int}}\cup\mathbf{I}_{\mathrm{bdr}}$. In the second estimate, the affine functions $A_j(\cdot,0)$ in \eqref{eq:5-prop:smooth-clamped-immersion} can be disregarded since only derivatives of $G$ occur. Thus, here and below, $C$ is a universal constant changing from line to line. 

  Let $u_0=\widetilde f_0-G$. By \Cref{prop:smooth-clamped-immersion}(a), we have $u_0=du_0=0$ on $\partial\Sigma$. By \Cref{rem:choice-of-parameters:resulting-estimates}, we have $1+1/p<4(\mu-1/p)<1+\alpha<2$, so using \Cref{lem:structure-for-A-and-F}(a), \Cref{prop:max-L^p-reg} applies to $2A(G)$ and yields that there exists $u^\ast\in\mathbb E_{1,\mu}(J_1)$ satisfying
  \begin{equation}\label{eq:3-proof-prop:lifespan-estimate}
    \begin{cases}
      \partial_tu^\ast+A(G).u^\ast=0&\text{a.e.\ on $[0,1]\times\Sigma$}\\
      u^\ast(0,\cdot)=u_0&\text{on $\Sigma$}\\
      u^\ast=du^\ast=0&\text{on $[0,1]\times\partial\Sigma$}.
    \end{cases}
  \end{equation}
  Moreover,
  \begin{equation}\label{eq:4-proof-prop:lifespan-estimate}
    \Vert u^\ast\Vert_{\mathbb E_{1,\mu}(J_1)}+\sup_{0\leq t\leq 1}\Vert u^\ast(t)\Vert_{\bm W^{4(\mu-1/p),p}(\Sigma)}\leq C\Vert u_0\Vert_{\bm W^{4(\mu-1/p),p}(\Sigma)}\leq C\varepsilon.
  \end{equation}
  Set $\bar f=G+u^\ast$. For $\tau_0\in(0,1]$ and $r>0$, define
  \begin{equation}\label{eq:5-proof-prop:lifespan-estimate}
    \mathbb B_r=\{v\in\mathbb E_{1,\mu}(J_{\tau_0})\mid v(0,\cdot)=0,\operatorname{tr}_{\partial\Sigma}v=\operatorname{tr}_{\partial\Sigma}dv=0,\text{ and }\Vert v\Vert_{\mathbb E_{1,\mu}(J_{\tau_0})}\leq r\}.
  \end{equation}
  \Cref{prop:BUC-estimate} with $v(0)=0$ and \eqref{eq:4-proof-prop:lifespan-estimate} yield
  \begin{equation}\label{eq:6-proof-prop:lifespan-estimate}
    \Vert u^\ast+v\Vert_{\mathbb E_{1,\mu}(J_{\tau_0})}+\sup_{0\leq t\leq\tau_0}\Vert u^\ast(t)+v(t)\Vert_{\bm W^{4(\mu-1/p),p}(\Sigma)}\leq C(r+\varepsilon)
  \end{equation}
  for all $v\in\mathbb B_r$, with a constant independent of $\tau_0$. We choose $r+\varepsilon$ sufficiently small such that $\bar f(t)+v(t)\in V_\mu$ for all $t\in J_{\tau_0}$.

  We next argue that for $\kappa>0$ chosen below, we have
  \begin{equation}\label{eq:8-proof-prop:lifespan-estimate}
    \Vert F(\bar f+v_1)-F(\bar f+v_2)\Vert_{\bm L^p_\mu((0,\tau_0),\bm L^p(\Sigma))}\leq C\tau_0^\kappa\Vert v_1-v_2\Vert_{\mathbb E_{1,\mu}(J_{\tau_0})}\quad\text{for all $v_1,v_2\in\mathbb B_r$}.
  \end{equation}
  To this end, set
  \begin{equation}
    \vartheta=\frac{\beta-(\mu-1/p)}{1-(\mu-1/p)}\quad\text{and}\quad\vartheta_j=\frac{\beta_j-(\mu-1/p)}{1-(\mu-1/p)}.
  \end{equation}
  Then $0<\vartheta<1$, $0\leq\vartheta_j<1$, and \eqref{eq:3-lem:lipschitz-estimate-for-F} yields that $\rho_j\vartheta+\vartheta_j<1$. Choose
  \begin{equation}
    0<\kappa<\min_{j=1,2,3}\{(1-(\mu-1/p))(1-\vartheta_j),(1-(\mu-1/p))(1-\rho_j\vartheta-\vartheta_j)\}.
  \end{equation}
  For $v_1,v_2\in\mathbb B_r$, write $h_{i,k}=(u^\ast+v_i)\circ\alpha_k$ and $w_k=(v_1-v_2)\circ\alpha_k$. For more details on the following two estimates, see \cite[Equation~(2.13) and the estimates thereafter]{lecronepruesswilke2014}. Applying interpolation estimates on $D_k$, Hölder's inequality in $t$ with respect to the measure $t^{p(1-\mu)}\dd t$, and $\eqref{eq:6-proof-prop:lifespan-estimate}$ and \Cref{prop:BUC-estimate} with $v_i(0)=0$ and $w(0)=0$, we get
  \begin{align}
    &\max_{k\in\mathbf I_{\mathrm{int}}\cup\mathbf I_{\mathrm{bdr}}}\Big\Vert\Vert h_{i,k}(t)\Vert_{W^{4\beta,p}(D_k)}^{\rho_j}\Vert w_k(t)\Vert_{W^{4\beta_j,p}(D_k)}\Big\Vert_{L^p_\mu(0,\tau_0)}\\
    &\quad\leq C\tau_0^{(1-(\mu-1/p))(1-\rho_j\vartheta-\vartheta_j)}
    \Big(\sup_{t\in J_{\tau_0}}\Vert h_i(t)\Vert_{W^{4(\mu-1/p),p}(\Sigma)}\Big)^{\rho_j(1-\vartheta)}
    \Vert h_i\Vert_{\mathbb E_{1,\mu}(J_{\tau_0})}^{\rho_j\vartheta}\\
    &\quad\hspace{3em}\cdot\Big(\sup_{t\in J_{\tau_0}}\Vert w(t)\Vert_{W^{4(\mu-1/p),p}(\Sigma)}\Big)^{1-\vartheta_j}
    \Vert w\Vert_{\mathbb E_{1,\mu}(J_{\tau_0})}^{\vartheta_j} \leq C\tau_0^\kappa(r+\varepsilon)^{\rho_j}\Vert w\Vert_{\mathbb E_{1,\mu}(J_{\tau_0})}.
  \end{align}
  Proceeding similarly, one estimates
  \begin{equation}
    \Vert w\Vert_{\bm L^p_\mu((0,\tau_0),\bm W^{4\beta_j,p}(\Sigma))}\leq C\tau_0^{(1-(\mu-1/p))(1-\vartheta_j)}\Vert w\Vert_{\mathbb E_{1,\mu}(J_{\tau_0})}\leq C\tau_0^{\kappa}\Vert w\Vert_{\mathbb E_{1,\mu}(J_{\tau_0})}.
  \end{equation}
  With \Cref{lem:lipschitz-estimate-for-F} and \Cref{rem:unlocalized-coordinate-norms}, \eqref{eq:8-proof-prop:lifespan-estimate} follows provided $r+\varepsilon\leq 1$. For $v\in\mathbb B_r$, applying the same interpolation estimates to $h=u^\ast+v$ and $0$, using \eqref{eq:6-proof-prop:lifespan-estimate}, we get that
  \begin{equation}\label{eq:F-relative-to-G}
    \Vert F(\bar f+v)-F(G)\Vert_{\bm L^p_\mu((0,\tau_0),\bm L^p(\Sigma))}\leq C\tau_0^\kappa(r+\varepsilon).
  \end{equation}
  Since $\mathbf W_{G,G}=0$, the residual
  \begin{equation}
    \mathcal R_G=F(G)-A(G).G=-\frac12\big(\Delta_{g_G}\mathbf H_G+\vert\mathbf A_G^0\vert^2\mathbf H_G\big)\mathbf n_G
  \end{equation}
  satisfies $\Vert\mathcal R_G\Vert_{\bm L^p(\Sigma)}\leq C\varepsilon$ by \Cref{prop:smooth-clamped-immersion}(c) and \eqref{eq:boundary-data-small}. For $v\in\mathbb B_r$, let $\mathcal T(v)$ be the solution $u$ of
  \begin{equation}\label{eq:7-proof-prop:lifespan-estimate}
    \begin{cases}
      \partial_tu+A(G).u=F(\bar f+v)-F(G)+(A(G)-A(\bar f+v)).(\bar f+v)+\mathcal R_G &\text{a.e.\ on $[0,\tau_0]\times\Sigma$}\\
      u(0,\cdot)=0&\text{on $\Sigma$}\\
      \operatorname{tr}_{\partial\Sigma}u=\operatorname{tr}_{\partial\Sigma}du=0&\text{on $[0,\tau_0]\times\partial\Sigma$}.
    \end{cases}
  \end{equation}
  By the preceding estimates and again \Cref{prop:max-L^p-reg} applied to $2A(G)$, this map is well-defined. By \eqref{eq:3-proof-prop:lifespan-estimate}, a fixed point $v=\mathcal T(v)$ yields a solution $f=\bar f+v$ of \eqref{eq:DeTurck-Willmore-flow}. More precisely, the maximal regularity estimate yields
  \begin{align}\label{eq:9-proof-prop:lifespan-estimate}
    \Vert\mathcal T(v)\Vert_{\mathbb E_{1,\mu}(J_{\tau_0})}\leq C_{\mathrm{MR}}\big(&\Vert F(\bar f+v)-F(G)\Vert_{\bm L^p_\mu((0,\tau_0),\bm L^p(\Sigma))}\nonumber\\
    &+\Vert(A(G)-A(\bar f+v)).(\bar f+v)\Vert_{\bm L^p_\mu((0,\tau_0),\bm L^p(\Sigma))}\nonumber\\
    &+\Vert\mathcal R_G\Vert_{\bm L^p_\mu((0,\tau_0),\bm L^p(\Sigma))}\big).
  \end{align}
  Set $a=1-\mu+1/p>0$, and write $\bar f+v=G+h$. First applying \eqref{eq:2-proof-prop:lifespan-estimate} at a fixed $k\in\mathbf{I}_{\mathrm{int}}\cup\mathbf{I}_{\mathrm{bdr}}$ and then taking $\Vert\cdot\Vert_{L^p_\mu(0,\tau_0)}$, \eqref{eq:6-proof-prop:lifespan-estimate} implies that
  \begin{equation}
    \Vert(A(G)-A(G+h)).(G+h)\Vert_{\bm L^p_\mu((0,\tau_0),\bm L^p(\Sigma))}\leq C(r+\varepsilon)(\tau_0^a+r+\varepsilon).
  \end{equation}
  Moreover, $\Vert\mathcal R_G\Vert_{\bm L^p_\mu((0,\tau_0),\bm L^p(\Sigma))}\leq C\tau_0^a\varepsilon$. Thus, \eqref{eq:F-relative-to-G} yields that
  \begin{equation}\label{eq:10-proof-prop:lifespan-estimate}
    \Vert\mathcal T(v)\Vert_{\mathbb E_{1,\mu}(J_{\tau_0})}\leq C\big((\tau_0^\kappa+\tau_0^a+r+\varepsilon)(r+\varepsilon)+\tau_0^a\varepsilon\big).
  \end{equation}

  For $v_1,v_2\in\mathbb B_r$, set $f_i=\bar f+v_i$. Using
  \begin{equation}
    (A(G)-A(f_1)).f_1-(A(G)-A(f_2)).f_2=(A(G)-A(f_1)).(v_1-v_2)+(A(f_2)-A(f_1)).f_2
  \end{equation}
  and \Cref{prop:BUC-estimate} applied to $v_1-v_2$, one gets that
  \begin{align}
    \Vert(A(G)-A(f_1)).f_1-(A(G)-A(f_2)).f_2\Vert_{\bm L^p_\mu((0,\tau_0),\bm L^p(\Sigma))}\leq C(r+\varepsilon+\tau_0^a)\Vert v_1-v_2\Vert_{\mathbb E_{1,\mu}(J_{\tau_0})}.
  \end{align}
  Combined with \eqref{eq:8-proof-prop:lifespan-estimate} and \Cref{prop:max-L^p-reg}, we obtain
  \begin{equation}\label{eq:11-proof-prop:lifespan-estimate}
    \Vert\mathcal T(v_1)-\mathcal T(v_2)\Vert_{\mathbb E_{1,\mu}(J_{\tau_0})}\leq C(\tau_0^\kappa+\tau_0^a+r+\varepsilon)\Vert v_1-v_2\Vert_{\mathbb E_{1,\mu}(J_{\tau_0})}.
  \end{equation}
  We now fix $r=\varepsilon$. First choose $\tau_0>0$ sufficiently small, and then $\varepsilon_{\ref{thm:lifespan-and-smoothing}}>0$ sufficiently small, only depending on universal constants. Equations \eqref{eq:6-proof-prop:lifespan-estimate}, \eqref{eq:10-proof-prop:lifespan-estimate}, and \eqref{eq:11-proof-prop:lifespan-estimate} then show that $\bar f+v\in V_\mu$, that $\mathcal T$ maps $\mathbb B_\varepsilon$ into itself, and that its Lipschitz constant is at most $1/2$. Thus, Banach's fixed-point theorem yields a fixed point $v\in\mathbb B_\varepsilon$.

  Set $f=G+u^\ast+v$ and $T=\tau_{\ref{thm:lifespan-and-smoothing}}=\tau_0$. By \eqref{eq:4-proof-prop:lifespan-estimate},
  \begin{equation}
    \Vert f-G\Vert_{\mathbb E_{1,\mu}(J_{\tau_0})}\leq C\varepsilon.
  \end{equation}
  Further, the trace estimates and \eqref{eq:4-proof-prop:lifespan-estimate} yield that
  \begin{align}
    \sup_{0\leq t\leq\tau_0}&\Vert f(t)-\widetilde f_0\Vert_{\bm W^{4(\mu-1/p),p}(\Sigma)}\\
    &\leq\sup_{0\leq t\leq\tau_0}\Vert v(t)\Vert_{\bm W^{4(\mu-1/p),p}(\Sigma)}
    +\sup_{0\leq t\leq\tau_0}\Vert u^\ast(t)\Vert_{\bm W^{4(\mu-1/p),p}(\Sigma)}
    +\Vert u_0\Vert_{\bm W^{4(\mu-1/p),p}(\Sigma)}\leq C\varepsilon.
  \end{align}
  The embedding in \eqref{eq:1-rem:choice-of-parameters:resulting-estimates} proves \eqref{eq:2-thm:lifespan-and-smoothing}. \Cref{prop:BUC-estimate} also yields $f\in C^0([0,\tau_0],C^{1+\alpha'}(\Sigma,\R^3))$, and $f(t)\in V_\mu$ ensures that each $f(t)$ is an immersion. Since $G$ is smooth and independent of time, \eqref{eq:1-thm:lifespan-and-smoothing} follows. Finally, $f=G=\widetilde f_0$ and $df=dG$ on $[0,\tau_0]\times\partial\Sigma$. Thus, $\partial_{\eta_f}f=\partial_{\eta_G}G=\mathbf N\circ G=\mathbf N\circ f$ there. Taking the normal component of \eqref{eq:DeTurck-Willmore-flow} cancels the tangential term and yields the PDE in \Cref{thm:lifespan-and-smoothing}.

    We finally prove \eqref{eq:curvature-smoothing}. Let $k\geq 5$ such that \eqref{eq:boundary-data-small} holds, and put $u=f-G$ and $I_\sigma=[\sigma,\tau_0]$ for $0<\sigma<\tau_0$. For noninteger $s>0$, we use the localized parabolic H\"older norm
  \begin{equation}
    \Vert h\Vert_{\bm{\mathcal H}^s(I)}\vcentcolon=\max_{j\in\mathbf I_{\mathrm{int}}\cup\mathbf I_{\mathrm{bdr}}}\Vert(\pi_jh)\circ\alpha_j\Vert_{C^{s/4,s}(I\times\mathbb X_j)}.
  \end{equation}
  By \Cref{prop:smooth-clamped-immersion}(c), we have
  \begin{equation}\label{eq:smoothing-reference-immersion}
    \max_{j\in\mathbf I_{\mathrm{int}}\cup\mathbf I_{\mathrm{bdr}}}\Vert G\circ\alpha_j-A_j(\cdot,0)\Vert_{C^k(D_j)}\leq C_k\varepsilon\quad\text{and}\quad\Vert\mathcal R_G\Vert_{\bm{\mathcal H}^{k-5+\alpha'}(I)}\leq C_k\varepsilon
  \end{equation}
  for all closed intervals $I\subset[0,\tau_0]$, where $\mathcal R_G=F(G)-A(G).G$ is the residual from above. For the second estimate, recall that $\mathcal R_G$ is independent of time, involves at most four derivatives of $G$, and vanishes for an affine immersion.

  We first establish a parabolic H\"older estimate away from $t=0$. Since $p>10/(3\alpha)>3$, choose
  \begin{equation}
    0<\gamma<\min\{\alpha',2-6/p\}.
  \end{equation}
  On $I_{\tau_0/8}$ the time weight is bounded above and below by universal positive constants. Thus, the estimate for $u$ in $\mathbb E_{1,\mu}(J_{\tau_0})$ and, using  $2+\gamma<4-6/p$, the parabolic Morrey embedding (see \cite[Theorem~3.9.1 and Example~3.7.4, with $\omega=(1,1,4)$, $s=4$ and $t=2+\gamma$]{amann2009}) applied on each chart yields
  \begin{equation}\label{eq:smoothing-first-embedding}
    \Vert u\Vert_{\bm W^{1,p}(I_{\tau_0/8};\bm L^p(\Sigma))\cap\bm L^p(I_{\tau_0/8};\bm W^{4,p}(\Sigma))}+\Vert u\Vert_{\bm{\mathcal H}^{2+\gamma}(I_{\tau_0/8})}\leq C\varepsilon.
  \end{equation}
  Write $\mathcal R(f)=F(f)-A(f).G$. We use that the terms in \eqref{eq:2-lem:structure-for-A-and-F} are affine in $D^3f_j$, where $f_j=f\circ\alpha_j$. Subtracting $\mathcal R(G)=\mathcal R_G$, and first collecting all terms containing $D^3u_j$, shows that there is a linear differential operator $B(t)$ of order at most three such that
  \begin{equation}\label{eq:smoothing-lower-order-operator}
    \mathcal R(f)=\mathcal R_G+B(t).u\quad\text{and}\quad(B(t).v)\circ\alpha_j=\sum_{|\nu|\leq3}b_{j,\nu}(t,x)D^\nu(v\circ\alpha_j).
  \end{equation}
  The coefficients are matrix-valued smooth functions of $Df_j$, $D^2f_j$ and derivatives of $G_j=G\circ\alpha_j$ of orders one through four, and do not involve $D^3f_j$. To see this last point, a term of the form $\mathcal T(Df_j)*D^2f_j*D^3f_j$ is first written as $\mathcal T(Df_j)*D^2f_j*D^3u_j+\mathcal T(Df_j)*D^2f_j*D^3G_j$. In the second term, subtraction of its value at $G_j$ and the mean value formula only produce coefficients multiplying $Du_j$ and $D^2u_j$. The other terms are treated in the same way. By \eqref{eq:smoothing-reference-immersion} and \eqref{eq:smoothing-first-embedding}, all coefficients of
  \begin{equation}
    L(t)=A(f(t))-B(t)
  \end{equation}
  have uniform $C^{\gamma/4,\gamma}$-bounds on $I_{\tau_0/8}\times D_j$. Its principal symbol is $\frac12(g_f^{ab}\xi_a\xi_b)^2\operatorname{id}_{\R^3}$, which is uniformly positive for $\xi\neq0$. The homogeneous clamped boundary operators satisfy the Lopatinskii--Shapiro condition. Choose a smooth $\zeta$ with $\chi_{I_{\tau_0/4}}\leq\zeta\leq\chi_{I_{\tau_0/8}}$, and set $w=\zeta u$. Then
  \begin{equation}\label{eq:smoothing-cutoff-problem}
    \begin{cases}
      \partial_tw+L(t).w=\zeta\mathcal R_G+\zeta'u&\text{a.e.\ on $I_{\tau_0/8}\times\Sigma$}\\
      w(\tau_0/8,\cdot)=0&\text{on $\Sigma$}\\
      w=dw=0&\text{on $I_{\tau_0/8}\times\partial\Sigma$},
    \end{cases}
  \end{equation}
  using \eqref{eq:smoothing-lower-order-operator} to rewrite the PDE in the first line. Its right-hand-side has $\bm{\mathcal H}^\gamma$-norm at most $C\varepsilon$. Proceeding by localization as in the proof of \Cref{prop:max-L^p-reg}, and using \cite[Theorem~4.9]{solonnikov1965mathphysics} on the localized versions of \eqref{eq:smoothing-cutoff-problem}, therefore yields that $w\in\bm{\mathcal H}^{4+\gamma}(I_{\tau_0/8})$ and
  \begin{equation}\label{eq:smoothing-first-schauder}
    \Vert u\Vert_{\bm{\mathcal H}^{4+\gamma}(I_{\tau_0/4})}\leq C\big(\Vert\mathcal R_G\Vert_{\bm{\mathcal H}^\gamma(I_{\tau_0/8})}+\Vert u\Vert_{\bm{\mathcal H}^{2+\gamma}(I_{\tau_0/8})}\big)\leq C\varepsilon.
  \end{equation}
  Here and below, we use that $w$ vanishes near the initial time, all initial compatibility conditions are satisfied, and we use uniqueness of solutions in $\mathbb E_{1,1}(I_{\tau_0/4})$.

  We next improve the H\"older exponent to $\alpha'$. Using $4+\gamma-3>\alpha'$ and \eqref{eq:smoothing-first-schauder}, \eqref{eq:smoothing-reference-immersion} and the smooth dependence of the coefficients, we obtain on $I_{\tau_0/4}$ that the coefficients of $A(f)$ have uniform $C^{\alpha'/4,\alpha'}$-bounds and
  \begin{equation}
    \Vert\mathcal R(f)\Vert_{\bm{\mathcal H}^{\alpha'}(I_{\tau_0/4})}\leq C\varepsilon.
  \end{equation}
  Localizing and applying \cite[Theorem~4.9]{solonnikov1965mathphysics} to the analog of \eqref{eq:smoothing-cutoff-problem}, only with a $\zeta$ such that $\chi_{I_{3\tau_0/8}}\leq\zeta\leq \chi_{I_{\tau_0/4}}$, yields
  \begin{equation}\label{eq:smoothing-base-schauder}
    \Vert u\Vert_{\bm{\mathcal H}^{4+\alpha'}(I_{3\tau_0/8})}\leq C\varepsilon.
  \end{equation}

  For $k>5$, choose equidistant times $3\tau_0/8=a_4<a_5<\dots<a_{k-1}<\tau_0/2$, suppose that $5\leq m\leq k-1$ and that $\Vert u\Vert_{\bm{\mathcal H}^{m-1+\alpha'}(I_{a_{m-1}})}\leq C_k\varepsilon$. The expression for $\mathcal R(f)$ contains at most three derivatives of $u$ and four derivatives of $G$. Consequently, \eqref{eq:smoothing-reference-immersion} and the induction hypothesis yield that
  \begin{equation}
    \Vert\mathcal R(f)\Vert_{\bm{\mathcal H}^{m-4+\alpha'}(I_{a_{m-1}})}\leq C_k\varepsilon.
  \end{equation}
  Again, localizing and applying \cite[Theorem~4.9]{solonnikov1965mathphysics} to \eqref{eq:smoothing-cutoff-problem}, only with $\zeta$ such that $\chi_{I_{a_m}}\leq\zeta\leq \chi_{I_{a_{m-1}}}$, yields that $\Vert u\Vert_{\bm{\mathcal H}^{m+\alpha'}(I_{a_m})}\leq C_k\varepsilon$. Thus, by induction, together with \eqref{eq:smoothing-base-schauder} for $k=5$, we get
  \begin{equation}\label{eq:smoothing-coordinate-estimate}
    \sup_{\tau_0/2\leq t\leq\tau_0}\Vert f(t)-G\Vert_{\bm C^{k-1+\alpha'}(\Sigma)}\leq C_k\varepsilon.
  \end{equation}
  Finally, \eqref{eq:smoothing-reference-immersion}, \eqref{eq:smoothing-coordinate-estimate} and \Cref{rem:unlocalized-coordinate-norms} imply that
  \begin{equation}
    \sup_{\tau_0/2\leq t\leq\tau_0}\max_{j\in\mathbf I_{\mathrm{int}}\cup\mathbf I_{\mathrm{bdr}}}\Vert f(t)\circ\alpha_j-A_j(\cdot,0)\Vert_{C^{k-1+\alpha'}(D_j)}\leq C_k\varepsilon,
  \end{equation}
  and therefore also \eqref{eq:curvature-smoothing}.
\end{proof}

\section{Existence and rigidity of stationary blow-ups}
\label{sec:pf-intro-results}

\subsection{Construction of blow-ups}
\label{sec:blow-up-construction}

In this section, we prove the following version of \Cref{thm:intro-stationary-blow-ups}. For the proof, we use Simon's monotonicity formula, see \Cref{app:monotonicity-formula} for more details.
\begin{theorem}[Existence of Willmore blow-ups]\label{prop:existence-of-blow-ups}
  For a maximal solution $f$ of \eqref{eq:Dirichlet-Willmore-flow}, exactly one of the following two cases occurs.
  \begin{enumerate}[(a)]
    \item $T=\infty$ and $\displaystyle\sup_{t\in [0,\infty)}\Vert \mathbf A_{f(t)}\Vert_{L^\infty(\Sigma)}<\infty$.
    \item There exist $t_j\nearrow T$, $r_j\searrow0$, $p_j\in\Sigma$ and a smooth proper immersion $\hat f\colon\hat\Sigma\to\R^3$ of a complete connected oriented surface, with $\hat p\in\hat\Sigma$ and $\hat f(\hat p)=0$, such that
    \begin{equation}\label{eq:pointed-convergence-blow-ups}
      \Big(\frac{f(t_j)-f(t_j,p_j)}{r_j},p_j\Big)\longrightarrow(\hat f,\hat p)
    \end{equation}
    in the sense of local smooth pointed convergence of \Cref{def:local-smooth-pointed-convergence}. Moreover,
    \begin{align}
      &\qquad \Delta_{g_{\hat f}}\mathbf H_{\hat f}+|\mathbf A_{\hat f}^0|^2\mathbf H_{\hat f}=0,\quad |\mathbf A_{\hat f}(\hat p)|>0,\\
      &\int_{\hat\Sigma}|\mathbf A_{\hat f}|^2\dd\vol_{g_{\hat f}}<\infty,\quad [\mathbf n_{\hat f}]_{\alpha,1}<\infty,\quad\text{and}\quad \mathrm{genus}(\hat\Sigma)\leq\mathrm{genus}(\Sigma).\label{eq:properties-stationary-blow-up}
    \end{align}
    Either $\partial\hat\Sigma=\emptyset$, or, after applying a Euclidean isometry, $\hat f\vert_{\partial\hat\Sigma}$ is an embedding onto $\R e_1$ and $\partial_{\eta_{\hat f}}\hat f=(0,-1,0)$ on $\partial\hat\Sigma$.
  \end{enumerate}
  In case (b), the density at infinity $\theta^2(\hat f,\infty)\vcentcolon=\lim_{R\to\infty}\frac{\int_{\hat f^{-1}(B_R(0))}\dd\vol_{g_{\hat f}}}{\pi R^2}$ exists and we have
  \begin{align}
      \theta^2(\hat f,\infty)+\frac{\W(\hat f)}{4\pi}&\leq\frac{\W(f_0)}{4\pi}+\frac1{2\pi}\sup_{z\in\R^3}\int_\Gamma\frac{\langle\mathbf N(p),p-z\rangle}{|p-z|^2}\dd\mathcal H^1(p)\quad\text{if $\partial\hat\Sigma=\emptyset$},\\
      \theta^2(\hat f,\infty)+\frac{\W(\hat f)}{4\pi} &\leq\frac{\W(f_0)}{4\pi}+\frac1{2\pi}\sup_{z\in\Gamma}\int_\Gamma\frac{\langle\mathbf N(p),p-z\rangle}{|p-z|^2}\dd\mathcal H^1(p)\quad\text{if $\partial\hat\Sigma\neq\emptyset$}.\label{eq:density-bound-stationary-blow-up}
  \end{align}
\end{theorem}
\begin{proof}
  Using \eqref{eq:energy-dissipation-Dirichlet}, $\W(f(t))$ is non-increasing. Moreover, by Gauss--Bonnet,
  \begin{equation}\label{eq:total-curvature-Dirichlet}
    \int_\Sigma|\mathbf A_{f(t)}|^2\dd\vol_{g_{f(t)}}=4\W(f(t))-4\pi\chi(\Sigma)-2\int_\Gamma\langle\bm\kappa_\Gamma,\mathbf N\rangle\dd\mathcal H^1\leq C.
  \end{equation}
  Here and below, constants may depend on the initial immersion and on the fixed boundary data, but not on $j$. Fix $K=1$ in \Cref{thm:lifespan-and-smoothing}. We first claim that, if (a) does not hold, then
  \begin{equation}\label{eq:curvature-unbounded-blow-ups}
    \sup_{t\in[0,T)}\Vert \mathbf A_{f(t)}\Vert_{L^\infty(\Sigma)}=\infty.
  \end{equation}
  Otherwise, necessarily $T<\infty$. By \Cref{rem:Af-Linfty-Lipschitz-est}, choosing $r>0$ sufficiently small, the initial immersions $r^{-1}f(s)$ and their boundary data satisfy the assumptions of \Cref{thm:lifespan-and-smoothing}, uniformly for $s\in[0,T)$. Using \Cref{prop:geometric-uniqueness-Dirichlet}, maximality of $f$ implies
  \begin{equation}
    T-s\geq r^4\tau_{\ref{thm:lifespan-and-smoothing}}\quad\text{for all $s\in[0,T)$},
  \end{equation}
  a contradiction as $s\nearrow T$. This proves \eqref{eq:curvature-unbounded-blow-ups}.

  Choose $\delta>0$ sufficiently small such that
  \begin{equation}\label{eq:choice-delta-blow-ups}
    2\delta\leq\varepsilon_{\ref{rem:Af-Linfty-Lipschitz-est}}\quad\text{and}\quad 2C_{\ref{rem:Af-Linfty-Lipschitz-est}}^2\delta\leq\varepsilon_{\ref{thm:lifespan-and-smoothing}}.
  \end{equation}
  By \eqref{eq:curvature-unbounded-blow-ups}, there exist $\tilde t_j\nearrow T$ and $p_j\in\Sigma$ such that
  \begin{equation}\label{eq:record-times-blow-ups}
      |\mathbf A_{f(\tilde t_j)}(p_j)|=\max_{0\leq t\leq\tilde t_j}\Vert \mathbf A_{f(t)}\Vert_{L^\infty(\Sigma)}\longrightarrow\infty,\quad\text{and}\quad r_j\vcentcolon=\frac{\delta}{|\mathbf A_{f(\tilde t_j)}(p_j)|}\searrow0,\quad \tilde t_j-r_j^4\tau_{\ref{thm:lifespan-and-smoothing}}>0.
  \end{equation}
  For $0\leq t\leq\tau_{\ref{thm:lifespan-and-smoothing}}$, set
  \begin{equation}\label{eq:rescaled-flows-blow-ups}
    f_j(t,x)\vcentcolon=\frac{f\big(\tilde t_j+r_j^4(t-\tau_{\ref{thm:lifespan-and-smoothing}}),x\big)-f(\tilde t_j,p_j)}{r_j}.
  \end{equation}
  These are again solutions of \eqref{eq:Dirichlet-Willmore-flow}, with translated and rescaled boundary data $\Gamma_j$ and $\mathbf N_j$, and
  \begin{equation}\label{eq:curvature-normalization-blow-ups}
    |\mathbf A_{f_j(\tau_{\ref{thm:lifespan-and-smoothing}})}(p_j)|=\delta=\max_{0\leq t\leq\tau_{\ref{thm:lifespan-and-smoothing}}}\Vert \mathbf A_{f_j(t)}\Vert_{L^\infty(\Sigma)}.
  \end{equation}
  Since the original boundary data are fixed and smooth, we have
  \begin{align}
    \mathfrak r_{\Gamma_j}&=r_j^{-1}\mathfrak r_\Gamma\longrightarrow\infty,\quad [\bm\tau_{\Gamma_j}]_{C^\alpha(\Gamma_j)}=r_j^\alpha[\bm\tau_\Gamma]_{C^\alpha(\Gamma)}\longrightarrow0,\text{ and }\\
    \Vert D_s^m\bm\kappa_{\Gamma_j}\Vert_{C^0(\Gamma_j)}&+\Vert D_s^{m+1}\mathbf N_j\Vert_{C^0(\Gamma_j)}=r_j^{m+1}\big(\Vert D_s^m\bm\kappa_\Gamma\Vert_{C^0(\Gamma)}+\Vert D_s^{m+1}\mathbf N\Vert_{C^0(\Gamma)}\big)\longrightarrow0\label{eq:boundary-scaling-blow-ups}
  \end{align}
  for all $m\in\N_0$. Thus, for $j$ sufficiently large, \Cref{rem:Af-Linfty-Lipschitz-est} and \eqref{eq:choice-delta-blow-ups} yield that $[\mathbf n_{f_j(0)}]_{\alpha,1}\leq2C_{\ref{rem:Af-Linfty-Lipschitz-est}}^{2-\alpha}\delta\leq\varepsilon_{\ref{thm:lifespan-and-smoothing}}$. Applying \Cref{thm:lifespan-and-smoothing} to $f_j(0)$ and using \Cref{prop:geometric-uniqueness-Dirichlet}, we obtain
  \begin{equation}\label{eq:curvature-smoothing-blow-ups}
    \sup_{\frac12\tau_{\ref{thm:lifespan-and-smoothing}}\leq t\leq\tau_{\ref{thm:lifespan-and-smoothing}}}\Vert\nabla_{g_{f_j(t)}}^m \mathbf A_{f_j(t)}\Vert_{L^\infty(\Sigma)}\leq C_m\quad\text{for all $m\in\N_0$}.
  \end{equation}
  More precisely, for each fixed $m$, choose $k\geq\max\{5,m+3\}$. By \eqref{eq:boundary-scaling-blow-ups}, \eqref{eq:boundary-data-small} holds for all sufficiently large $j$, with the same $K=1$ and $\varepsilon=\varepsilon_{\ref{thm:lifespan-and-smoothing}}$. The remaining finitely many terms can be included in $C_m$. In particular, the time interval in \eqref{eq:curvature-smoothing-blow-ups} is independent of $m$.

  By \Cref{lem:area-bound-Dirichlet}, \eqref{eq:energy-dissipation-Dirichlet} and \eqref{eq:total-curvature-Dirichlet}, we also have
  \begin{align}
    \int_{f_j(t)^{-1}(B_R(z))}\dd\vol_{g_{f_j(t)}}&\leq\big(\W(f(0))+C(\Gamma,\mathbf N)\big)R^2\quad\text{for all $R>0$, $z\in\R^3$},\\
    \int_\Sigma|A_{f_j(t)}|^2\dd\vol_{g_{f_j(t)}}&\leq C\quad\text{for all $0\leq t\leq\tau_{\ref{thm:lifespan-and-smoothing}}$}.\label{eq:area-and-total-curvature-blow-ups}
  \end{align}
  Using \cite[Equation~(29)]{kuwertschaetzle2001}, and \eqref{eq:curvature-smoothing-blow-ups}, we obtain, on the time interval in \eqref{eq:curvature-smoothing-blow-ups},
  \begin{align}
    \big|\partial_t|\mathbf A_{f_j}|^2\big|&\leq C\big(|\mathbf A_{f_j}|\,|\nabla^4 \mathbf A_{f_j}|+|\mathbf A_{f_j}|^3|\nabla^2 \mathbf A_{f_j}|+|\mathbf A_{f_j}|^2|\nabla \mathbf A_{f_j}|^2+|\mathbf A_{f_j}|^6\big)\leq C.
  \end{align}
  Taking without loss of generality $C>1$, set $\hat\tau\vcentcolon=\tau_{\ref{thm:lifespan-and-smoothing}}-\min\{\tau_{\ref{thm:lifespan-and-smoothing}}/2,\delta^2/(2C)\}$.
  Then \eqref{eq:curvature-normalization-blow-ups} implies
  \begin{equation}\label{eq:nontrivial-curvature-late-interval}
    |\mathbf A_{f_j(t)}(p_j)|^2\geq\frac{\delta^2}2\quad\text{for all $t\in[\hat\tau,\tau_{\ref{thm:lifespan-and-smoothing}}]$}.
  \end{equation}
  Choose $\tau_j\in[\hat\tau,\tau_{\ref{thm:lifespan-and-smoothing}}]$ such that
  \begin{align}
    &(\tau_{\ref{thm:lifespan-and-smoothing}}-\hat\tau)\int_\Sigma|\nabla\W(f_j(\tau_j))|^2\dd\vol_{g_{f_j(\tau_j)}} \leq\int_{\hat\tau}^{\tau_{\ref{thm:lifespan-and-smoothing}}}\int_\Sigma|\nabla\W(f_j(t))|^2\dd\vol_{g_{f_j(t)}}\dd t\\
    &\quad=\W\big(f(\tilde t_j+r_j^4(\hat\tau-\tau_{\ref{thm:lifespan-and-smoothing}}))\big)-\W(f(\tilde t_j))\longrightarrow0.\label{eq:almost-stationary-slices-blow-ups}
  \end{align}
  The last convergence follows since $\tilde t_j\to T$ and $\tilde t_j+r_j^4(\hat\tau-\tau_{\ref{thm:lifespan-and-smoothing}})\to T$, and since $\lim_{t\nearrow T}\W(f(t))=\inf_{t\in(0,T)}\W(f(t))$ exists. Put
  \begin{equation}
    t_j\vcentcolon=\tilde t_j+r_j^4(\tau_j-\tau_{\ref{thm:lifespan-and-smoothing}})\longrightarrow T.
  \end{equation}
  After passing to a subsequence, $t_j\nearrow T$. We now replace $f_j(t)$ by $f_j(t)-f_j(\tau_j,p_j)$, without changing the notation. This is a time-independent translation, preserves all preceding estimates, and implies
  \begin{equation}\label{eq:def-fj}
    f_j(\tau_j)=\frac{f(t_j)-f(t_j,p_j)}{r_j},\quad\text{and thus}\quad f_j(\tau_j,p_j)=0.
  \end{equation}

  We can therefore apply \Cref{prop:compactness-with-boundary} to $(f_j(\tau_j),p_j)$. This implies \eqref{eq:pointed-convergence-blow-ups} for a smooth proper immersion of a complete connected oriented surface. By \eqref{eq:nontrivial-curvature-late-interval} and local smooth convergence,
  \begin{equation}
    |\mathbf A_{\hat f}(\hat p)|\geq\frac{\delta}{\sqrt2},\quad \Vert \mathbf A_{\hat f}\Vert_{L^\infty(\hat\Sigma)}\leq\delta,\quad \int_{\hat\Sigma}|\mathbf A_{\hat f}|^2\dd\vol_{g_{\hat f}}\leq C.
  \end{equation}
  For each fixed $\ell\in\N$, using the sets $\hat\Omega_\ell$ in \Cref{def:local-smooth-pointed-convergence} and \eqref{eq:almost-stationary-slices-blow-ups}, we have
  \begin{equation}
    \int_{\hat\Omega_\ell}|\nabla\W(\hat f)|^2\dd\vol_{g_{\hat f}}\leq\liminf_{j\to\infty}\int_\Sigma|\nabla\W(f_j(\tau_j))|^2\dd\vol_{g_{f_j(\tau_j)}}=0.
  \end{equation}
  Thus $\hat f$ is Willmore. Moreover, $\mathrm{genus}(\hat\Sigma)=\sup_{\ell\in\N}\mathrm{genus}(\hat\Omega_\ell) = \sup_{\ell\in\N}\mathrm{genus}(\Omega_\ell)\leq\mathrm{genus}(\Sigma)$ with $\Omega_\ell=\hat\phi_\ell(\hat\Omega_\ell)\subset\Sigma$ as in \Cref{def:local-smooth-pointed-convergence}.

  If $\partial\hat\Sigma\neq\emptyset$, using \eqref{eq:boundary-scaling-blow-ups} one gets that its image is contained in a single affine line and that its conormal is constant. The restriction $\hat f\vert_{\partial\hat\Sigma}$ is a proper embedding and its image is both open and closed in this line. Thus, up to applying a rotation and reflection, we have $\hat f(\partial\hat\Sigma)=a_0+\R e_1$ for some $a_0\in\R^3$ and $\partial_{\eta_{\hat f}}\hat f\equiv -e_2$. In both cases, \Cref{rem:Af-Linfty-Lipschitz-est} with $r=1$ and \eqref{eq:choice-delta-blow-ups} imply $[\mathbf n_{\hat f}]_{\alpha,1}<\infty$. This proves \eqref{eq:properties-stationary-blow-up}.

  It remains to prove \eqref{eq:density-bound-stationary-blow-up}, for which we follow the argument used in the case without boundary in \cite[Equation~(A.23)]{kuwertschaetzle2004}. For the case $\partial\hat\Sigma\neq\emptyset$ in \eqref{eq:density-bound-stationary-blow-up}, the boundary terms in \Cref{lem:monotonicity-with-boundary} must be treated carefully. By \eqref{eq:area-and-total-curvature-blow-ups} and local smooth convergence, $\hat f$ has quadratic area growth and finite Willmore energy. Since its boundary is either empty or an affine line with constant conormal, \eqref{eq:integrable-boundary-potential} holds for every $z\in\R^3$. Hence \Cref{lem:monotonicity-with-boundary} implies that the density at infinity exists.

  Take any $z_0\notin \hat f(\hat\Sigma)$. Using \Cref{rem:lsc-monotonicity-formula}, \eqref{eq:monotonicity-with-boundary} and \eqref{eq:limits-monotonicity-with-boundary} with $\theta^2(\hat f,z_0)=0$ and $\theta^2(f(t_j),\infty)=0$, we get
  \begin{align}
    \pi\theta^2(\hat f,\infty)+\frac{\W(\hat f)}{4}&+\frac12\int_{\partial\hat\Sigma}\frac{\langle\partial_{\eta_{\hat f}}\hat f,\hat f-z_0\rangle}{|\hat f-z_0|^2}\dd s\leq \sup_{R>0} \liminf_{j\to\infty} \mathcal M_{f_j(\tau_j)}(R,z_0) \\
    &\leq \liminf_{j\to\infty}\sup_{R>0}\mathcal M_{f(t_j)}(r_j R,f(t_j,p_j)+r_j z_0)\label{eq:application-of-lsc}\\
    &\leq\frac {\W(f_0)}4+\frac12\sup_{z\in\R^3}\int_\Gamma\frac{\langle\mathbf N(p),p-z\rangle}{|p-z|^2}\dd\mathcal H^1(p),
  \end{align}
  where the second estimate also uses \eqref{eq:def-fj} and the behavior of $\mathcal M$ with respect to scaling and translation. This already shows \eqref{eq:density-bound-stationary-blow-up} in the case where $\partial\hat\Sigma=\emptyset$. For $\partial\hat\Sigma\neq\emptyset$, using \eqref{eq:monotonicity-with-boundary} and \eqref{eq:limits-monotonicity-with-boundary}, we estimate for fixed $j\in\N$  
  \begin{equation}\label{eq:application-of-lsc-2}
    \sup_{R>0}\mathcal M_{f(t_j)}(r_j R,f(t_j,p_j)+r_j z_0) \leq \frac{\W(f(t_j))}{4} + \frac12\mathbf B_{\Gamma,\mathbf N}(f(t_j,p_j)+r_jz_0).
  \end{equation}
  We will argue that there exists $q\in\Gamma$ such that
  \begin{equation}\label{eq:bdry-jump-application}
    \limsup_{j\to\infty} \frac12\mathbf B_{\Gamma,\mathbf N}(f(t_j,p_j)+r_jz_0) \leq \frac{1}{2} \mathbf B_{\Gamma,\mathbf N}(q) + \frac{\pi}{2} \frac{\langle e_2,z_0-a_0\rangle}{\operatorname{dist}(z_0,a_0+\R e_1)},
  \end{equation}
  using the normalization $\hat f(\partial\hat\Sigma)=a_0+\R e_1$ and $\partial_{\eta_{\hat f}}\hat f=-e_2$ from above. Once this is established, using
  \begin{equation}
    \frac12\int_{\partial\hat\Sigma}\frac{\langle\partial_{\eta_{\hat f}}\hat f,\hat f-z_0\rangle}{|\hat f-z_0|^2}\dd s = \frac12\int_{\R}\frac{\langle e_2,  z_0-a_0\rangle}{|(s,0,0)+a_0-z_0|^2}\dd s = \frac{\pi}{2}\frac{\langle e_2,z_0-a_0\rangle}{\operatorname{dist}(z_0,a_0+\R e_1)},
  \end{equation}
  we have that \eqref{eq:density-bound-stationary-blow-up} follows by combining \eqref{eq:application-of-lsc}, \eqref{eq:application-of-lsc-2} and \eqref{eq:bdry-jump-application}. 
  
  To prove \eqref{eq:bdry-jump-application}, fix $\hat x\in\partial\hat\Sigma$ with $\hat f(\hat x)=a_0$, and set $b_j\vcentcolon=\gamma(\hat\phi_j(\hat x))$. By \eqref{eq:def-fj} and \eqref{eq:definition-pointed-convergence}, after passing to a subsequence, we get
  \begin{equation}\label{eq:asymptotics-bj-f(t_j,p_j)}
    b_j\in\Gamma,\quad \frac{b_j-f(t_j,p_j)}{r_j}\longrightarrow a_0,\quad b_j\longrightarrow q\in\Gamma,\quad\text{and}\quad f(t_j,p_j)\longrightarrow q.
  \end{equation}
  Moreover, \eqref{eq:Dirichlet-Willmore-flow} and \eqref{eq:definition-pointed-convergence} yield that
  \begin{align}
    \mathbf N(q)&=\lim_{j\to\infty}\mathbf N(b_j)=\lim_{j\to\infty}\partial_{\eta_{f_j(\tau_j)}}f_j(\tau_j,\hat\phi_j(\hat x))=\partial_{\eta_{\hat f}}\hat f(\hat x)=-e_2,\text{ and}\\
    T_q\Gamma&=\lim_{j\to\infty}T_{b_j}\Gamma=d\hat f(T_{\hat x}\partial\hat\Sigma)=\R e_1.
  \end{align}

  Let $\pi_\Gamma$ be the nearest point projection on a tubular neighborhood of $\Gamma$, and set $z_j\vcentcolon=f(t_j,p_j)+r_jz_0$. Then
  \begin{equation}
    D\pi_\Gamma(b_j)=P_{T_{b_j}\Gamma}\longrightarrow P_{\R e_1}\quad\text{and}\quad
    |\pi_\Gamma(z_j)-b_j-D\pi_\Gamma(b_j)(z_j-b_j)|\leq C|z_j-b_j|^2\leq Cr_j^2,
  \end{equation}
  and hence, using \eqref{eq:asymptotics-bj-f(t_j,p_j)},
  \begin{equation}
    \frac{z_j-\pi_\Gamma(z_j)}{r_j}
    =(I-D\pi_\Gamma(b_j))\frac{f(t_j,p_j)+r_jz_0-b_j}{r_j}+O(r_j)
    \longrightarrow (I-P_{\R e_1})(z_0-a_0)=\vcentcolon w,
  \end{equation}
  and in particular $|w|=\operatorname{dist}(z_0,a_0+\R e_1)>0$. Thus, for all sufficiently large $j$,
  \begin{align}
    \pi_\Gamma(z_j)&\longrightarrow q,\quad 0<d_j\vcentcolon=|z_j-\pi_\Gamma(z_j)|=r_j(|w|+o(1))\longrightarrow0,\\
    v_j&\vcentcolon=\frac{z_j-\pi_\Gamma(z_j)}{d_j}\in(T_{\pi_\Gamma(z_j)}\Gamma)^\bot,\quad |v_j|=1,\quad
    v_j\longrightarrow\frac{w}{|w|}\in(T_q\Gamma)^\bot.
  \end{align}
  Applying \Cref{lem:boundary-integral-jump} with $z_j=\pi_\Gamma(z_j)+d_j v_j$, we obtain
  \begin{equation}
    \lim_{j\to\infty}\frac12\mathbf B_{\Gamma,\mathbf N}(z_j)
    =\frac12\mathbf B_{\Gamma,\mathbf N}(q)-\frac{\pi}{2}\left\langle\mathbf N(q),\frac{w}{|w|}\right\rangle
    =\frac12\mathbf B_{\Gamma,\mathbf N}(q)+\frac{\pi}{2}\frac{\langle e_2,z_0-a_0\rangle}{\operatorname{dist}(z_0,a_0+\R e_1)}.
  \end{equation}
  This proves \eqref{eq:bdry-jump-application}.
\end{proof}

\begin{remark}[Subconvergence in case (a)]\label{rem:subconvergence-in-case-a}
  In the setting of \Cref{prop:existence-of-blow-ups}, suppose that $T=\infty$ and $M=\max\{1,\sup_{t\geq 0}\Vert\mathbf A_{f(t)}\Vert_{L^\infty(\Sigma)}\}<\infty$. Proceeding as in the proof of \Cref{prop:existence-of-blow-ups}(b), we get the following subconvergence result.
  
  Let $p_0\in\partial\Sigma$ and $\tilde t_j\nearrow \infty$. With $\delta=\delta(\Gamma,\mathbf N)$ sufficiently small, for $r_j\equiv r=\delta/M$,
  \begin{equation}
    f_j(t,x)=\frac{f(\tilde t_j+r_j^4(t-\tau_{\ref{thm:lifespan-and-smoothing}}),x)-\gamma(p_0)}{r_j}.
  \end{equation}
  \Cref{rem:Af-Linfty-Lipschitz-est,thm:lifespan-and-smoothing,prop:geometric-uniqueness-Dirichlet} yield suitable analogs of \eqref{eq:curvature-smoothing-blow-ups} and \eqref{eq:area-and-curvature-compactness}. Moreover,  similarly as in \eqref{eq:almost-stationary-slices-blow-ups}, we can choose $\tau_j\in[\tau_{\ref{thm:lifespan-and-smoothing}}/2,\tau_{\ref{thm:lifespan-and-smoothing}}]$ such that
  \begin{align}
    &\frac12\tau_{\ref{thm:lifespan-and-smoothing}}\int_\Sigma|\nabla\W(f_j(\tau_j))|^2\dd\vol_{g_{f_j(\tau_j)}} \leq\int_{\frac12\tau_{\ref{thm:lifespan-and-smoothing}}}^{\tau_{\ref{thm:lifespan-and-smoothing}}}\int_\Sigma|\nabla\W(f_j(t))|^2\dd\vol_{g_{f_j(t)}}\dd t\\
    &\quad=\W\big(f(\tilde t_j-\frac12 r_j^4\tau_{\ref{thm:lifespan-and-smoothing}})\big)-\W(f(\tilde t_j))\longrightarrow0.\label{eq:point-picking-2}
  \end{align}
  Since $p_0\in\partial\Sigma$, the boundary conditions in \eqref{eq:Dirichlet-Willmore-flow} yield that $f_j(\tau_j,p_0)=0$, and we can apply \Cref{prop:compactness-with-boundary} to $(f_j(\tau_j),p_0)$ to get \eqref{eq:pointed-convergence-blow-ups} after passing to a subsequence, for a smooth proper immersion of a complete connected oriented surface $\hat\Sigma$ with $\hat p\in\partial\hat\Sigma$, and $\mathrm{genus}(\hat\Sigma)\leq\mathrm{genus}(\Sigma)$. By the choice $r=\delta/M$, \eqref{eq:total-curvature-Dirichlet} and local smooth convergence,
  \begin{equation}
    \Vert \mathbf A_{\hat f}\Vert_{L^\infty(\hat\Sigma)}\leq\delta,\quad\text{and}\quad \int_{\hat\Sigma}|\mathbf A_{\hat f}|^2\dd\vol_{g_{\hat f}}\leq C,
  \end{equation}
  and $\hat f$ is Willmore by \eqref{eq:point-picking-2}, i.e.\ satisfies $\nabla\W(\hat f)=0$ on $\hat\Sigma$. Finally, arguing as in \cite[Lemma~4.3]{kuwertschaetzle2001}, one obtains that either $\hat\Sigma$ is noncompact, or one may take $\hat\Sigma=\Sigma$, and in the latter case, up to passing to a subsequence, that there are diffeomorphisms $\Phi_j\colon\Sigma\to\Sigma$ such that $f_j(\tau_j)\circ\Phi_j\to\hat f$ in $C^\infty(\Sigma)$. Moreover, in both cases 
  \begin{equation}
    \hat f\vert_{\partial\hat\Sigma}\colon \partial\hat\Sigma\to \frac{1}{\delta/M}(\Gamma-\gamma(p_0))
  \end{equation}
  is an embedding which is surjective if $\hat\Sigma=\Sigma$ is compact. If $\hat\Sigma$ is noncompact, $\hat f(\partial\hat\Sigma)$ contains at least the (translated and rescaled) component of $\Gamma$ containing $\gamma(p_0)$.
\end{remark}

\subsection{Rigidity: A gap result for Dirichlet Willmore blow-ups}

This section's main result \Cref{lem:rigidity-bdry-blow-ups} will yield the blow-up classification \Cref{prop:intro-boundary-rigidity} in the introduction. The underlying idea is an adaptation of the Pohozaev-type identity used in \cite{palmer2000,dallacqua2012}. In order to apply it to proper immersions of \emph{noncompact} surfaces, we first prove an $\varepsilon$-boundary-regularity estimate in the spirit of \cite[Theorem~2.10]{kuwertschaetzle2001} for Willmore immersions whose boundary is a line, and with constant conormal field. 

For $\varepsilon$-regularity estimates, the proof strategy employed in \cite{kuwertschaetzle2001} does not immediately translate to boundary estimates for Willmore surfaces, as the underlying energy estimates would involve boundary terms which cannot simply be treated in the same way as the surface integrals. We thus localize in order to apply suitable elliptic regularity estimates on the (local) graph representations. To guarantee existence of such local graph representations with uniform radius and $C^{1}$-bounds, we impose an additional smallness assumption on the unit normal $\mathbf n_f$ which is motivated by the following remark.

\begin{remark}\label{rem:flat-graph-rep}
  Let $\hat\Sigma$ be a complete oriented surface with $\partial\hat\Sigma\neq\emptyset$, and $f\colon\hat\Sigma\to\R^3$ a smooth immersion satisfying $f(\partial\hat\Sigma) = \R e_1$ and $\partial_{\eta_{f}}f = (0,-1,0)$ on $\partial\hat\Sigma$.
  Using \Cref{rem:Af-Linfty-Lipschitz-est,prop:ex-of-graph-charts} and setting $\varepsilon_{\ref{rem:flat-graph-rep}}=\min\{\varepsilon_{\ref{rem:Af-Linfty-Lipschitz-est}}, \varepsilon_{\ref{prop:ex-of-graph-charts}}/C_{\ref{rem:Af-Linfty-Lipschitz-est}}^2\}$, one can show the following. Suppose that, for $0<\varepsilon\leq\varepsilon_{\ref{rem:flat-graph-rep}}$, $q\in\hat\Sigma$ and $r>0$, we have either
  \begin{equation}\label{eq:1-rem:flat-graph-rep}
    (4r)^\alpha[\mathbf n_f]_{\alpha,4r,q}\leq\varepsilon\quad\text{or}\quad 4r\Vert\mathbf A_f\Vert_{L^\infty(f^{-1}(B_{4r}(f(q))))}\leq\varepsilon \leq\varepsilon_{\ref{rem:flat-graph-rep}}
  \end{equation}
  If $q\in\partial\hat\Sigma$, or $d_{g_f}(q,\partial\hat\Sigma)\geq2\sqrt{2}r$, we have a local graph representation as in \eqref{eq:2-prop:ex-of-graph-charts} over $Q_r^+$, respectively $Q_r$, satisfying \eqref{eq:3-prop:ex-of-graph-charts} and, in the second case in \eqref{eq:1-rem:flat-graph-rep}, also
  \begin{equation}
    \frac{\Vert u_q\Vert_{L^\infty(D_{r,q})}}{r}+\Vert\nabla u_q\Vert_{L^\infty(D_{r,q})}+r\Vert D^2u_q\Vert_{L^\infty(D_{r,q})}\leq C_{\ref{rem:flat-graph-rep}}\varepsilon.
  \end{equation}
  where $C_{\ref{rem:flat-graph-rep}}>0$ is a universal constant.
\end{remark}

\begin{proposition}[$\varepsilon$-boundary-regularity for flat-clamped Willmore surfaces]
  \label{prop:eps-bdry-reg}
  Let $\hat\Sigma$ be a complete oriented surface with $\partial\hat\Sigma\neq\emptyset$, $f\colon\hat\Sigma\to\R^3$ a smooth immersion with $[\mathbf n_f]_{\alpha,1} < \varepsilon_{\ref{rem:flat-graph-rep}}$, satisfying 
  \begin{equation}
    \begin{cases}
      \Delta_{g_f}\mathbf H_{f}+|\mathbf A^0_f|^2\mathbf H_f=0&\text{on $\hat\Sigma$}\\
      f(\partial\hat\Sigma) = \R e_1\\
      \partial_{\eta_{f}}f = (0,-1,0)&\text{on $\partial\hat\Sigma$}.
    \end{cases}
  \end{equation}
  There exist universal constants $\sigma_0>0$ and $(C_k)_{k\geq 0}$ such that, if for some $p\in\R e_1$ and $R>8$
  \begin{equation}
    \int_{f^{-1}(B_R(p))} |\mathbf A_f|^2\dd\vol_{g_f} \leq\sigma^2<\sigma_0^2,
  \end{equation}
  then 
  \begin{equation}\label{eq:1-prop:eps-bdry-reg}
    \Vert \nabla^k_{g_f} \mathbf A_f \Vert_{L^\infty(f^{-1}(B_{\frac12 R}(p)))} \leq \frac{C_k}{R^{k+1}} \sigma\quad\text{for all $k\geq 0$}.
  \end{equation}
\end{proposition}

Our proof of \Cref{prop:eps-bdry-reg} relies on the local graph representations in \Cref{rem:flat-graph-rep}, and elliptic regularity theory applied to the graph Willmore operator, also see \cite{schaetzle2010,gulyak2026}. Revisiting the proof, we only need that $f$ can be locally represented as a graph over uniform rectangles $Q_r$ (resp.\ $Q_r^+$) in the tangent plane with uniform Lipschitz control --- so the $[\mathbf n_f]_{\alpha,1}$ assumption could be replaced by a weaker assumption, such as suitable control of the modulus of continuity of the Gauss map, cf.\ \cite[Lemma~V.10, Corollary~V.11]{lanmartinoriviere2026}. However, the statement of \Cref{prop:eps-bdry-reg} is sufficient for our purposes.

\begin{lemma}\label{lem:graphical-eps-bdry-regularity}
  There exists a universal constant $\sigma_0>0$ with the following property. Let $\mathfrak r>0$, $u\in W^{4,2}(Q_\mathfrak r^+)$ satisfy $\Delta_{g_{f_u}}\mathbf H_{f_u}+|\mathbf A^0_{f_u}|^2\mathbf H_{f_u}=0$ in $Q_\mathfrak r^+$, $u=\nabla u=0$ on $Q_\mathfrak r^+\cap (\R\times\{0\})$, and
  \begin{equation}\label{eq:2-lem:graphical-eps-bdry-regularity}
    |\nabla u|\leq C_{\ref{rem:flat-graph-rep}}\varepsilon_{\ref{rem:flat-graph-rep}}
    ,\quad \int_{Q_\mathfrak r^+} |\mathbf A_{f_u}|^2\dd\vol_{g_{f_u}} \leq \sigma^2<\sigma_0^2,
  \end{equation}
  where $f_u(x)=(x,u(x))$ denotes the associated graph. Then we have $u\in C^\infty(Q_{\frac23\mathfrak r}^+)$ and, for all $k\geq 2$, there are universal constants $C_k>0$ only depending on $k$ such that
  \begin{equation}\label{eq:1-lem:graphical-eps-bdry-regularity}
    \Vert \nabla^{k-2}_{g_{f_u}}\mathbf A_{f_u}\Vert_{L^\infty(Q_{\frac23\mathfrak r}^+)}+\Vert \nabla^k u\Vert_{L^\infty(Q_{\frac23 \mathfrak r}^+)} \leq \frac{C_k}{\mathfrak r^{k-1}}\sigma.
  \end{equation}
\end{lemma}
\begin{remark}\label{rem:graphical-eps-bdry-regularity}
    Notice that in \eqref{eq:2-lem:graphical-eps-bdry-regularity} we only require an upper bound for $\vert\nabla u\vert$, and no $\mathfrak r^\alpha$-decay for higher-order Hölder seminorms, by relying on the elliptic regularity estimates \cite{dongkim2011} for operators with BMO-small-coefficients instead of the classical theory \cite{agmondouglisnirenberg1959}. The key observation is that the smallness of the second fundamental form in $L^2$ bounds the squared BMO-seminorm for the coefficients of the elliptic Willmore graph-operator. And thus one obtains uniform elliptic regularity estimates for the homogeneous boundary conditions $u=\nabla u=0$ not depending on any a-priori decay of higher order derivatives. Also compare the related graphical $\varepsilon$-boundary-regularity estimate in \cite[Proposition~3.1]{schaetzle2010} for \emph{minimizers} where $u$ is assumed to be controlled in $C^{1+\alpha}\cap W^{2,2}$.
\end{remark}
\begin{proof}[Proof of \Cref{lem:graphical-eps-bdry-regularity}]
  Replacing $u$ by $\mathfrak r^{-1}u(\mathfrak r\cdot)$ and using the scaling invariance of the Willmore equation, of the assumptions, and of the claims, it suffices to prove the lemma for $\mathfrak r=1$. Following \cite[Section~2]{gulyak2026}, there exist smooth functions $A^4_\alpha\colon\R^2\to\R$ for $|\alpha|=4$, $b^{2,3}_{\alpha,\beta}\colon\R^2\to\R$ for $|\alpha|=2$, $|\beta|=3$, and $b^{2,2,2}_{\alpha,\beta,\gamma}\colon\R^2\to\R$ for $|\alpha|=|\beta|=|\gamma|=2$ such that $u$ satisfies
  \begin{equation}\label{eq:1-pf-lem:graphical-eps-bdry-regularity}
    \begin{cases}
      A^4_\alpha(\nabla u)\partial^\alpha u = b^{2,3}_{\alpha,\beta}(\nabla u) \partial^\alpha u \partial^\beta u + b^{2,2,2}_{\alpha,\beta,\gamma}(\nabla u) \partial^\alpha u\partial^\beta u\partial^\gamma u &\text{on $Q_1^+$}\\
      u=\nabla u=0&\text{on $(-1,1)\times\{0\}$}
    \end{cases}
  \end{equation}
  where, for $q,\xi=(\xi^1,\xi^2)\in\R^2$,
  \begin{equation}\label{eq:2-pf-lem:graphical-eps-bdry-regularity}
    A^4_{(i,j,k,l)}(q) \xi^i\xi^j\xi^k\xi^l\geq \frac{|\xi|^4}{(1+|q|^2)^\frac{5}{2}}.
  \end{equation}
  Moreover, letting $E_r(x)=(x+Q_r)\cap Q_1^+$, an elementary estimate involving the Poincar\'e and Hölder inequalities and \eqref{eq:2-lem:graphical-eps-bdry-regularity} yields for a universal constant $C>0$ that
  \begin{equation}\label{eq:3-pf-lem:graphical-eps-bdry-regularity}
    \sup_{\substack{x\in Q_1^+,r>0\\\text{s.th. } x+Q_r\subset Q_1}} \fint_{E_r(x)}\Big\vert A^4_\alpha(\nabla u(y)) - \fint_{E_r(x)} A^4_\alpha(\nabla u(z))\dd z\Big\vert\dd y \leq C \Vert \nabla A^4_{\alpha}\Vert_{L^\infty(\bar B_{C\varepsilon_{\ref{rem:flat-graph-rep}}}(0))}\sigma_0.
  \end{equation}
  Let $0<r<R\leq 7/8$. By \eqref{eq:2-lem:graphical-eps-bdry-regularity}, \eqref{eq:2-pf-lem:graphical-eps-bdry-regularity} and \eqref{eq:3-pf-lem:graphical-eps-bdry-regularity}, choosing $\sigma_0$ sufficiently small only depending on universal constants, \cite[Theorem~8]{dongkim2011}, using the equivalent formulation of BMO with squares, can be applied to $u$ multiplied with a suitable cut-off function, on a suitable universally chosen smooth domain between $Q_{7/8}^+$ and $Q_1^+$ to obtain that, for a universal constant $C>0$,
  \begin{align}
    \Vert u\Vert_{W^{4,2}(Q_r^+)} &\leq C\big(\Vert A^4_\alpha(\nabla u)\partial^\alpha u\Vert_{L^2(Q_R^+)} + (R-r)^{-4}\Vert u\Vert_{L^2(Q_R^+)}\big)\\
    &\leq C\big(\Vert A^4_\alpha(\nabla u)\partial^\alpha u\Vert_{L^2(Q_R^+)} + (R-r)^{-4}\sigma \big) \label{eq:4-pf-lem:graphical-eps-bdry-regularity}
  \end{align}
  using the homogeneous boundary conditions and Poincar\'e in the last step. For the $L^2(Q_R^+)$-norm of the right-hand-side of \eqref{eq:1-pf-lem:graphical-eps-bdry-regularity}, using the gradient bound in \eqref{eq:2-lem:graphical-eps-bdry-regularity}, we estimate
  \begin{align}
    \Vert b^{2,3}_{\alpha,\beta}(\nabla u) \partial^\alpha u \partial^\beta u\Vert_{L^2(Q_R^+)} &\leq C \Vert D^2 u\Vert_{L^6(Q_R^+)}\Vert D^3 u\Vert_{L^3(Q_R^+)},\\
    \Vert b^{2,2,2}_{\alpha,\beta,\gamma}(\nabla u) \partial^\alpha u\partial^\beta u\partial^\gamma u\Vert_{L^2(Q_R^+)} &\leq C \Vert D^2 u\Vert_{L^6(Q_R^+)}^3.
  \end{align}
  Moreover, by Gagliardo--Nirenberg and \eqref{eq:2-lem:graphical-eps-bdry-regularity}, we obtain for a universal constant $C>0$
  \begin{align}
    \Vert D^2 u\Vert_{L^6(Q_R^+)} &\leq C(\sigma^{\frac23}\Vert u\Vert_{W^{4,2}(Q_R^+)}^\frac13 + R^{-\frac23}\sigma) ,\quad \Vert D^3 u\Vert_{L^3(Q_R^+)} \leq C(\sigma^{\frac13}\Vert u\Vert_{W^{4,2}(Q_R^+)}^\frac23 + R^{-\frac43}\sigma).
  \end{align}
  Altogether, \eqref{eq:4-pf-lem:graphical-eps-bdry-regularity} yields that 
  \begin{equation}
    \Vert u\Vert_{W^{4,2}(Q_r^+)}\leq C\sigma \Vert u\Vert_{W^{4,2}(Q_R^+)} + C(R-r)^{-4}\sigma.
  \end{equation}
  Choosing $\sigma_0<\frac{1}{2C}$, a standard iteration argument yields that for a universal constant $C>0$
  \begin{equation}
    \Vert u\Vert_{W^{4,2}(Q_{\frac34}^+)} \leq C\sigma.
  \end{equation}
  In particular, using $W^{4,2}(Q_{\frac34}^+)\hookrightarrow C^{2+\alpha}(Q_{\frac34}^+)$, one can proceed to prove the higher-order estimates on $Q_{\frac23}^+$ with a similar strategy as in \cite{gulyak2026}, using \cite[Section~9]{agmondouglisnirenberg1959}.
\end{proof}

\begin{proof}[Proof of \Cref{prop:eps-bdry-reg}]
  Without loss of generality, assume $\sigma>0$. We first argue that, for $0<r<s\leq R$, we have the implication
  \begin{equation}\label{eq:local-cover-eps-bdry-reg}
    r\Vert \mathbf A_f\Vert_{L^\infty(f^{-1}(B_s(p)))}\leq\varepsilon_{\ref{rem:flat-graph-rep}} \implies \Vert \mathbf A_f\Vert_{L^\infty(f^{-1}(B_{s-r}(p)))}\leq\frac{C_{\mathrm{reg}}\sigma}{r}
  \end{equation}
  where $C_{\mathrm{reg}}>0$ is a universal constant. Indeed, for any $z\in f^{-1}(B_{s-r}(p))$ with $d_{g_f}(z,\partial\hat\Sigma)<r/16$, use \Cref{rem:flat-graph-rep} at a nearest boundary point $q_z$ with graph radius $r/8$. Then $z\in U_{r/16,q_z}$, with graph coordinate in $Q_{r/16}^+\subset Q_{r/12}^+$, so \Cref{lem:graphical-eps-bdry-regularity} applies with $\mathfrak r=r/8$. If $d_{g_f}(z,\partial\hat\Sigma)\geq r/16$, apply \cite[Theorem~2.10 and Remark~2.11]{kuwertschaetzle2001} on the intrinsic ball $B^{g_f}_{r/32}(z)$. All these charts and the neighbourhoods used above lie in $f^{-1}(B_s(p))$. Thus, after choosing $\sigma_0$ smaller if necessary, we get \eqref{eq:local-cover-eps-bdry-reg}.

  Moreover, arguing analogously as for \eqref{eq:local-cover-eps-bdry-reg}, for $r_0=1$, using that $[\mathbf n_f]_{\alpha,1}<\varepsilon_{\ref{rem:flat-graph-rep}}$ and the first case in \eqref{eq:1-rem:flat-graph-rep}, we get
  \begin{equation}
    \Vert \mathbf A_f\Vert_{L^\infty(f^{-1}(B_{R-1}(p)))}\leq\frac{C_{\mathrm{reg}}\sigma}{r_0}.
  \end{equation}
  Inductively suppose that, for $j\in\N$, $r_0=1<\cdots<r_{j-1}$ are chosen such that we have
  \begin{equation}
    \sum_{i=0}^{j-1}r_i < \frac18R\quad\text{and}\quad \Vert \mathbf A_f\Vert_{L^\infty(f^{-1}(B_{R-\sum_{i=0}^{j-1}r_i}(p)))} \leq \frac{C_{\mathrm{reg}}\sigma}{r_{j-1}}
  \end{equation}
  where $C_{\mathrm{reg}}$ is chosen as above. We may assume that $\sigma\leq\sigma_0<\varepsilon_{\ref{rem:flat-graph-rep}}/(2C_{\mathrm{reg}})$, so $\lambda\vcentcolon=\varepsilon_{\ref{rem:flat-graph-rep}}/(C_{\mathrm{reg}}\sigma)>2$. Then there exists $r_j>r_{j-1}$ which is maximal with respect to the property that
  \begin{equation}
    r_j\leq r_{j-1} \frac{\varepsilon_{\ref{rem:flat-graph-rep}}}{C_{\mathrm{reg}}\sigma}\quad\text{and}\quad \sum_{i=0}^{j}r_i \leq \frac14 R.
  \end{equation}
  In particular, we obtain that
  \begin{equation}
    r_j\Vert \mathbf A_f\Vert_{L^\infty(f^{-1}(B_{R-\sum_{i=0}^{j-1}r_i}(p)))} \leq r_j \frac{C_{\mathrm{reg}}\sigma}{r_{j-1}} \leq \varepsilon_{\ref{rem:flat-graph-rep}},
  \end{equation}
  and thus \eqref{eq:local-cover-eps-bdry-reg}, with $s=R-\sum_{i=0}^{j-1}r_i$ and $r=r_j$, yields
  \begin{equation}
    \Vert \mathbf A_f\Vert_{L^\infty(f^{-1}(B_{R-\sum_{i=0}^{j}r_i}(p)))} \leq \frac{C_{\mathrm{reg}}\sigma}{r_j}.
  \end{equation}
  This process terminates at the first $J$ with $\sum_{j=0}^Jr_j\geq R/8$, since every step before a possible final truncated step multiplies the radius by $\lambda>2$. Moreover, we have $\sum_{j=0}^Jr_j\leq R/4$ and $r_J\geq R/16$. Indeed, if the final step is not truncated, then $\sum_{j=0}^Jr_j\leq2r_J$, and otherwise $r_J=R/4-\sum_{j=0}^{J-1}r_j>R/8$. Thus,
  \begin{equation}
    \Vert \mathbf A_f\Vert_{L^\infty(f^{-1}(B_{\frac34R}(p)))}\leq \Vert \mathbf A_f\Vert_{L^\infty(f^{-1}(B_{R-\sum_{i=0}^Jr_i}(p)))} \leq \frac{C_{\mathrm{reg}}\sigma}{r_J} \leq \frac{16C_{\mathrm{reg}}\sigma}{R}.
  \end{equation}
  Choosing $\sigma_0$ sufficiently small such that $4C_{\mathrm{reg}}\sigma_0\leq\varepsilon_{\ref{rem:flat-graph-rep}}$, applying \eqref{eq:local-cover-eps-bdry-reg} with $s=3R/4$ and $r=R/4$ yields \eqref{eq:1-prop:eps-bdry-reg} for $k=0$. Moreover, arguing as in the proof of $\eqref{eq:local-cover-eps-bdry-reg}$ and using \Cref{lem:graphical-eps-bdry-regularity} and \cite[Theorem~3.5]{kuwertschaetzle2001}, we also get \eqref{eq:1-prop:eps-bdry-reg} for all $k\geq 1$.
\end{proof}

With the $\varepsilon$-boundary-regularity result, we can adapt an observation of \cite{dallacqua2012} to prove the main result of this section which in particular implies \Cref{prop:intro-boundary-rigidity}, i.e.\ that our blow-ups of the Dirichlet Willmore flow satisfy the additional boundary condition $\mathbf H_f=0$. 

\begin{proposition}[Rigidity for boundary-blow-ups] \label{lem:rigidity-bdry-blow-ups}
  Let $\hat\Sigma$ be a complete, oriented, connected surface with $\partial\hat\Sigma\neq\emptyset$, and $f\colon\hat\Sigma\to\R^3$ a smooth proper immersion with finite density at infinity, satisfying 
  \begin{equation}\label{eq:1-lem:rigidity-of-blow-ups}
    \begin{cases}
      \Delta_{g_f}\mathbf H_{f}+|\mathbf A^0_f|^2\mathbf H_f=0&\text{on $\hat\Sigma$}\\
      f(\partial\hat\Sigma) = \R e_1\\
      \partial_{\eta_{f}}f = (0,-1,0)&\text{on $\partial\hat\Sigma$}.
    \end{cases}
  \end{equation}
  Furthermore, assume that
  \begin{equation}
    \int_{\hat\Sigma}|\mathbf A_f|^2\dd\vol_{g_f}<\infty \quad\text{and}\quad [\mathbf n_f]_{\alpha,1}<\infty
  \end{equation} 
  for some $\alpha\in(0,1)$. Then $\mathbf H_f=0$ on $\partial\hat\Sigma$ and the reflection $\mathcal Rf\colon D(\hat\Sigma)\to\R^3$ of $f$ about $\R e_1$ is a smooth Willmore immersion of the complete oriented connected doubling\footnote{See \cite[Theorem~9.29 and Example~9.32]{lee2013}.} $D(\hat\Sigma)$. 

  Furthermore, if $\mathrm{genus}(\hat\Sigma)=0$, $f\vert_{\partial\hat\Sigma}\colon\partial\hat\Sigma\to \R e_1$ is a diffeomorphism, and the density at infinity satisfies
  \begin{equation}
    \theta^2(f,\infty)\vcentcolon=\lim_{R\to\infty} \frac{\int_{f^{-1}(B_R(0))}\dd\vol_{g_f}}{\pi R^2} < 1,
  \end{equation} 
  then $f$ parametrizes a flat half-plane.
\end{proposition}
\begin{proof}
  After a fixed dilation of $f$, without loss of generality we can assume that $[\mathbf n_f]_{\alpha,1}<\varepsilon_{\ref{rem:flat-graph-rep}}$, as required in \Cref{prop:eps-bdry-reg}. For any compact Lipschitz domain $\Omega\subset\hat\Sigma$ (not necessarily connected) with $\partial\Omega\cap\partial\hat\Sigma\neq\emptyset$, testing the Willmore equation with $f+e_2$, the scaling- and translation invariance of $\W$ yields as in \cite[Equation~(2.1)]{dallacqua2012} that
  \begin{equation}\label{eq:2-pf-lem:rigidity-of-blow-ups}
    \int_{\partial\Omega} 2\partial_{\eta_\Omega}\mathbf H_f\langle f+e_2,\mathbf n_f\rangle-\partial_{\eta_{\Omega}}(\mathbf H_f\langle f+e_2,\mathbf n_f\rangle ) - \frac12 \mathbf H_f^2\langle f+e_2,\partial_{\eta_\Omega}f\rangle \dd s_\Omega = 0
  \end{equation}
  where $\dd s_\Omega\vcentcolon= \dd \vol_{\iota_{\partial\Omega}^*g_f}$ and $\eta_{\Omega}$ denotes the outer unit conormal of $\Omega$ in $(\hat\Sigma,g_f)$. Proceeding as in \cite[Equations~(2.3)-(2.5)]{dallacqua2012}, and using that \eqref{eq:1-lem:rigidity-of-blow-ups} yields $\langle f+e_2,\mathbf n_f\rangle=0$ on $\partial\hat\Sigma$ and $\mathbf A_f(s_f,\cdot)=0$ for a $g_f$-unit boundary tangent $s_f$, and thus $\mathbf H_f= \mathbf A_f(\eta_f,\eta_f)$ on $\partial\hat\Sigma$, we obtain that 
  \begin{equation}
    \partial_{\eta_\Omega}(\mathbf H_f\langle f+e_2,\mathbf n_f\rangle ) = -\langle f+e_2,\partial_{\eta_\Omega}f\rangle \mathbf H_f^2\quad\text{ on $\partial\Omega\cap\partial\hat\Sigma$.}
  \end{equation}
  Thus, using \eqref{eq:2-pf-lem:rigidity-of-blow-ups} and again that $\langle f+e_2,\mathbf n_f\rangle=0$ on $\partial\hat\Sigma$, we have
  \begin{equation}\label{eq:3-pf-lem:rigidity-of-blow-ups}
    \int_{\partial\Omega\cap\partial\hat\Sigma} \mathbf H_f^2 \dd s_{\Omega} = 2\int_{\partial\Omega\setminus\partial\hat\Sigma} \partial_{\eta_\Omega} \mathbf H_f\langle f+e_2,\mathbf n_f\rangle + \mathbf H_f \mathbf A_f(f+e_2,\eta_\Omega) - \frac12 \mathbf H_f^2\langle f+e_2,\partial_{\eta_\Omega}f\rangle \dd s_\Omega.
  \end{equation}
  We use \Cref{prop:eps-bdry-reg} to show that the right-hand side tends to zero along a suitable sequence of domains $\Omega_j$. We first argue that there exist radii $\rho_j\in(j/2,j)$ such that, taking $\Omega_j=f^{-1}(B_{\rho_j}(0))$, $\Omega_j$ is a Lipschitz domain, and we have that
  \begin{align}
    &\qquad\int_{\partial\Omega_j\setminus\partial\hat\Sigma} 1 \dd s_{\Omega_j} \leq C j, \quad |f+e_2|\leq 2 j\text{ on $\partial\Omega_j$}, \\
    & \max_{\partial\Omega_j\setminus\partial\hat\Sigma} j^{k+1}|\nabla_{g_f}^k \mathbf A_f| = o(1) \quad\text{as $j\to\infty$, for all $k\in\N_0$}.\label{eq:1-pf-lem:rigidity-bdry-blow-ups}
  \end{align}
  This is achieved by choosing $\rho_j\in(j/2,j)$ by Sard's theorem to be a regular value of $|f|$ such that, denoting by $\mathcal H_{g_f}^1$ the one-dimensional Hausdorff measure with respect to the metric $g_f$, we have
  \begin{align}
    \mathcal H_{g_f}^1(|f|^{-1}(\rho_j))&\leq \frac{2}{j}\int_{j/2}^j \mathcal H_{g_f}^1(|f|^{-1}(t))\dd t \\
    &\leq \frac{2}{j} \int_{f^{-1}(B_j(0))} \big|\nabla_{g_f}|f|\big|\dd\vol_{g_f} \leq \frac{2}{j} \int_{f^{-1}(B_j(0))} 1\dd\vol_{g_f}
  \end{align}
  where we used the co-area formula in the second estimate. As the density at infinity is finite we get $\int_{f^{-1}(B_j(0))} 1\dd\vol_{g_f}\leq Cj^2$. So choosing $\Omega_j=f^{-1}(B_{\rho_j}(0))$, we have $\partial\Omega_j\setminus\partial\hat\Sigma\subset |f|^{-1}(\rho_j)$, and as $\rho_j$ is a regular value of $|f|$, $\Omega_j\setminus\partial\hat\Sigma$ is a smooth domain in $\hat\Sigma\setminus\partial\hat\Sigma$. Thus, the first line in \eqref{eq:1-pf-lem:rigidity-bdry-blow-ups} holds. For the second line, let $z\in\partial\Omega_j\setminus\partial\hat\Sigma$. If $\mathrm{dist}(f(z),\R e_1)\leq j/16$, apply \Cref{prop:eps-bdry-reg} on the ball of radius $j/4$ centred at the projection of $f(z)$ onto the boundary line. Otherwise use \cite[Theorem~3.5]{kuwertschaetzle2001} on $B_{j/32}(f(z))$. These balls avoid $B_{j/8}(0)$ and their concentric balls of half the radius contain $f(z)$. Thus the second line in \eqref{eq:1-pf-lem:rigidity-bdry-blow-ups} follows from 
  \begin{equation}
    \int_{f^{-1}(\R^3\setminus B_{j/8}(0))}| \mathbf A_f|^2\dd\vol_{g_f}\longrightarrow0\quad\text{as $j\to\infty$}.
  \end{equation}
  Next, by \eqref{eq:1-pf-lem:rigidity-bdry-blow-ups}, we can apply \eqref{eq:3-pf-lem:rigidity-of-blow-ups} with $\Omega=\Omega_j$. Then letting $j\to\infty$ along this sequence yields $\mathbf H_f=0$ on $\partial\hat\Sigma$. 

  Now consider the doubling $D(\hat\Sigma)$ of $\hat\Sigma$ as in \cite[Example~9.32]{lee2013}, where $D(\hat\Sigma) = \hat\Sigma\times\{1,2\} / \sim$ where $(p_1,i_1)\sim(p_2,i_2)$ iff.\ $(p_1,i_1)=(p_2,i_2)$ or $p_1=p_2\in\partial\hat\Sigma$. By an abuse of notation, we identify $\partial\hat\Sigma=\{[(p,1)]\mid p\in\partial\hat\Sigma\}$ and assume that the smooth structure on the doubling is induced by reflecting isothermal boundary charts of $f$. 
  Then the reflection of $f$ about $\R e_1$ is given by 
  \begin{equation}
    \mathcal Rf([(p,1)])=f(p),\quad\mathcal Rf([(p,2)])=\begin{pmatrix}
      1&0&0\\0&-1&0\\0&0&-1
    \end{pmatrix}f(p)\quad\text{for $p\in\hat\Sigma$}
  \end{equation}
  and using the boundary conditions in \eqref{eq:1-lem:rigidity-of-blow-ups}, one checks that $\mathcal Rf\in C^1\cap W^{2,2}_{\mathrm{loc}}(D(\hat\Sigma))$, and we equip $D(\hat\Sigma)$ with the pullback metric of $\mathcal R f$. Moreover, for each $p\in\partial\hat\Sigma\subset D(\hat\Sigma)$, the induced metric of $\mathcal Rf$ admits isothermal charts $\phi\colon (-1,1)^2\to D(\hat \Sigma)$ with $\phi(0,0)=p$ and $\phi((-1,1)\times\{0\})=\phi((-1,1)^2)\cap\partial\hat\Sigma$
  such that
  \begin{equation}
    \mathcal Rf(\phi(x,-y))=\begin{pmatrix}
      1&0&0\\0&-1&0\\0&0&-1
    \end{pmatrix} f(\phi(x,y))\quad\text{for all $x\in(-1,1)$, $y\in[0,1)$}.
  \end{equation}
  Thus, using $\mathbf H_f=0$ on $\partial\hat\Sigma$, \cite[Theorem~4.2]{kuwertlamm2021} yields that $\mathcal Rf$ is smooth and Willmore critical. In particular, since $\theta^2(\mathcal Rf,\infty)=2\theta^2(f,\infty)$, \cite[Lemma~4.3]{kuwertschaetzle2004} and \cite{palmer2000,dallacqua2012} yield the last part of the claim, using that also $\mathrm{genus}(D(\hat\Sigma))=0$ if $\mathrm{genus}(\hat\Sigma)=0$ and $\partial\hat\Sigma\cong\R$.
\end{proof}

\subsection{Stability and proof of the main result}

Let $\Sigma$ be a compact connected oriented surface with non-empty boundary, $\Gamma\subset\R^3$ a smooth compact one-dimensional submanifold and $\mathbf N$ a smooth unit normal field along $\Gamma$. Fix a smooth embedding $\gamma\colon\partial\Sigma\to\Gamma$ with $\gamma(\partial\Sigma)=\Gamma$.

We first sketch how classical arguments using a \L ojasiewicz--Simon gradient inequality and parabolic Schauder theory yield full convergence for $t\to\infty$ in the setting of \Cref{rem:subconvergence-in-case-a} if $\hat\Sigma$ is compact. Showing all details however is beyond the scope of this article. To this end, we first note that adapting the arguments of \cite{dallacquaruppschaetzleschlierf2026} to the much simpler setting of our Dirichlet boundary conditions, one can prove the following \L ojasiewicz--Simon inequality. 

\begin{theorem}\label{thm:loja}
  Let $\bar{f}\colon \Sigma\to\R^3$ be a smooth immersion satisfying $\bar f=\gamma$ and $\partial_{\eta_{\bar f}}\bar f=\mathbf N\circ\gamma$ on $\partial\Sigma$. Then there exist $\theta\in(0,\frac12]$, $C>0$, and $\sigma>0$ such that, for $f\in W^{4,2}(\Sigma,\R^3)$ with $\|f-\bar{f}\|_{W^{4,2}(\Sigma)}<\sigma$ satisfying $f=\gamma$ and $\partial_{\eta_{f}}f=\mathbf N\circ\gamma$ on $\partial\Sigma$, we have 
    \begin{equation}\label{eq:1-thm:loja}
        |\W(f)-\W(\bar{f})|^{1-\theta} \leq C \|\nabla\W(f)\|_{L^2(\dd\vol_{g_f})}.
    \end{equation}
\end{theorem}

As an application of \Cref{thm:loja}, using standard parabolic Schauder theory \cite{solonnikov1965mathphysics} and arguing as in \cite{chillfasangovaschaetzle2009}, one can show the following stability result. Also cf.~\cite[Lemma~7.9]{rupp2023}, \cite[Lemma~5.8]{rupp2024}, \cite[Lemma~6.2]{ruppscharrerschlierf2026}.

\begin{lemma}\label{lem:stability-loja-schauder-application}
  Let $\hat{f}\colon \Sigma\to\R^3$ be a Willmore immersion with $\hat f=\gamma$ and $\partial_{\eta_{\hat f}}\hat f=\mathbf N\circ\gamma$ on $\partial\Sigma$. Let $\alpha\in(0,1)$, $k\in \N$ with $k\geq 4$, and $\delta>0$. Then there exists $\varepsilon=\varepsilon(\hat{f},\alpha,k,\delta)>0$ such that if $f\colon [0,T)\times \Sigma\to\R^3$ solves \eqref{eq:Dirichlet-Willmore-flow}, and
		\begin{enumerate}[(i)]
			\item there exists a diffeomorphism $\Phi_0\colon\Sigma\to\Sigma$ with $\Vert {f_0\circ\Phi_0-\hat{f}}\Vert_{C^{k+\alpha}(\Sigma)}<\varepsilon$, and
			\item $\W(f(t))\geq \W(\hat{f})$ for all $t$ such that there exists a diffeomorphism $\Phi(t) \colon\Sigma\to\Sigma$ with the property that $\Vert {f(t)\circ \Phi(t) - \hat{f}}\Vert_{C^k(\Sigma)}\leq \delta$,
		\end{enumerate}
		then we may take $T=\infty$. Moreover, as $t\to\infty$, the flow converges smoothly after reparametrization by some diffeomorphisms $\tilde{\Phi}(t)\colon\Sigma\to\Sigma$ to a Willmore immersion $f_\infty$, such that $f_\infty\vert_{\partial\Sigma}$ is a reparametrization of $\gamma$, and $\partial_{\eta_{f_\infty}}f_\infty=\mathbf N\circ f_\infty$ on $\partial\Sigma$, and such that $\Vert {f_{\infty}-\hat{f}}\Vert_{C^k(\Sigma)}\leq\delta$.
\end{lemma}

\begin{corollary}[Full convergence for compact $\hat f$]\label{cor:full-convergence}
  In the setting of \Cref{rem:subconvergence-in-case-a}, suppose that $\hat\Sigma=\Sigma$ is compact. Then as $t\to\infty$, the flow $f(t)$ converges smoothly after reparametrization to a Willmore immersion $f_\infty$, such that $f_\infty\vert_{\partial\Sigma}$ is a reparametrization of $\gamma$, and $\partial_{\eta_{f_\infty}}f_\infty=\mathbf N\circ f_\infty$ on $\partial\Sigma$.
\end{corollary}
\begin{proof}
  By \Cref{rem:subconvergence-in-case-a}, for $j$ sufficiently large, $f_j(\tau_j)\circ\Phi_j$ satisfies the assumption (i) in \Cref{lem:stability-loja-schauder-application}. Assumption (ii) is clear using that $\partial_t\W(f(t,\cdot))\leq 0$ and $\W(f(t_j,\cdot))=\W(f_j(\tau_j)\circ\Phi_j)\to \W(\hat f)$ for $j\to\infty$. Thus, the statement is an immediate consequence of \Cref{lem:stability-loja-schauder-application}.
\end{proof}

\begin{proof}[Proof of \Cref{thm:main-result}]
  If $f(t)=f_0$ for all $t\in[0,T)$, then \eqref{eq:Dirichlet-Willmore-flow} and maximality yield $T=\infty$ and (b) with $f_\infty=f_0$. Otherwise, by \eqref{eq:energy-dissipation-Dirichlet}, there exists $s\in(0,T)$ such that
  \begin{equation}\label{eq:strict-energy-main-result}
    E\vcentcolon=\W(f(s))<\W(f_0)\leq\min\Big\{8\pi-2\sup_{z\in\R^3}\mathbf B_{\Gamma,\mathbf N}(z),4\pi-2\sup_{z\in\Gamma}\mathbf B_{\Gamma,\mathbf N}(z)\Big\}.
  \end{equation}
  Apply \Cref{prop:existence-of-blow-ups} to $f(s+\cdot)$. Suppose that case (b) of that theorem occurs, and let $\hat f\colon\hat\Sigma\to\R^3$ and $\hat p\in\hat\Sigma$ be as there. Then
  \begin{equation}
    \mathrm{genus}(\hat\Sigma)=0,\quad\int_{\hat\Sigma}|\mathbf A_{\hat f}|^2\dd\vol_{g_{\hat f}}<\infty,\quad[\mathbf n_{\hat f}]_{\alpha,1}<\infty,\quad\text{and}\quad |\mathbf A_{\hat f}(\hat p)|>0.
  \end{equation}
  First, suppose that $\partial\hat\Sigma=\emptyset$. By \Cref{prop:existence-of-blow-ups} applied to $f(s+\cdot)$ and \eqref{eq:strict-energy-main-result},
  \begin{equation}\label{eq:interior-density-main-result}
    \theta^2(\hat f,\infty)+\frac{\W(\hat f)}{4\pi}\leq\frac E{4\pi}+\frac1{2\pi}\sup_{z\in\R^3}\mathbf B_{\Gamma,\mathbf N}(z)<2.
  \end{equation}
  If $\hat\Sigma$ were compact, the exhaustion in \Cref{def:local-smooth-pointed-convergence} would yield $\hat\Omega_j=\hat\Sigma$ for all sufficiently large $j$. Then $\Omega_j=\hat\phi_j(\hat\Sigma)$ would be non-empty, open and closed in $\Sigma$, and hence $\Omega_j=\Sigma$. By \Cref{def:local-smooth-pointed-convergence}, this would imply $\partial\Sigma=\hat\phi_j(\partial\hat\Sigma)=\emptyset$, a contradiction. Thus $\hat\Sigma$ is noncompact. Moreover, \eqref{eq:monotonicity-with-boundary}, \eqref{eq:limits-monotonicity-with-boundary} and \eqref{eq:interior-density-main-result} yield
  \begin{equation}
    \#\hat f^{-1}(z)=\theta^2(\hat f,z)\leq\theta^2(\hat f,\infty)+\frac{\W(\hat f)}{4\pi}<2\quad\text{for all $z\in\hat f(\hat\Sigma)$}.
  \end{equation}
  Thus $\hat f$ is an embedding. Fix $x_0\notin\hat f(\hat\Sigma)$, and set
  \begin{equation}
    I_{x_0}(x)\vcentcolon=\frac{x-x_0}{|x-x_0|^2}\quad\text{and}\quad \bar\Sigma\vcentcolon=I_{x_0}(\hat f(\hat\Sigma))\cup\{0\}.
  \end{equation}
  By \cite[Lemma~4.3, Equations~(4.16) and~(4.20)]{kuwertschaetzle2004}, $\bar\Sigma$ is a smooth compact Willmore surface with
  \begin{align}
    \mathrm{genus}(\bar\Sigma)&=\mathrm{genus}(\hat\Sigma)=0,\quad\theta^2(\hat f,\infty)\geq1,\text{ and}\\
    \W(\bar\Sigma)&=\W(\hat f)+4\pi\leq\W(\hat f)+4\pi\theta^2(\hat f,\infty)<8\pi.
  \end{align}
  By \cite[Section~5, p.~51]{bryant1984}, $\bar\Sigma$ is a round sphere. Since $0\in\bar\Sigma$, its inverse image $\hat f(\hat\Sigma)=I_{x_0}^{-1}(\bar\Sigma\setminus\{0\})$ is a plane, contradicting $|\mathbf A_{\hat f}(\hat p)|>0$.

  If $\partial\hat\Sigma\neq\emptyset$, then \eqref{eq:density-bound-stationary-blow-up} applied to $f(s+\cdot)$ and \eqref{eq:strict-energy-main-result} yield
  \begin{equation}
    \theta^2(\hat f,\infty)\leq\theta^2(\hat f,\infty)+\frac{\W(\hat f)}{4\pi}\leq\frac E{4\pi}+\frac1{2\pi}\sup_{z\in\Gamma}\mathbf B_{\Gamma,\mathbf N}(z)<1.
  \end{equation}
  Together with \eqref{eq:properties-stationary-blow-up} and the boundary conditions in \Cref{prop:existence-of-blow-ups}(b), \Cref{lem:rigidity-bdry-blow-ups} implies that $\hat f$ parametrizes a flat half-plane. This again contradicts $|\mathbf A_{\hat f}(\hat p)|>0$. Consequently, \Cref{prop:existence-of-blow-ups}(a) and smoothness on $[0,s]\times\Sigma$ imply
  \begin{equation}\label{eq:bounded-curvature-main-result}
    T=\infty\quad\text{and}\quad M\vcentcolon=\max\Big\{1,\sup_{t\geq0}\Vert\mathbf A_{f(t)}\Vert_{L^\infty(\Sigma)}\Big\}<\infty.
  \end{equation}
  Fix $p_0\in\partial\Sigma$ and apply \Cref{rem:subconvergence-in-case-a}. Without loss of generality disregarding the fixed translation and dilation in \Cref{rem:subconvergence-in-case-a}, we get a complete connected oriented surface $\hat\Sigma$, a smooth proper Willmore immersion $\hat f\colon\hat\Sigma\to\R^3$, and $\hat p\in\partial\hat\Sigma$ such that, after passing to a subsequence,
  \begin{align}
    &\quad t_j\nearrow\infty,\quad(f(t_j),p_0)\longrightarrow(\hat f,\hat p)\quad\text{in the sense of \Cref{def:local-smooth-pointed-convergence}, and}\label{eq:subconvergence-main-result}\\
    &\mathrm{genus}(\hat\Sigma)=0,\quad\Vert\mathbf A_{\hat f}\Vert_{L^\infty(\hat\Sigma)}\leq M,\quad\hat f\vert_{\partial\hat\Sigma}\colon\partial\hat\Sigma\to\Gamma\text{ is an embedding}.
  \end{align}
  By \eqref{eq:Dirichlet-Willmore-flow} and \eqref{eq:definition-pointed-convergence},
  \begin{equation}
    \partial_{\eta_{\hat f}}\hat f=\mathbf N\circ\hat f\quad\text{on $\partial\hat\Sigma$}.
  \end{equation}
  If $\hat\Sigma$ is noncompact, this proves (a). Otherwise, \Cref{rem:subconvergence-in-case-a} yields that $\hat\Sigma=\Sigma$ and diffeomorphisms $\Phi_j\colon\Sigma\to\Sigma$ such that
  \begin{equation}
    f(t_j)\circ\Phi_j\longrightarrow\hat f\quad\text{in $C^\infty(\Sigma,\R^3)$}.
  \end{equation}
  Applying \Cref{cor:full-convergence} yields (b). The boundary conditions for $f_\infty$ follow from \eqref{eq:Dirichlet-Willmore-flow}, and the restriction $f_\infty\vert_{\partial\Sigma}$ is an embedding onto $\Gamma$ by \Cref{rem:subconvergence-in-case-a}, and hence a reparametrization of $\gamma$.

  To see that (a) and (b) are mutually exclusive, suppose that (b) holds, including the stationary case considered above. Then there are diffeomorphisms $\Phi_t\colon\Sigma\to\Sigma$ such that
  \begin{equation}
    f(t)\circ\Phi_t\longrightarrow f_\infty\text{ in }C^\infty(\Sigma,\R^3),\quad\operatorname{diam}_{g_{f(t)}}(\Sigma)=\operatorname{diam}_{g_{f(t)\circ\Phi_t}}(\Sigma)\leq C\quad\text{for all sufficiently large $t$}.
  \end{equation}
  For any limit as in (a), \eqref{eq:intrinsic-exhaustion-pointed-convergence} with $R>C$ yields that $\Omega_j=\Sigma$ for all sufficiently large $j$. Thus the sets $\hat\Omega_j=\hat\phi_j^{-1}(\Sigma)$ are compact and connected for all sufficiently large $j$. Since $\hat\Omega_j$ is both open and closed in $\hat\Omega_{j+1}$, they coincide, and $\hat\Sigma=\bigcup_j\hat\Omega_j$ is compact, a contradiction.

  Finally, suppose that $\W(f_0)\leq4\pi$. The stationary case has already been treated, and otherwise \eqref{eq:strict-energy-main-result} yields that $E<4\pi$. Put $L\vcentcolon=\mathcal H^1(\Gamma)$. By \eqref{eq:monotonicity-with-boundary}, \eqref{eq:limits-monotonicity-with-boundary} and \eqref{eq:energy-dissipation-Dirichlet}, for $t\geq s$ and $z\in f(t,\Sigma)\setminus\Gamma$,
  \begin{equation}
    4\pi\leq4\pi\theta^2(f(t),z)\leq\W(f(t))+2\mathbf B_{\Gamma,\mathbf N}(z)\leq E+\frac{2L}{\operatorname{dist}(z,\Gamma)},
  \end{equation}
  also cf.\ \cite[First equation in proof of Theorem~4.1]{novagapozzetta2020} where such an argument is already used. Fix $q\in\Gamma$. It follows that
  \begin{equation}
    f(t,\Sigma)\subset\overline B_R(q)\quad\text{for all $t\geq s$}\quad\text{with}\quad R\vcentcolon=\sup_{p\in\Gamma}|p-q|+\frac{2L}{4\pi-E}<\infty.
  \end{equation}
  For any limit in (a), \eqref{eq:definition-pointed-convergence} therefore implies $\hat f(\hat\Sigma)\subset\overline B_R(q)$. Since $\hat f$ is proper, $\hat\Sigma=\hat f^{-1}(\overline B_R(q))$ is compact, a contradiction.
\end{proof}

\section*{Acknowledgments}
This research was funded in whole, or in part, by the Austrian Science Fund (FWF), grant number \href{https://doi.org/10.55776/PAT8894924}{10.55776/PAT8894924}.

\appendix
\crefalias{section}{appendix}%
\crefalias{subsection}{appendix}

\section{Existence of local graph charts and compactness}
\label{app:graph-rep}

\subsection{Local graph charts with estimates}

For simplicity, assume that $\Sigma$ is an orientable surface and $f\in C^{1+\alpha}(\Sigma,\R^3)$ is an immersion. We assume that $(\Sigma,d_{g_f})$ is complete, with the boundary included. Denote by $\mathbf n_f\colon\Sigma\to\S^2\subset\R^3$ the outer unit normal induced by the orientation of $\Sigma$.

Set $Q_{r}=(-r, r)^2$ and $Q_{r}^+=(-r, r)\times[0,r)$. Let $p\in \Sigma$ and consider a Euclidean isometry $A_p$ of $\R^3$ with 
\begin{equation}\label{eq:choice-Ap-1}
  A_p(0)=f(p)\quad\text{and}\quad DA_p.e_3=\mathbf n_f(p).
\end{equation}
If in addition $p\in\partial\Sigma$, we choose $A_p$ such that 
\begin{equation}\label{eq:choice-Ap-2}
  DA_p.e_2+\partial_{\eta_f} f(p)=0.
\end{equation}
For $\pi\colon((x^{(1)},\dots,x^{(3)})\in\R^3\mapsto(x^{(1)},x^{(2)})\in\R^2)$, let $U_{r,p}\subset\Sigma$ be the connected component of $(\pi\circ A_p^{-1}\circ f)^{-1}(Q_r)$ containing $p$.

This section's main result \Cref{prop:ex-of-graph-charts} is a quantitative version of the varifold-result \cite{bournivolkmann2016} for immersed $C^{1+\alpha}$ surfaces which does not a priori assume a small-mass-type bound on the scale on which local graphs should exist, and instead proves it after choosing a suitable ``graphical sheet''. One may also view it as a $C^{1+\alpha}$-generalization, to surfaces \emph{with boundary}, of the construction used in \cite{langer1985} in the case of closed surfaces. We define
\begin{equation}\label{eq:def-nf-alpha-seminorm}
  [\mathbf n_f]_{\alpha,r,p}\vcentcolon= \sup_{x,y\in U_{r,p} \text{ s.th.\ } f(x)\neq f(y)} \frac{|\mathbf n_f(x)-\mathbf n_f(y)|}{|f(x)-f(y)|^\alpha},\quad [\mathbf n_f]_{\alpha,r}\vcentcolon= \sup_{p\in\Sigma}\ [\mathbf n_f]_{\alpha,r,p}.
\end{equation}

\begin{definition}\label{def:rGamma}
  Let $\Gamma\subset\R^3$ be a smooth (not necessarily connected) one-dimensional submanifold. We define $\mathfrak r_\Gamma\in[0,\infty]$ to be the supremum of all $r>0$ with the property that, for all $0<s\leq r$ and $p,q\in\Gamma$, we have
  \begin{equation}
    d_\Gamma(p,q)\geq s\Longrightarrow|p-q|\geq\frac12s
  \end{equation}
  where $d_\Gamma\colon\Gamma\times\Gamma\to[0,\infty]$ denotes the intrinsic distance on $\Gamma$. We use the convention $\sup\emptyset=0$.
\end{definition}

\begin{proposition}\label{prop:ex-of-graph-charts}
  There exist universal constants $C_{\ref{prop:ex-of-graph-charts}}$, $\varepsilon_{\ref{prop:ex-of-graph-charts}}>0$ with the following property. Let $f\in C^{1+\alpha}(\Sigma,\R^3)$ be an immersion such that $\Gamma=f(\partial\Sigma)$ is a smooth submanifold of $\R^3$ with unit tangent field $\bm\tau_\Gamma$. Moreover, suppose that $p\in\Sigma$ and $r>0$ satisfy
  \begin{equation}\label{eq:1-prop:ex-of-graph-charts}
    r^\alpha [\mathbf n_f]_{\alpha,r,p} + r^\alpha [\bm\tau_\Gamma]_{C^\alpha(\Gamma)} \leq\varepsilon\leq \varepsilon_{\ref{prop:ex-of-graph-charts}}.
  \end{equation}
  With $A_p$ as in \eqref{eq:choice-Ap-1} and \eqref{eq:choice-Ap-2}, $P_p\vcentcolon=\pi\circ A_p^{-1}\circ f\colon\Sigma\to\R^2$ satisfies the following.
  \begin{enumerate}[(a)]
    \item If $p\in\partial\Sigma$ and $16r<\mathfrak r_\Gamma$, then $D_{r,p}\vcentcolon=P_p(U_{r,p})$ is the intersection of $Q_r$ with the super-level-set of a $C^{1+\alpha}$-function on $(-r,r)$, and $P_p\vert_{U_{r,p}}\colon U_{r,p}\to D_{r,p}$ is a $C^{1+\alpha}$-diffeomorphism.
    \item If $p\in\Sigma$ with $d_{g_f}(p,\partial\Sigma)\geq 2\sqrt{2}r$, then $D_{r,p}\vcentcolon=P_p(U_{r,p})$ satisfies $D_{r,p}=Q_r$, and $P_p\vert_{U_{r,p}}\colon U_{r,p}\to D_{r,p}$ is a $C^{1+\alpha}$-diffeomorphism.
  \end{enumerate}
  In both cases, denoting $\Psi_p=(P_p\vert_{U_{r,p}})^{-1}\colon D_{r,p}\to U_{r,p}$, we have
  \begin{equation}\label{eq:2-prop:ex-of-graph-charts}
    f(\Psi_p(x)) = A_p(x,u_p(x))\quad\text{for all $x\in D_{r,p}$}
  \end{equation}
  for some $u_p\in C^{1+\alpha}(D_{r,p},\R)$ which satisfies the estimate
  \begin{equation}\label{eq:3-prop:ex-of-graph-charts}
    \frac{\Vert u_p\Vert_{C^0(D_{r,p})}}{r}+\Vert \nabla u_p\Vert_{C^0(D_{r,p})} + r^\alpha [\nabla u_p]_{C^\alpha(D_{r,p})} \leq C_{\ref{prop:ex-of-graph-charts}} \varepsilon.
  \end{equation}
  Moreover, we have
  \begin{equation}\label{eq:4-prop:ex-of-graph-charts}
    d_{g_f}(z,p)\geq r\text{ for all $z\in\Sigma\setminus U_{r,p}$}\quad\text{and}\quad d_{g_f}(z,p)<(\sqrt{2}+C_{\ref{prop:ex-of-graph-charts}}\varepsilon)r\text{ for all $z\in U_{r,p}$},
  \end{equation}
  and
  \begin{align}
    \Big\vert\int_{U_{r,p}}\dd\vol_{g_f} - 2r^2\Big\vert &\leq C_{\ref{prop:ex-of-graph-charts}}\varepsilon r^2\quad\text{in case (a), and}\\
    4r^2&\leq \int_{U_{r,p}}\dd\vol_{g_f}\leq \big(4+C_{\ref{prop:ex-of-graph-charts}}\varepsilon \big)r^2\quad\text{in case (b)}.\label{eq:5-prop:ex-of-graph-charts}
  \end{align}
\end{proposition}
\begin{proof} 
  Without loss of generality, assume that $A_p=\operatorname{id}$ and $r=1$. With this normalization, $P=\pi\circ f$ and \eqref{eq:choice-Ap-1}, \eqref{eq:choice-Ap-2} yield that $\mathbf n_f(p)=e_3$, $T_{f(p)}\Gamma=\R e_1$, and $\partial_\eta f(p)=-e_2$ if $p\in\partial\Sigma$. We start by showing that, for $\varepsilon_{\ref{prop:ex-of-graph-charts}}$ sufficiently small, we have that
  \begin{equation}\label{eq:1-pf-prop:ex-of-graph-charts}
    dP_x\text{ is invertible with } \Vert (dP_x)^{-1}\Vert_{L(\R^2,(T_x\Sigma,g_f))}\leq 2\quad\text{for all $x\in U_{1,p}$ with $\vert f(x)\vert\leq 8$}.
  \end{equation}
  Indeed, using that $P=\pi\circ f$ and that $\mathbf n_f(x)\perp df_x(T_x\Sigma)$, for $x\in U_{1,p}$ and $v\in T_x\Sigma$ with $\vert v\vert_{g_f}=1$, we have
  \begin{align}
    \vert dP_x.v\vert^2 &= \vert df_x.v\vert ^2 - \langle df_x.v,e_3\rangle^2 = \vert df_x.v\vert ^2 - \langle (I-\mathbf n_f(x)\cdot\mathbf n_f(x)^t).df_x.v,e_3\rangle^2\\
    &\geq \vert df_x.v\vert ^2 - \vert df_x.v\vert ^2(1-\langle \mathbf n_f(x),e_3\rangle^2) =  \vert df_x.v\vert^2\big(1-\frac12\vert\mathbf n_f(x)-e_3\vert^2\big)^2=\big(1-\frac12\vert\mathbf n_f(x)-e_3\vert^2\big)^2.
  \end{align}
  After shrinking $\varepsilon_{\ref{prop:ex-of-graph-charts}}$ if necessary, $\mathbf n_f(p)=e_3$ and \eqref{eq:1-prop:ex-of-graph-charts} yield $\vert \mathbf n_f(x)-e_3\vert\leq 1$ for all $x\in U_{1,p}$ with $\vert f(x)\vert\leq 8$ and thus, we obtain \eqref{eq:1-pf-prop:ex-of-graph-charts}.

  First consider the case (a), i.e.\ $p\in\partial\Sigma$ and $\mathfrak r_\Gamma>16$. Let $\gamma\colon(-16,16)\to\Gamma$ be a unit-speed curve with $\gamma(0)=f(p)$ and $\gamma'(0)=e_1$. Then \eqref{eq:1-prop:ex-of-graph-charts} with $r=1$ yields that
  \begin{equation}\label{eq:2-pf-prop:ex-of-graph-charts}
    \vert \gamma'(s)-e_1\vert \leq \varepsilon |s|^\alpha \quad\text{for all $s\in(-16,16)$}.
  \end{equation}
  Thus, choosing $\varepsilon_{\ref{prop:ex-of-graph-charts}}$ sufficiently small, there exists $\varphi=(\varphi^{(1)},\varphi^{(2)})\in C^{1+\alpha}((-1,1),\R^2)$ such that, for a universal constant $C>0$,
  \begin{equation}\label{eq:3-pf-prop:ex-of-graph-charts}
    \gamma((-16,16))\cap (-1,1)^3 = \{(x,\varphi(x))\mid x\in(-1,1)\}\quad\text{and}\quad \Vert \varphi\Vert_{C^{1+\alpha}(-1,1)} \leq C\varepsilon < \frac12.
  \end{equation}
  Moreover, there exists $\eta\colon(-1,1)\to\partial\Sigma\cap U_{1,p}$ such that $P(\eta(x))=(x,\varphi^{(1)}(x))$ for all $x\in(-1,1)$. By \eqref{eq:3-pf-prop:ex-of-graph-charts}, we have $\vert f(\eta(x))\vert \leq \sqrt{2}$ for all $x\in(-1,1)$. Our next goal is to construct a $C^{1+\alpha}$-map $\Psi\colon D=\{(x,y)\in Q_1\mid \varphi^{(1)}(x)\leq y<1\}\to U_{1,p}$ such that $P\circ\Psi=\operatorname{id}$. 

  Since $\vert f\circ\eta\vert\leq \sqrt{2}$, \eqref{eq:1-pf-prop:ex-of-graph-charts} and the inverse function theorem yield that $P$ is locally invertible in a neighborhood of $\eta(x)$, for all $x\in(-1,1)$. Fix any $x\in(-1,1)$. Using $P(\eta(x))=(x,\varphi^{1}(x))$, we can apply the local inverse of $P$ at $\eta(x)$ to the vertical segment $y\mapsto (x,\varphi^{1}(x)+y)$, there exist a maximal $T_x\in (0,1-\varphi^{(1)}(x)]$ and a $C^{1+\alpha}$-curve $\tilde\gamma_x\colon [0,T_x)\to U_{1,p}$ such that 
  \begin{equation}
    \tilde\gamma_x(0)=\eta(x),\quad\text{and}\quad P(\tilde\gamma_x(y))=(x,\varphi^{1}(x)+y) \text{ and } \vert f(\tilde\gamma_x(y))\vert\leq 8 \text{ for all $y\in[0,T_x)$}.
  \end{equation}
  Since $\tilde\gamma_x'(y)=(dP_{\tilde\gamma_x(y)})^{-1}.(0,1)$, \eqref{eq:1-pf-prop:ex-of-graph-charts} and \eqref{eq:3-pf-prop:ex-of-graph-charts} yield that $\mathcal L_{(\Sigma,g_f)}(\tilde\gamma_x)\leq 4$ and thus
  \begin{equation}\label{eq:4-pf-prop:ex-of-graph-charts}
    \vert f(\tilde\gamma_x(t))-f(p)\vert \leq \mathcal L_{(\Sigma,g_f)}(\tilde\gamma_x) + \vert (x,\varphi(x))-f(p)\vert \leq4+\sqrt2<8.
  \end{equation}
  For the sake of contradiction, suppose $T_x<1-\varphi^{(1)}(x)$. Then $q=\lim_{y\nearrow T_x}\tilde\gamma_x(y)\in U_{1,p}$. If $q\in \partial\Sigma$, $\mathfrak r_\Gamma>16$ and \eqref{eq:4-pf-prop:ex-of-graph-charts} yield $d_\Gamma(f(p),f(q))<16$, which contradicts \eqref{eq:3-pf-prop:ex-of-graph-charts}. Therefore, $q\notin \partial\Sigma$ and, by \eqref{eq:4-pf-prop:ex-of-graph-charts}, $\vert f(q)\vert\leq 8$. Thus \eqref{eq:1-pf-prop:ex-of-graph-charts} shows that $P$ is locally invertible at $q$ which contradicts the maximality of $T_x$. This yields that $T_x=1-\varphi^{(1)}(x)$, for all $x\in (-1,1)$. A similar argument using local inverses of $P$ also shows that the trajectories of $\gamma_x$ and $\gamma_{x'}$ do not intersect for $x\neq x'$.

  By setting $\Psi(z)=\tilde\gamma_{z_x}(z_y-\varphi^{(1)}(z_x))$ for all $z=(z_x,z_y)\in D$, this construction yields an injective continuous map $\Psi\colon D\to U_{1,p}$ with $P\circ\Psi=\operatorname{id}$ such that $P$ is locally invertible at each point in $\Psi(D)$. Thus, $\Psi$ is a homeomorphism onto its image. Since its image is also relatively closed in $U_{1,p}$, we get that $\Psi(D)=U_{1,p}$ which finishes the proof of case (a).

  Next we turn to case (b), i.e.\ $p\in\Sigma$ with $d_{g_f}(p,\partial\Sigma)>2\sqrt{2}$. By the inverse function theorem, there exists $r_0\in(0,1)$ sufficiently small such that $P\vert_{U_{r_0,p}}$ is a $C^{1+\alpha}$-diffeomorphism onto its image. So for each $v\in\R^2$ with $\vert v\vert =1$, we can choose a maximal curve $\tilde\gamma_v$ with $\tilde\gamma_v(0)=p$, $P\circ\tilde\gamma_v(y)=y\cdot v$, and $\vert f\circ\tilde\gamma_v\vert \leq 8$. By \eqref{eq:1-pf-prop:ex-of-graph-charts}, we again get $d_{g_f}(\tilde\gamma_v(y),p)< 2\sqrt{2}$ for all $y$ in the maximal existence interval, and thus $\gamma_v$ cannot terminate at a boundary point since $d_{g_f}(p,\partial\Sigma)>2\sqrt{2}$ by assumption. Then case (b) follows analogously to the argument used in case (a) above.
  
  In both cases, let $\Psi=(P\vert_{U_{1,p}})^{-1}\colon D_{1,p}\to U_{1,p}$. Setting $u(x)\vcentcolon=f^{(3)}(\Psi(x))$ yields $u\in C^{1+\alpha}(D_{1,p})$ and \eqref{eq:2-prop:ex-of-graph-charts}. With $F\colon D_{1,p}\to\R^3$, $F(x)=(x,u(x))$, we thus have $\mathbf n_F=\mathbf n_f\circ \Psi$ and
  \begin{equation}\label{eq:5-pf-prop:ex-of-graph-charts}
    \nabla u = - \frac{\pi.\mathbf n_F}{\langle\mathbf n_F,e_3\rangle}  = - \frac{\pi.\big(\mathbf n_f\circ\Psi\big)}{\langle \mathbf n_f\circ\Psi,e_3\rangle} .
  \end{equation}
  Revisiting the above construction of $\Psi$, we have $\vert n_f-e_3\vert\leq 1$ on $U_{1,p}$. Therefore, $\langle \mathbf n_f,e_3\rangle \geq 1/2$ so that \eqref{eq:5-pf-prop:ex-of-graph-charts} yields that $F$ is Lipschitz on $D_{1,p}$ with a universal Lipschitz constant. So using $\vert n_F(x)-n_F(y)\vert\leq [\mathbf n_f]_{\alpha,1,p}\vert F(x)-F(y)\vert^\alpha$ for $x,y\in D_{1,p}$, \eqref{eq:5-pf-prop:ex-of-graph-charts} and \eqref{eq:1-prop:ex-of-graph-charts} yield \eqref{eq:3-prop:ex-of-graph-charts} for a suitable universal constant $C_{\ref{prop:ex-of-graph-charts}}>0$.

  The first part of \eqref{eq:4-prop:ex-of-graph-charts} is due to the fact that $\pi$ is a contraction. In the interior case (b), the second estimate in \eqref{eq:4-prop:ex-of-graph-charts} is an immediate consequence of \eqref{eq:3-prop:ex-of-graph-charts}. For case (a), we also use that \eqref{eq:3-pf-prop:ex-of-graph-charts} yields that any $z\in D_{1,p}$ can be connected to $p$ by a piecewise $C^{1}$-curve in $D_{1,p}$ of length bounded by $\sqrt 2 + C\varepsilon$ for some universal constant $C>0$. Similarly, the second part of \eqref{eq:5-prop:ex-of-graph-charts} follows immediately from \eqref{eq:3-prop:ex-of-graph-charts} and $D_{1,p}=Q_1$ in case (b). For the boundary case (a), we additionally use that \eqref{eq:3-pf-prop:ex-of-graph-charts} yields $\vert \mathcal H^2(D_{1,p})-2\vert \leq C\varepsilon$ for some universal constant $C>0$.
\end{proof}
\begin{remark}\label{rem:Af-Linfty-Lipschitz-est}
  Arguing similarly as in the proof of \Cref{prop:ex-of-graph-charts}, one can show that there exist universal constants $\varepsilon_{\ref{rem:Af-Linfty-Lipschitz-est}}>0$ and $C_{\ref{rem:Af-Linfty-Lipschitz-est}}>1$ with the following property. Suppose that $f\in C^2(\Sigma,\R^3)$ is an immersion such that $\Gamma=f(\partial\Sigma)$ is a smooth submanifold of $\R^3$ with unit tangent field $\bm\tau_\Gamma$ and curvature vector $\bm\kappa_\Gamma=D_s\bm\tau_\Gamma$ such that, for $0<\varepsilon\leq\varepsilon_{\ref{rem:Af-Linfty-Lipschitz-est}}$, $r>0$ with $16r<\mathfrak r_\Gamma$ and $p\in\Sigma$, we have
  \begin{equation}
    r\big(\Vert \bm\kappa_\Gamma\Vert_{C^0(\Gamma)} + \Vert A_f\Vert_{L^\infty(U_{r,p})}\big) \leq \varepsilon\leq \varepsilon_{\ref{rem:Af-Linfty-Lipschitz-est}}.
  \end{equation} 
  Then 
  \begin{equation}
    \operatorname{diam}(f(U_{r,p})) \leq C_{\ref{rem:Af-Linfty-Lipschitz-est}} r\quad\text{and}\quad \vert \mathbf n_f(x)-\mathbf n_f(y)\vert \leq C_{\ref{rem:Af-Linfty-Lipschitz-est}} \Vert A_f\Vert_{L^\infty(U_{r,p})} \vert f(x)-f(y)\vert\text{ for all $x,y\in U_{r,p}$}.
  \end{equation}
  In particular, $r^\alpha[\mathbf n_f]_{\alpha,r,p}\leq C_{\ref{rem:Af-Linfty-Lipschitz-est}}^{2-\alpha}\varepsilon$ for all $\alpha\in(0,1)$. 
\end{remark}

\subsection{Compactness for immersions with boundary}
\label{app:compactness-with-boundary}

We use the following notion of convergence, including convergence up to the boundary. Completeness below always refers to the induced intrinsic metric with the boundary included.

\begin{definition}\label{def:local-smooth-pointed-convergence}
  Let $f_j\colon\Sigma_j\to\R^3$ and $\hat f\colon\hat\Sigma\to\R^3$ be smooth immersions, with $p_j\in\Sigma_j$ and $\hat p\in\hat\Sigma$. We say that $(f_j,p_j)\to(\hat f,\hat p)$ in the sense of \emph{local smooth pointed convergence} if there exist relatively open sets $\hat\Omega_j\Subset\hat\Omega_{j+1}\subset\hat\Sigma$ and diffeomorphisms $\hat\phi_j\colon\hat\Omega_j\to\Omega_j\subset\Sigma_j$ onto relatively open subsets such that
  \begin{align}
    &\hat\Sigma=\bigcup_{j\in\N}\hat\Omega_j,\quad \hat\phi_j(\partial\hat\Sigma\cap\hat\Omega_j)=\partial\Sigma_j\cap\Omega_j,\\
    &f_j\circ\hat\phi_j\longrightarrow\hat f\text{ in }C^\infty_{\mathrm{loc}}(\hat\Sigma,\R^3),\quad p_j\in\Omega_j,\quad\hat\phi_j^{-1}(p_j)\longrightarrow\hat p,\label{eq:definition-pointed-convergence}
  \end{align}
  and, for every $R<\infty$,
  \begin{equation}\label{eq:intrinsic-exhaustion-pointed-convergence}
    \{x\in\Sigma_j\mid d_{g_{f_j}}(x,p_j)<R\}\subset\Omega_j\quad\text{for all sufficiently large $j$}.
  \end{equation}
  The convergence in \eqref{eq:definition-pointed-convergence} is understood in charts of manifolds with boundary, on compact subsets including boundary points.
\end{definition}

\begin{proposition}[Compactness with boundary]\label{prop:compactness-with-boundary}
  Let $f_j\colon\Sigma_j\to\R^3$ be smooth proper immersions of connected oriented surfaces with boundary. Suppose that $\Gamma_j\subset\R^3$ are smooth compact one-dimensional submanifolds, $\mathbf N_j$ smooth unit normal fields along $\Gamma_j$, and
  \begin{equation}\label{eq:boundary-conditions-compactness}
    f_j\vert_{\partial\Sigma_j}\colon\partial\Sigma_j\to\Gamma_j\text{ is an embedding onto $\Gamma_j$},\quad \partial_{\eta_{f_j}}f_j=\mathbf N_j\circ f_j\text{ on $\partial\Sigma_j$}.
  \end{equation}
  Assume that $p_j\in\Sigma_j$, $f_j(p_j)=0$, and that, for all $R>0$ and $m\in\N_0$,
  \begin{equation}\label{eq:area-and-curvature-compactness}
    \int_{f_j^{-1}(B_R(0))}\dd\vol_{g_{f_j}}\leq C(R),\quad\Vert\nabla_{g_{f_j}}^m \mathbf A_{f_j}\Vert_{L^\infty(\Sigma_j)}\leq C_m.
  \end{equation}
  Furthermore, suppose that
  \begin{equation}\label{eq:boundary-geometry-compactness}
      \inf_{j\in\N}\mathfrak r_{\Gamma_j}>0,\text{ and }\Vert D_s^m\bm\kappa_{\Gamma_j}\Vert_{C^0(\Gamma_j)}+\Vert D_s^{m+1}\mathbf N_j\Vert_{C^0(\Gamma_j)}\leq C_m\quad\text{for all $j\in\N$ and $m\in\N_0$}.
  \end{equation}
  All derivatives of the boundary data are ambient arclength derivatives. 

  Then there exist a subsequence and a smooth proper immersion $\hat f\colon\hat\Sigma\to\R^3$ of a complete connected oriented surface, with $\hat p\in\hat\Sigma$ and $\hat f(\hat p)=0$, such that $(f_j,p_j)\to(\hat f,\hat p)$ in the sense of \Cref{def:local-smooth-pointed-convergence}. The bounds in \eqref{eq:area-and-curvature-compactness} also hold for $\hat f$. Moreover, $\hat f\vert_{\partial\hat\Sigma}$ is an embedding, and
  \begin{equation}\label{eq:boundary-alternative-compactness}
    \partial\hat\Sigma=\emptyset\quad\Longleftrightarrow\quad d_{g_{f_j}}(p_j,\partial\Sigma_j)\longrightarrow\infty.
  \end{equation}
\end{proposition}

This is the boundary version of \cite[Theorem~1.3]{breuning2015}. A detailed proof is beyond the scope of this article, so we only indicate the modifications of the argument used in \cite[Theorem~1.3]{breuning2015}. By \Cref{rem:Af-Linfty-Lipschitz-est,prop:ex-of-graph-charts}, \eqref{eq:area-and-curvature-compactness} and \eqref{eq:boundary-geometry-compactness}, there are interior and boundary graph charts with a uniform radius and uniform bounds for all derivatives. Flatten the graphical boundary in the boundary charts. The area lower bound in \eqref{eq:5-prop:ex-of-graph-charts} and the upper bound in \eqref{eq:area-and-curvature-compactness} bound the number of charts in each ambient ball. Proceeding as in \cite[Section~7]{breuning2015}, use the limits of the graph systems and transition maps, now also on half-charts, and boundary-preserving reparametrizations on collars. Hereby, the first condition in \eqref{eq:boundary-geometry-compactness} ensures that injectivity of the boundary maps is preserved in the limit.

\section{Function spaces and parabolic theory on surfaces with boundary}

The following hypothesis is motivated by the more general concept of \emph{uniformly regular Riemannian manifolds}, for instance see \cite[Definition~3.3]{amann2025}. Here we use the notation $Q_r=(-r,r)^2$ and $Q_r^+=(-r,r)\times[0,r)$ for $r>0$.

\begin{hypothesis}[Uniformly regular surfaces with boundary]\label{hyp:unif-reg}
  Assume that $\Sigma$ is a compact surface with non-empty boundary and a corresponding \emph{good atlas}. That is, for a universal constant $\rho^\ast>0$,  there are $N_1,N_2\in\N$ and mappings $\alpha_j\colon D_j\vcentcolon=Q_{\rho^\ast} \to \Sigma$ for $j\in \mathbf I_{\mathrm{int}}=\{1,\dots,N_1\}$ and $\alpha_j\colon D_j\vcentcolon=Q_{16\rho^\ast}^+ \to \Sigma$ for $j\in \mathbf I_{\mathrm{bdr}}=\{N_1+1,\dots,N_1+N_2\}$ which are diffeomorphic onto their respective images, and let $\Omega_j=\alpha_j(D_j)$ for $j\in \mathbf I_{\mathrm{int}}\cup \mathbf I_{\mathrm{bdr}}$.  The inverses are denoted by $\varphi_j$. Moreover, we assume that the following holds:
  \begin{enumerate}[(a)]
    \item $\displaystyle 1\leq \sum_{j\in \mathbf I_{\mathrm{int}}} \chi_{\alpha_j(Q_{\rho^\ast/2})} + \sum_{j\in \mathbf I_{\mathrm{bdr}}} \chi_{\alpha_j(Q_{12\rho^\ast}^+)} \leq \sum_{j\in \mathbf I_{\mathrm{int}}\cup \mathbf I_{\mathrm{bdr}}} \chi_{\Omega_j}\leq Z^\ast$ for a universal constant $Z^\ast\in\N$.\label{it:1-hyp:good-atlas}
    \item $|\partial^\beta (\varphi_j\circ\alpha_k)|\leq C_{\vert\beta\vert}$ on $\varphi_k(\Omega_j\cap\Omega_k)$ for any multiindex $\beta\in\N_0^2$.\label{it:2-hyp:good-atlas}
  \end{enumerate}
\end{hypothesis}

This concept plays a crucial role in the analysis of the Willmore flow in \Cref{thm:lifespan-and-smoothing}, as it allows us to obtain universal constants in global a-priori estimates for the linearized system after localizing via the charts of the \emph{good atlas} and the \emph{subordinate localization system} we introduce below.

\subsection{Function spaces and localization system}

\begin{remark}[Subordinate localization system]
\label{rem:local-sys}
  Assuming \Cref{hyp:unif-reg}, we can follow the construction for the \emph{localization system} in \cite[Lemma~3.2]{amann2013}.
  Choose nonnegative smooth coordinate cutoffs $\zeta_j$ equal to one on $Q_{\rho^\ast/2}$ resp.\ $Q_{12\rho^\ast}^+$, with supports in $Q_{3\rho^\ast/4}$ resp.\ $Q_{25\rho^\ast/2}^+$, for $j\in \mathbf I_{\mathrm{int}}$ resp.\ $j\in \mathbf I_{\mathrm{bdr}}$. Set $\pi_j=\zeta_j/(\sum_\ell\zeta_\ell^2)^{1/2}$, so that $\sum_j\pi_j^2=1$. 
 
  Using \Cref{hyp:unif-reg}\eqref{it:1-hyp:good-atlas} and \eqref{it:2-hyp:good-atlas}, we have
  \begin{equation}\label{eq:unif-bounds-pij}
    \Vert\partial^\beta (\pi_j\circ\alpha_k)\Vert_{C^0(\varphi_k(\Omega_j\cap\Omega_k))} \leq C_{|\beta|}\quad\text{for all multiindices $\beta\in\N_0^2$}. 
  \end{equation}
  This allows defining a \emph{reference metric} $\bar g$ on $\Sigma$ by $\bar g=\sum_ {j\in \mathbf I_{\mathrm{int}}\cup \mathbf I_{\mathrm{bdr}}}\pi_j^2 \varphi_j^*g_{\R^2}$.
\end{remark}

Throughout this section, assume \Cref{hyp:unif-reg} and fix the subordinate localization system from \Cref{rem:local-sys}.

\begin{definition}\label{lem:unif-equivalent-norms}
  Let $I_r=(-r,r)$ for $r>0$. For smooth $u$, we define by localization and pullback
  \begin{align}
    \Vert u\Vert_{\bm C^s(\Sigma)} &\vcentcolon= \max_{j\in \mathbf I_{\mathrm{int}}\cup \mathbf I_{\mathrm{bdr}}} \Vert (\pi_j u)\circ\alpha_j\Vert_{C^s(D_j)}\quad\text{and}\\
    \Vert u\Vert_{\bm C^s(\partial\Sigma)} &\vcentcolon= \max_{{j\in \mathbf I_{\mathrm{bdr}}}} \Vert (\pi_j u)\circ\alpha_j(\cdot,0)\Vert_{C^s(I_{16\rho^\ast})}\quad\text{for $0\leq s<\infty$};\\
    \Vert u\Vert_{\bm W^{s,p}(\Sigma)} &\vcentcolon= \max_{{j\in \mathbf I_{\mathrm{int}}\cup \mathbf I_{\mathrm{bdr}}}} \Vert (\pi_j u)\circ\alpha_j\Vert_{W^{s,p}(D_j)}\quad\text{and}\\
    \Vert u\Vert_{\bm W^{s,p}(\partial\Sigma)} &\vcentcolon= \max_{{j\in \mathbf I_{\mathrm{bdr}}}} \Vert (\pi_j u)\circ\alpha_j(\cdot,0)\Vert_{W^{s,p}(I_{16\rho^\ast})}\quad\text{for $s\geq 0$ and $1\leq p\leq\infty$}.
  \end{align}
  The \emph{Sobolev space} $W^{k,p}(\Sigma)$ is defined as the closure of $C^\infty(\Sigma)$ in $L^1(\Sigma)$ with respect to $\Vert\cdot\Vert_{W^{k,p}(\Sigma)}$ defined above.  With $p,q\geq 1$ and $s>0$, the \emph{Besov space} 
  \begin{equation}
    B^s_{p,q}(\Sigma)\vcentcolon= \begin{cases}
      (W^{k,p}(\Sigma),W^{k+1,p}(\Sigma))_{s-k,q}&\text{if $k<s<k+1$},\\
      (W^{k,p}(\Sigma),W^{k+2,p}(\Sigma))_{\frac12,q}&\text{if $s=k+1$}.
    \end{cases}
  \end{equation}
  Here $(\cdot,\cdot)_{\theta,r}$ denotes the real interpolation functor for $0<\theta<1$ and $r\geq 1$. The \emph{Sobolev--Slobodeckii} spaces are defined by $W^{s,p}(\Sigma)=B^s_{p,p}(\Sigma)$ for $0<s\notin\N$ and $1<p<\infty$. These function spaces are analogously defined on $\partial\Sigma$.
\end{definition}

As a direct consequence of \Cref{lem:unif-equivalent-norms}, the embedding and trace theorems of the Euclidean counterparts of the above function spaces generalize to the spaces on $\Sigma$ through the \emph{good atlas} in \Cref{hyp:unif-reg}, see \cite[Theorems~4.4, 4.5 and 4.7]{amann2025}.

\begin{remark}
  We point out a conceptual difference in the approach to defining the Sobolev-type spaces in \Cref{lem:unif-equivalent-norms} and \cite[Theorems~4.4, 4.5 and 4.7]{amann2025}. All Sobolev-norms in \cite{amann2025} are 
  equivalent (with uniform constants) not to the $\ell^\infty$ norm of the vector of localizations as used in \Cref{lem:unif-equivalent-norms}, but to the $\ell^p$-norm. However, this makes Sobolev and Hölder norms not uniformly comparable, i.e.\ one cannot uniformly estimate $\Vert u\Vert_{\tilde{\bm W}^{s,p}(\Sigma)}$ by $\Vert u\Vert_{\tilde{\bm C}^s(\Sigma)}$, without controlling the mass of $\Sigma$ (or equivalently, controlling $N_1$, $N_2$ and $Z^\ast$). Note that, since the $\ell^p$ and $\ell^\infty$ norms on $\R^{N_1+N_2}$ are equivalent, the function spaces are equal as sets and all non-quantitative properties immediately translate from \cite{amann2025} to our definition. In particular, we get $B^s_{p,p}(\Sigma)=W^{s,p}(\Sigma)$ and $B^s_{p,p}(\partial\Sigma)=W^{s,p}(\partial\Sigma)$ for $s>0$, $s\notin\N$, and $1\leq p<\infty$. 
  Nevertheless, we also obtain the following \emph{uniform} constants in the Sobolev, Morrey, trace and Gagliardo--Nirenberg inequalities as an immediate application of the Euclidean results.
\end{remark}

\begin{remark}[Unlocalized coordinate norms]\label{rem:unlocalized-coordinate-norms}
  Assume \Cref{hyp:unif-reg} and fix the localization system from \Cref{rem:local-sys}. For the norms and parameters in \Cref{lem:unif-equivalent-norms}, we have
  \begin{align}\label{eq:1-rem:unlocalized-coordinate-norms}
    \max_{k\in\mathbf I_{\mathrm{int}}\cup\mathbf I_{\mathrm{bdr}}}\Vert u\circ\alpha_k\Vert_{\mathcal F(D_k)}\leq C\Vert u\Vert_{\bm{\mathcal F}(\Sigma)},\quad\text{and}\quad\max_{k\in\mathbf I_{\mathrm{bdr}}}\Vert u\circ\alpha_k(\cdot,0)\Vert_{\mathcal F(I_{16\rho^\ast})}\leq C\Vert u\Vert_{\bm{\mathcal F}(\partial\Sigma)}
  \end{align}
  where $\mathcal F$ denotes any of the function spaces $C^s$, $W^{s,p}$, or $B^s_{p,q}$. The constants $C$ depend only on the respective function space parameters and the universal constants in \Cref{hyp:unif-reg}, independently of $N_1$, $N_2$. To see this, one uses that
  \begin{equation}
    u\circ\alpha_k=\sum_{j:\Omega_k\cap\operatorname{supp}\pi_j\neq\emptyset}(\pi_j^2u)\circ\alpha_k.
  \end{equation}
  Moreover, the fact that $\operatorname{dist}\big(\varphi_j(\operatorname{supp}\pi_j),\partial D_j\setminus(\R\times\{0\})\big)\geq\rho^\ast/4$ for all $j\in \mathbf I_{\mathrm{int}}\cup\mathbf I_{\mathrm{bdr}}$ and \Cref{hyp:unif-reg}(a) together yield
  \begin{equation}
    \#\{j:\Omega_k\cap\operatorname{supp}\pi_j\neq\emptyset\}\leq CZ^\ast \quad\text{for a universal constant $C>0$}. 
  \end{equation}
  Combined with \Cref{hyp:unif-reg}(b) and \eqref{eq:unif-bounds-pij}, we get \eqref{eq:1-rem:unlocalized-coordinate-norms}.
\end{remark}

\begin{proposition}\label{prop:embeddings-and-interpolation-results}
  \begin{enumerate}[(a)]
    \item\label{it:Sobolev-prop:embeddings-and-interpolation-results} 
      (Sobolev) For $1\leq p_0,p_1<\infty$ and $0\leq s_0<s_1$ such that $s_1-2/p_1\geq s_0-2/p_0$ (with a strict inequality unless $s_0,s_1\in\N$), there exists a constant $C$ only depending on $s_0,s_1,p_0,p_1$ and the universal constants in \Cref{hyp:unif-reg} such that
      \begin{equation}
        \Vert u\Vert_{\bm W^{s_0,p_0}(\Sigma)}\leq C \Vert u\Vert_{\bm W^{s_1,p_1}(\Sigma)}\quad\text{for all $u\in W^{s_1,p_1}(\Sigma)$}.
      \end{equation}
    \item\label{it:Morrey-prop:embeddings-and-interpolation-results} 
      (Morrey) For $1\leq p<\infty$ and $0\leq \alpha <s-2/p$, there exists a constant $C$ only depending on $s,\alpha,p$ and the universal constants in \Cref{hyp:unif-reg} such that 
      \begin{equation}
        \Vert u\Vert_{\bm C^\alpha(\Sigma)} \leq C \Vert u\Vert_{\bm W^{s,p}(\Sigma)}\quad\text{for all $u\in W^{s,p}(\Sigma)$}.
      \end{equation}
    \item\label{it:tracethm-prop:embeddings-and-interpolation-results} 
      (Trace theorem) Denote by $\bar\eta$ the unit outer conormal at the boundary of $(\Sigma,\bar g)$, with $\bar g$ from \Cref{rem:local-sys}. For $1<p<\infty$, $j\in\N_0$ and $s$ such that $s-j-\frac{1}{p}>0$ is not an integer, we have that $\partial_{\bar\eta}^j\colon W^{s,p}(\Sigma)\to W^{s-j-\frac1p,p}(\partial\Sigma)$ is a retraction.
    \item\label{it:gagliardonirenberg-prop:embeddings-and-interpolation-results}
      (Gagliardo--Nirenberg) Let $0\leq s_0<s<s_1$, $p,p_0,p_1\in [1,\infty)$ (where $p_i>1$ if $s_i\in\N$) and $\theta\in(0,1)$ such that 
      \begin{equation}
        s-\frac2p = (1-\theta)s_0+\theta s_1 - 2\Big(\frac{1-\theta}{p_0}+\frac\theta{p_1}\Big)\quad\text{and}\quad \frac1p \leq \frac{1-\theta}{p_0}+\frac\theta{p_1}.
      \end{equation}
      Then there exists a constant $C$ only depending on $s,s_0,s_1,p,p_0,p_1,\theta$ and the universal constants in \Cref{hyp:unif-reg} such that
      \begin{equation}
        \Vert u\Vert_{\bm W^{s,p}(\Sigma)} \leq C \Vert u\Vert_{\bm W^{s_0,p_0}(\Sigma)}^{1-\theta} \Vert u\Vert_{\bm W^{s_1,p_1}(\Sigma)}^\theta\quad\text{for all $u\in W^{s_0,p_0}(\Sigma)\cap W^{s_1,p_1}(\Sigma)$}.
      \end{equation}
  \end{enumerate}
\end{proposition}

\begin{lemma}[Whitney-type extension]\label{lem:whitney-extension}
  Fix $0<\alpha<1$ and an even nonnegative $\zeta\in C_c^\infty((-1,1))$ with $\int_\R\zeta=1$. For $u\in C^\alpha(\R,\R^N)$, define
  \begin{equation}\label{eq:def-whitney-extension}
    \mathcal Eu(x,y)\vcentcolon=y\int_\R\zeta(t)u(x-yt)\dd t\quad\text{for $x\in\R$ and $y\geq0$}.
  \end{equation}
  Then $\mathcal Eu\in C^{1+\alpha}(\R\times[0,2],\R^N)$, and
  \begin{equation}\label{eq:whitney-extension-traces}
    \mathcal Eu(x,0)=0,\quad \partial_x\mathcal Eu(x,0)=0\quad\text{and}\quad \partial_y\mathcal Eu(x,0)=u(x).
  \end{equation}
  Moreover, we have
  \begin{equation}\label{eq:whitney-extension-estimates}
    \Vert\mathcal Eu\Vert_{C^{1+\alpha}(\R\times[0,2])}\leq C\Vert u\Vert_{C^\alpha(\R)}\quad\text{and}\quad \Vert\mathcal Eu\Vert_{C^l(\R\times[0,2])}\leq C_l\Vert u\Vert_{C^{l-1}(\R)}\quad\text{for $l\geq2$ and smooth $u$}.
  \end{equation}
  Here, $C$ depends at most on $\alpha$ and $\zeta$, and $C_l$ depends at most on $l$ and $\zeta$. Smooth data yield an extension which is smooth up to $y=0$. Assuming \Cref{hyp:unif-reg}, there is also a linear operator $\mathcal E_{\partial\Sigma}\colon C^\alpha(\partial\Sigma,\R^N)\to C^{1+\alpha}(\partial\Sigma\times[0,1],\R^N)$ such that
  \begin{equation}\label{eq:boundary-whitney-extension-traces}
    \mathcal E_{\partial\Sigma}u(p,0)=0\quad\text{and}\quad \partial_y\mathcal E_{\partial\Sigma}u(p,0)=u(p),
  \end{equation}
  with the analog estimates to \eqref{eq:whitney-extension-estimates}. The norms on the cylinder are defined by the boundary charts and the localization system as in \Cref{lem:unif-equivalent-norms}, and the constants additionally depend on the universal constants in \Cref{hyp:unif-reg}.
\end{lemma}
\begin{proof}[Proof sketch]
  The first part of the statement is standard, so we only comment on the second part on $\mathcal E_{\partial\Sigma}$. Fix $\chi_0\in C_c^\infty(I_{16\rho^\ast})$ equal to one on $I_{13\rho^\ast}$. Extend $u_j(x)\vcentcolon=\chi_0(x)u(\alpha_j(x,0))$ by zero to $\R$, and set
  \begin{equation}\label{eq:def-boundary-whitney-extension}
    \mathcal E_{\partial\Sigma}u(p,y)\vcentcolon=\sum_{j\in\mathbf I_{\mathrm{bdr}}}\pi_j^2(p)\mathcal Eu_j(\varphi_j^{(1)}(p),y).
  \end{equation}
  Since $\operatorname{supp}\pi_j\cap\partial\Sigma\subset\alpha_j(I_{25\rho^\ast/2}\times\{0\})$ and $\sum_{j\in\mathbf I_{\mathrm{bdr}}}\pi_j^2=1$ on $\partial\Sigma$, \eqref{eq:whitney-extension-traces} yields \eqref{eq:boundary-whitney-extension-traces}. The estimates follow from \eqref{eq:whitney-extension-estimates}, \eqref{eq:unif-bounds-pij} and the bounded overlap and transition estimates in \Cref{hyp:unif-reg}.
\end{proof}

\subsection{Weighted parabolic Sobolev spaces and maximal regularity}
\label{app:max-reg-unif-surfaces}

For $1<p<\infty$, $\mu\in(\frac1p,1]$, $T>0$ and $J_T=[0,T]$ and some Banach space $X$, we define the weighted Lebesgue and Sobolev spaces
\begin{align}
  L^p_\mu(J_T;X)&=\{u\in L^1_{\mathrm{loc}}((0,T);X)\mid t^{1-\mu}u\in L^p(J_T;X)\},\\
  \Vert u\Vert_{L^p_\mu(J_T;X)}&=\Vert t^{1-\mu}u\Vert_{L^p(J_T;X)},\\
  W^{1,p}_\mu(J_T;X)&=\{u\in W^{1,1}_{\mathrm{loc}}((0,T);X)\mid u,\partial_tu\in L^p_\mu(J_T;X)\},\\
  \Vert u\Vert_{W^{1,p}_\mu(J_T;X)}&=\Vert u\Vert_{L^p_\mu(J_T;X)}+\Vert\partial_tu\Vert_{L^p_\mu(J_T;X)}. \label{eq:parabolic-norms}
\end{align}
These are Banach spaces equipped with their respective norms. For functions on $\Sigma$, we use the following localization convention throughout the manuscript.

We assume \Cref{hyp:unif-reg}, fixing the \emph{good atlas} therein, and in particular work with \Cref{lem:unif-equivalent-norms}, and the localization system \Cref{rem:local-sys}. Set
\begin{equation}
  \mathbb X_j=\mathbb R^2\quad\text{for }j\in \mathbf I_{\mathrm{int}}\quad\text{and}\quad\mathbb X_j=\mathbb R\times[0,\infty)\quad\text{for }j\in \mathbf I_{\mathrm{bdr}}.
\end{equation}
For a function $u$ on $J_T\times\Sigma$, put $u_j(t,x)=\pi_j(\alpha_j(x)) u(t,\alpha_j(x))$. We identify $u_j$ with its extension by zero from $D_j$ to $\mathbb X_j$. For $r\in\{0,1\}$ and $s\geq0$, define
\begin{equation}\label{eq:localized-parabolic-norms}
  \Vert u\Vert_{\bm W^{r,p}_\mu(J_T;\bm W^{s,p}(\Sigma))}\vcentcolon=\max_{j\in \mathbf I_{\mathrm{int}}\cup\mathbf I_{\mathrm{bdr}}}\Vert u_j\Vert_{W^{r,p}_\mu(J_T;W^{s,p}(\mathbb X_j))}\quad\text{for $u\in W^{r,p}_\mu(J_T;W^{s,p}(\Sigma))$},
\end{equation}
where $W^{0,p}_\mu=L^p_\mu$. It is easy to see that the norms in \eqref{eq:localized-parabolic-norms} are equivalent to those from \eqref{eq:parabolic-norms}.
\begin{remark}\label{rem:prop-parabolic-norms}
  The spaces and (localized) norms defined above have the following properties.
  \begin{enumerate}[(a)]
    \item \label{it:1-rem:prop-parabolic-norms}
      Arguing as in \Cref{rem:unlocalized-coordinate-norms}, one finds for $u\in W^{r,p}_\mu(J_T;W^{s,p}(\Sigma))$ and $k\in\mathbf{I}_{\mathrm{int}}\cup\mathbf{I}_{\mathrm{bdr}}$ that
      \begin{equation}
        \Vert u(t,\alpha_k(x))\Vert_{W^{r,p}_\mu(J_T;W^{s,p}(D_k))} \leq C \Vert u\Vert_{\bm W^{r,p}_\mu(J_T;\bm W^{s,p}(\Sigma))}
      \end{equation}
      where $C>0$ is a constant only depending on the universal constants in \Cref{hyp:unif-reg} (not $N_1$ and $N_2$), and on the function space parameters.
    \item  \label{it:2-rem:prop-parabolic-norms}
      Using that $\mu p>1$ and Hölder's inequality, we have
      \begin{equation}
        \max_{j\in \mathbf I_{\mathrm{int}}\cup\mathbf I_{\mathrm{bdr}}}\Vert u_j\Vert_{L^1(J_T;W^{s,p}(\mathbb X_j))}\leq\Big(\frac{p-1}{\mu p-1}\Big)^{\frac{p-1}{p}}T^{\mu-\frac1p}\Vert u\Vert_{\bm L^p_\mu(J_T;\bm W^{s,p}(\Sigma))}.
      \end{equation}
      In particular, for $u\in W^{1,p}_\mu(J_T;L^p(\Sigma))$, there exists a representative $u\in C^0(J_T;L^p(\Sigma))$. In \eqref{eq:localized-parabolic-norms}, using $u_j(t)=u_j(0)+\int_0^t\partial_tu_j(r)\dd r$, we obtain
      \begin{equation}\label{eq:localized-small-time-Lp}
        \Vert u\Vert_{\bm L^p_\mu(J_T;\bm L^p(\Sigma))}\leq C_{p,\mu}\big(T^{1-\mu+\frac1p}\Vert u(0)\Vert_{\bm L^p(\Sigma)}+T\Vert\partial_tu\Vert_{\bm L^p_\mu(J_T;\bm L^p(\Sigma))}\big).
      \end{equation}
  \end{enumerate}
\end{remark}
Define
\begin{equation}
  X_0=L^p(\Sigma)\quad\text{and}\quad X_1=\{u\in W^{4,p}(\Sigma)\mid u=du=0\text{ on }\partial\Sigma\}.
\end{equation}
With \eqref{eq:localized-parabolic-norms}, define
\begin{align}
  \mathbb E_{1,\mu}(J_T)&=W^{1,p}_\mu(J_T;X_0)\cap L^p_\mu(J_T;X_1),\quad \Vert u\Vert_{\mathbb E_{1,\mu}(J_T)}=\Vert\partial_tu\Vert_{\bm L^p_\mu(J_T;\bm L^p(\Sigma))}+\Vert u\Vert_{\bm L^p_\mu(J_T;\bm W^{4,p}(\Sigma))}.
\end{align}
For the trace spaces, set $X_{0,j}=L^p(\mathbb X_j)$ and
\begin{equation}
  X_{1,j}=\begin{cases}
    W^{4,p}(\mathbb R^2)&\text{if }j\in \mathbf I_{\mathrm{int}},\\
    \{v\in W^{4,p}(\mathbb X_j)\mid v=\nabla v=0\text{ on }\mathbb R\times\{0\}\}&\text{if }j\in \mathbf I_{\mathrm{bdr}}.
  \end{cases}
\end{equation}
Let $\mathbb E_{1,\mu}^j(J_T)=W^{1,p}_\mu(J_T;X_{0,j})\cap L^p_\mu(J_T;X_{1,j})$ with $\Vert u\Vert_{\mathbb E_{1,\mu}^j(J_T)}=\Vert\partial_tu\Vert_{L^p_\mu(J_T;L^p(\mathbb X_j))}+\Vert u\Vert_{\bm L^p_\mu(J_T;\bm W^{4,p}(\mathbb X_j))}$, and write
\begin{equation}
  \mathcal X_{\mu,j}=(X_{0,j},X_{1,j})_{\mu-\frac1p,p}\quad\text{and}\quad\mathcal X_\mu(\Sigma)=\{u\in L^p(\Sigma)\mid(\pi_j u)\circ\alpha_j\in\mathcal X_{\mu,j}\text{ for all }j\in \mathbf I_{\mathrm{int}}\cup\mathbf I_{\mathrm{bdr}}\},
\end{equation}
equipped with
\begin{equation}\label{eq:localized-initial-trace-norm}
  \Vert u\Vert_{\bm{\mathcal X}_\mu(\Sigma)}=\max_{j\in \mathbf I_{\mathrm{int}}\cup\mathbf I_{\mathrm{bdr}}}\Vert(\pi_j u)\circ\alpha_j\Vert_{\mathcal X_{\mu,j}}.
\end{equation}
In particular, if $4(\mu-1/p)\notin\mathbb N$ and $4(\mu-1/p)>1+1/p$, the classical Euclidean trace-space identification for each $\mathbb E_{1,\mu}^j(J_T)$ yields
\begin{equation}\label{eq:localized-trace-identification}
  \mathcal X_\mu(\Sigma)=\{u\in W^{4(\mu-1/p),p}(\Sigma)\mid u=du=0\text{ on }\partial\Sigma\},
\end{equation}
with uniform equivalence to $\Vert\cdot\Vert_{\bm W^{4(\mu-1/p),p}(\Sigma)}$ as in \Cref{lem:unif-equivalent-norms}.
\begin{proposition}\label{prop:BUC-estimate}
  Let $0<T_0\leq T$. Then $\mathbb E_{1,\mu}(J_T)\hookrightarrow C^0(J_T,\mathcal X_\mu(\Sigma))$, and there exists a constant $C>0$ only depending on $\mu$, $p$ and $T_0$ such that
  \begin{equation}
    \max_{t\in J_T}\Vert u(t)\Vert_{\bm{\mathcal X_\mu}(\Sigma)}\leq C\Vert u\Vert_{\mathbb E_{1,\mu}(J_T)}\quad\text{for all }u\in\mathbb E_{1,\mu}(J_T).
  \end{equation}
  If $u(0)=0$, the estimate holds with a constant only depending on $\mu$ and $p$, independently of $T$. Further, $u\in C^0((0,T],\mathcal X_1(\Sigma))$.
\end{proposition}
\begin{proof}[Proof sketch]
  Apply the Euclidean result (see \cite[Proposition~3.1(b)]{pruesssimonett2004}) to $u_j=(\pi_j u)\circ\alpha_j$ with the pair $(X_{0,j},X_{1,j})$ on each $\mathbb X_j$. Taking the maximum over $j$ yields the result since the  maximum over $t$ and the maximum over $j$ commute. Continuity follows from the corresponding local continuity and finiteness of the atlas. The last assertion follows similarly, using \cite[Proposition~3.1(c)]{pruesssimonett2004}.
\end{proof}

Suppose that $a_j\in C^\alpha(D_j,\mathbb R^{2\times2})$ and $c_{j,\nu}\in C^\alpha(D_j)$ for $|\nu|\leq3$ satisfy
\begin{equation}
  \Vert(a_j^{ab})-(\delta^{ab})\Vert_{C^\alpha(D_j)}<\frac12\quad\text{and}\quad\Vert c_{j,\nu}\Vert_{C^\alpha(D_j)}\leq M
\end{equation}
for all $j\in \mathbf I_{\mathrm{int}}\cup\mathbf I_{\mathrm{bdr}}$. Consider an operator $A\colon X_1\to X_0$ given in local coordinates by
\begin{equation}\label{eq:linear-operator-in-local-coords}
  (A.u)\circ\alpha_j=a_j^{ab}a_j^{cd}\partial_{abcd}^4(u\circ\alpha_j)+\sum_{|\nu|\leq3}c_{j,\nu}\partial^\nu(u\circ\alpha_j)\quad\text{on }D_j.
\end{equation}

\begin{lemma}[Maximal regularity]\label{prop:max-L^p-reg}
  There exists a constant $C_{\mathrm{MR}}>0$ only depending on the universal constants in \Cref{hyp:unif-reg}, $\alpha$, $M$, $p>2$, $\mu\in(\frac1p,1]$ and an upper bound for $T>0$ such that the following holds. For every $f\in L^p_\mu(J_T;L^p(\Sigma))$ and $u_0\in\mathcal X_\mu(\Sigma)$, there exists a unique solution $u\in\mathbb E_{1,\mu}(J_T)$ of
  \begin{equation}\label{eq:linear-max-reg-problem}
    \begin{cases}
      \partial_tu+A.u=f&\text{a.e.\ in }(0,T)\times\Sigma,\\
      u(0)=u_0&\text{on }\Sigma,\\
      u=du=0&\text{on }(0,T)\times\partial\Sigma,
    \end{cases}
  \end{equation}
  and
  \begin{equation}\label{eq:linear-max-reg-estimate}
    \Vert u\Vert_{\mathbb E_{1,\mu}(J_T)}\leq C_{\mathrm{MR}}\big(\Vert f\Vert_{\bm L^p_\mu(J_T;\bm L^p(\Sigma))}+\Vert u_0\Vert_{\bm{\mathcal X}_\mu(\Sigma)}\big).
  \end{equation}
  All parabolic norms are understood in the sense of \eqref{eq:localized-parabolic-norms}. In the range of \eqref{eq:localized-trace-identification}, the initial-data norm can equivalently be replaced by $\Vert u_0\Vert_{\bm W^{4(\mu-1/p),p}(\Sigma)}$.
\end{lemma}
\begin{proof}
  We first establish the a-priori estimate by means of localization via $\pi_j$ and the corresponding Euclidean maximal regularity results in \cite[Section~4]{denkhieberpruess2007} and \cite[Theorem~3.2]{pruesssimonett2004}. Recall the above definition of $\mathbb E_{1,\mu}^j(J_T)$. For each $j$, choose a $C^\alpha$-extension $\widetilde a_j$ to $\mathbb X_j$ agreeing with $a_j$ on a neighborhood of $\operatorname{supp}(\pi_j\circ\alpha_j)$, with uniform $C^\alpha$-bounds and $\Vert\widetilde a_j-\operatorname{id}\Vert_{C^0}<1/2$. Set $\mathcal A_jv=\widetilde a_j^{ab}\widetilde a_j^{cd}\partial_{abcd}^4v$ with domain $X_{1,j}$. These operators are uniformly normally elliptic, and the clamped boundary operators satisfy the Lopatinskii--Shapiro condition. The Euclidean results yield
  \begin{align}
    &\qquad\partial_t+\mathcal A_j \colon \mathbb E_{1,\mu}^{j}(J_T)\to L^p_\mu(J_T;L^p(\mathbb X_j))\times\mathcal X_{\mu,j} \text{ is an isomorphism with}\\
    &\Vert v\Vert_{\mathbb E_{1,\mu}^{j}(J_T)}\leq C_{\mathrm{MR}}\big(\Vert \partial_tv+\mathcal A_jv\Vert_{L^p_\mu(J_T;L^p(\mathbb X_j))}+\Vert v(0)\Vert_{\mathcal X_{\mu,j}}\big)\quad\text{for all }v\in \mathbb E_{1,\mu}^{j}(J_T).\label{eq:euclidean-max-reg-estimate}
  \end{align}
  Let $u$ be a solution of \eqref{eq:linear-max-reg-problem} and set $u_j=(\pi_j u)\circ\alpha_j$ and $f_j=(\pi_j f)\circ\alpha_j$. Then
  \begin{equation}
    \partial_tu_j+\mathcal A_ju_j=f_j+R_j u\quad\text{and}\quad u_j(0)=(\pi_j u_0)\circ\alpha_j,
  \end{equation}
  and $u_j\in \mathbb E_{1,\mu}^{j}(J_T)$. Here, $R_j$ is made up of the commutator with $\pi_j$ and the lower-order terms of $A$, and using \Cref{rem:prop-parabolic-norms}\eqref{it:1-rem:prop-parabolic-norms}, one finds
  \begin{equation}
    \Vert R_j u\Vert_{L^p_\mu(J_T;L^p(\mathbb X_j))} \leq C_M\Vert u\Vert_{\bm L^p_\mu(J_T;\bm W^{3,p}(\Sigma))}.
  \end{equation}
  Applying \eqref{eq:euclidean-max-reg-estimate}, we obtain
  \begin{equation}\label{eq:localized-max-reg-commutator}
    \Vert u\Vert_{\mathbb E_{1,\mu}(J_T)}\leq C\big(\Vert u_0\Vert_{\bm{\mathcal X}_\mu(\Sigma)}+\Vert f\Vert_{\bm L^p_\mu(J_T;\bm L^p(\Sigma))}+C_M\Vert u\Vert_{\bm L^p_\mu(J_T;\bm W^{3,p}(\Sigma))}\big).
  \end{equation}
  Using the Euclidean Gagliardo--Nirenberg inequality and Hölder's inequality with respect to $t^{p(1-\mu)}\dd t$, we have for each $j$ that
  \begin{equation}\label{eq:interpolation-estimate}
    \Vert u_j\Vert_{L^p_\mu(J_T;W^{3,p}(\mathbb X_j))}\leq C\Vert u_j\Vert_{L^p_\mu(J_T;W^{4,p}(\mathbb X_j))}^{\frac34}\Vert u_j\Vert_{L^p_\mu(J_T;L^p(\mathbb X_j))}^{\frac14}.
  \end{equation}
  Taking the maximum over $j$, applying Young's inequality and absorbing the fourth-order term in \eqref{eq:localized-max-reg-commutator}, we get
  \begin{equation}
    \Vert u\Vert_{\mathbb E_{1,\mu}(J_T)}\leq C\big(\Vert u_0\Vert_{\bm{\mathcal X}_\mu(\Sigma)}+\Vert f\Vert_{\bm L^p_\mu(J_T;\bm L^p(\Sigma))}\big)+C_M\Vert u\Vert_{\bm L^p_\mu(J_T;\bm L^p(\Sigma))}.
  \end{equation}
  Using \eqref{eq:localized-small-time-Lp}, we can choose $T_0\in(0,1]$ such that the last time-derivative term can be absorbed for all $T\leq T_0$. Since $\mathcal X_\mu(\Sigma)\hookrightarrow L^p(\Sigma)$ uniformly, this proves \eqref{eq:linear-max-reg-estimate} for $T\leq T_0$, with a constant independent of $T$. In particular, uniqueness holds on these intervals.

  Existence on a sufficiently short interval can also be obtained using \eqref{eq:euclidean-max-reg-estimate}. For $g\in L^p_\mu(J_T;L^p(\Sigma))$, using \eqref{eq:euclidean-max-reg-estimate}, let $v_j\in\mathbb E_{1,\mu}^j(J_T)$ solve
  \begin{equation}
    \partial_tv_j+\mathcal A_jv_j=(\pi_jg)\circ\alpha_j\quad\text{and}\quad v_j(0)=0
  \end{equation}
  and put $Qg=\sum_{j\in\mathbf{I}_{\mathrm{int}}\cup\mathbf{I}_{\mathrm{bdr}}}\pi_j(v_j\circ\varphi_j)$. Using \eqref{eq:euclidean-max-reg-estimate} and \Cref{rem:prop-parabolic-norms}\eqref{it:1-rem:prop-parabolic-norms}, we have $\Vert Qg\Vert_{\mathbb E_{1,\mu}(J_T)}\leq C\Vert g\Vert_{\bm L^p_\mu(J_T;\bm L^p(\Sigma))}$. Furthermore, $Qg$ solves
  \begin{equation}
    (\partial_t+A)Qg=g+Kg,
  \end{equation}
  where $Kg$ is a sum of localized terms of order at most three in $v_j$. Since $v_j(0)=0$, \eqref{eq:localized-small-time-Lp} on $\mathbb X_j$ and \eqref{eq:interpolation-estimate} yield
  \begin{equation}
    \Vert v_j\Vert_{L^p_\mu(J_T;W^{3,p}(\mathbb X_j))}\leq CT^{\frac14}\Vert v_j\Vert_{\mathbb E_{1,\mu}^{j}(J_T)}.
  \end{equation}
  Thus, $\Vert Kg\Vert_{\bm L^p_\mu(J_T;\bm L^p(\Sigma))}\leq C_MT^{1/4}\Vert g\Vert_{L^p_\mu(J_T;L^p(\Sigma))}$. After reducing $T_0$ if necessary, $\operatorname{id}+K$ is invertible by a Neumann series argument, and $Q(\operatorname{id}+K)^{-1}$ solves the problem with zero initial data.

  For general $u_0\in\mathcal X_\mu(\Sigma)$, let $z_j$ solve $\partial_tz_j+\mathcal A_jz_j=0$ with initial data $(\pi_j u_0)\circ\alpha_j$, and put $z=\sum_{j\in\mathbf{I}_{\mathrm{int}}\cup\mathbf{I}_{\mathrm{bdr}}}\pi_j(z_j\circ\varphi_j)$. Then $z(0)=u_0$ and, using \eqref{eq:euclidean-max-reg-estimate},
  \begin{equation}
    \Vert z\Vert_{\mathbb E_{1,\mu}(J_T)}+\Vert(\partial_t+A)z\Vert_{\bm L^p_\mu(J_T;\bm L^p(\Sigma))}\leq C\Vert u_0\Vert_{\bm{\mathcal X}_\mu(\Sigma)}.
  \end{equation}
  Applying the previous argument with right-hand-side $f-(\partial_t+A)z$, and adding $z$ proves existence for $T\leq T_0$.

  Finally, choose $T_0$ such that the short-time result holds for both $\mu$ and $1$. By \Cref{prop:BUC-estimate}, applied with $\mu=1$ on $[T_0/2,T_0]$, we have
  \begin{align}
    \Vert u(T_0)\Vert_{\mathcal X_1(\Sigma)}&\leq C(T_0)\Vert u\Vert_{\mathbb E_{1,1}([T_0/2,T_0])}\leq C(T_0,\mu,p)\Vert u\Vert_{\mathbb E_{1,\mu}(J_{T_0})}.
  \end{align}
  We can therefore continue the solution on intervals of length $T_0$ using the unweighted result. On these intervals, $t\mapsto t^{1-\mu}$ is bounded above and below by constants depending only on $T_0$, $\mu$ and the prescribed upper bound for $T$. The number of intervals is bounded in terms of the same quantities. Applying the estimates on each interval and concatenating them proves \eqref{eq:linear-max-reg-estimate} and uniqueness for every $T>0$. To obtain constants depending only on an upper bound for $T$, extend $f$ by zero to the end of the last full interval and then restrict the resulting solution to $J_T$. 
\end{proof}

\subsection{Geometric uniqueness for the Dirichlet Willmore flow}
\label{app:geometric-uniqueness-Dirichlet}

We state the geometric uniqueness statement used in \Cref{sec:blow-up-construction}. A detailed proof goes beyond the scope of this article. 

\begin{proposition}[Geometric uniqueness]\label{prop:geometric-uniqueness-Dirichlet}
  Let $\Sigma$ be a compact oriented surface with boundary, fix the parameters $p$ and $\mu$ as in \Cref{sec:lifespan-estimate}, and let $\gamma\colon\partial\Sigma\to\Gamma$ be a smooth embedding onto a smooth compact one-dimensional submanifold, with a smooth unit normal field $\mathbf N$ along $\Gamma$. Suppose that, for $i\in\{1,2\}$,
  \begin{equation}
    f_i\in C^0([0,T_i),W^{4(\mu-1/p),p}_{\mathrm{imm}}(\Sigma,\R^3))\cap C^\infty((0,T_i)\times\Sigma,\R^3)
  \end{equation}
  are solutions of
  \begin{equation}
    \begin{cases}
      \langle\partial_t f_i,\mathbf n_{f_i}\rangle=-\frac12\big(\Delta_{g_{f_i}}\mathbf H_{f_i}+|\mathbf A_{f_i}^0|^2\mathbf H_{f_i}\big)&\text{on $(0,T_i)\times\Sigma$}\\
      f_i(t,\cdot)\vert_{\partial\Sigma}\colon\partial\Sigma\to\Gamma\text{is a diffeomorphism, and }\partial_{\eta_{f_i}}f_i=\mathbf N\circ f_i&\text{on $(t,x)\in [0,T_i)\times\partial\Sigma$},
    \end{cases}
  \end{equation}
  and $f_1(0,\Xi_0(x))=f_2(0,x)$ for all $x\in\Sigma$ and for a diffeomorphism $\Xi_0\in W^{4(\mu-1/p),p}(\Sigma,\Sigma)$. Then there exists a family of diffeomorphisms
  \begin{equation}
    \Xi\in C^0([0,\min\{T_1,T_2\}),W^{4(\mu-1/p),p}(\Sigma,\Sigma))\cap C^\infty((0,\min\{T_1,T_2\})\times\Sigma,\Sigma)
  \end{equation}
  with
  \begin{align}
      &\hspace{2em}\Xi(0)=\Xi_0,\quad\Xi(t,\partial\Sigma)=\partial\Sigma,\text{ and}\\
      &f_1(t,\Xi(t,x))=f_2(t,x)\quad\text{for all $0\leq t<\min\{T_1,T_2\}$ and $x\in\Sigma$}.
  \end{align}
  If $f_2$ is a maximal solution with $\partial_tf_2-\langle \partial_tf_2,\mathbf n_{f_2}\rangle\mathbf n_{f_2}\equiv 0$, then $T_2\geq T_1$.
\end{proposition}

To prove this statement, one uses a common smooth reference immersion, the DeTurck trick with Dirichlet boundary conditions, and uniqueness and continuous dependence on the data for the associated quasilinear PDE, in the setting of weighted maximal $L^p$-regularity, see \Cref{prop:max-L^p-reg} and \cite{lecronepruesswilke2014}. Assuming rotational symmetry, more details are given in \cite[Appendix~C]{schlierf2024na}, and for analogous geometric results for elastic flows, cf.\ \cite[Theorems~5.4 and~5.5]{garckemenzelpluda2020}, \cite[Appendix~C]{dallacqualinpozzi2021}, and \cite[Proposition~2.16, Remark~3.15 and Theorem~3.18]{mantegazzapludapozzetta2021}. 

\section{Simon's monotonicity formula for immersions with boundary}
\label{app:monotonicity-formula}

The following proposition contains the version of Simon's monotonicity formula going back to \cite{simon1993} when adapted to our setting of immersed surfaces with (Dirichlet) boundary. For Dirichlet boundary conditions, it is already used in \cite{novagapozzetta2020,schlierf2025acv,ruppscharrer2025}, and \cite{volkmann2016} adapts it to surfaces with free boundary on the unit sphere. Also see \cite{ruppscharrer2023} for an adaptation to the Canham--Helfrich energy.

\begin{proposition}[Simon's monotonicity formula with boundary]\label{lem:monotonicity-with-boundary}
  Let $f\colon\Sigma\to\R^3$ be a smooth proper immersion of an oriented surface with possibly empty boundary, with $\W(f)<\infty$ and quadratic area growth, i.e.\ $\sup_{R>0}R^{-2}\int_{f^{-1}(B_R(0))}\dd\vol_{g_f}<\infty$. Fix $z\in\R^3$ such that
  \begin{equation}\label{eq:integrable-boundary-potential}
    \int_{\partial\Sigma}\frac{|\langle\partial_{\eta_f}f,f-z\rangle|}{|f-z|^2}\dd s<\infty.
  \end{equation}
  For $r>0$, set
  \begin{align}
    \mathcal M_f(r,z)&\vcentcolon=\frac{\int_{f^{-1}(B_r(z))}\dd\vol_{g_f}}{r^2}+\frac1{16}\int_{f^{-1}(B_r(z))}\mathbf H_f^2\dd\vol_{g_f}\\
    &\quad+\frac1{2r^2}\int_{f^{-1}(B_r(z))}\mathbf H_f\langle\mathbf n_f,f-z\rangle\dd\vol_{g_f}\\
    &\quad+\frac12\int_{\partial\Sigma\cap f^{-1}(B_r(z))}\Big(\frac1{|f-z|^2}-\frac1{r^2}\Big)\langle\partial_{\eta_f}f,f-z\rangle\dd s.
  \end{align}
  Then, for all $0<r<R$, we have
  \begin{equation}\label{eq:monotonicity-with-boundary}
    \mathcal M_f(r,z)+\int_{f^{-1}(B_R(z)\setminus B_r(z))}\Big|\frac14\mathbf H_f\mathbf n_f+\frac{(f-z)^\perp}{|f-z|^2}\Big|^2\dd\vol_{g_f}=\mathcal M_f(R,z).
  \end{equation}
  Here, $(f-z)^\perp=\langle f-z,\mathbf n_f\rangle\mathbf n_f$. The densities
  \begin{equation}
      \theta^2(f,z)\vcentcolon=\lim_{r\searrow0}\frac{\int_{f^{-1}(B_r(z))}\dd\vol_{g_f}}{\pi r^2},\quad \theta^2(f,\infty)\vcentcolon=\lim_{R\to\infty}\frac{\int_{f^{-1}(B_R(z))}\dd\vol_{g_f}}{\pi R^2}
  \end{equation}
  exist and are finite, and $\theta^2(f,\infty)$ does not depend on $z$. Moreover,
  \begin{align}
    \lim_{r\searrow0}\mathcal M_f(r,z)&=\pi\theta^2(f,z)\text{,}\quad
    \lim_{R\to\infty}\mathcal M_f(R,z)=\pi\theta^2(f,\infty)+\frac{\W(f)}4 +\frac12\int_{\partial\Sigma}\frac{\langle\partial_{\eta_f}f,f-z\rangle}{|f-z|^2}\dd s.\label{eq:limits-monotonicity-with-boundary}
  \end{align}
  Consequently, letting $R\to\infty$ and $r\searrow0$ in \eqref{eq:monotonicity-with-boundary}, we obtain
  \begin{align}
    &\pi\theta^2(f,z)+\int_{\Sigma\setminus f^{-1}(\{z\})}\Big|\frac14\mathbf H_f\mathbf n_f+\frac{(f-z)^\perp}{|f-z|^2}\Big|^2\dd\vol_{g_f}\\
    &\qquad=\pi\theta^2(f,\infty)+\frac{\W(f)}4+\frac12\int_{\partial\Sigma}\frac{\langle\partial_{\eta_f}f,f-z\rangle}{|f-z|^2}\dd s.\label{eq:global-monotonicity-with-boundary}
  \end{align}
\end{proposition}
\begin{proof}
  For \eqref{eq:monotonicity-with-boundary}, see \cite[Section~2 and Appendix~B]{novagapozzetta2020}, and \cite[Equation~(C.1)]{schlierf2025acv} in the finite-area case and \cite[Equations~(A.3), (A.4)]{kuwertschaetzle2004} in the case without boundary. In \eqref{eq:limits-monotonicity-with-boundary}, the limit $r\searrow 0$ follows as in \cite[Lemma~2.1]{novagapozzetta2020}. For the limit $R\to\infty$, the mixed term tends to zero by the estimate in \cite[Equation~(A.13)]{kuwertschaetzle2004}, which only uses finite Willmore energy and quadratic area growth. The boundary term converges by \eqref{eq:integrable-boundary-potential}, and thus the density at infinity exists as in \cite[Equation~(A.14)]{kuwertschaetzle2004}. 
\end{proof}

\begin{remark}[Lower semi-continuity of $\mathcal M$]\label{rem:lsc-monotonicity-formula}
  Let $(f_j,p_j)\to(\hat f,\hat p)$ in the sense of Definition~\ref{def:local-smooth-pointed-convergence}, and suppose that $\hat f$ and all $f_j$ satisfy the hypotheses of Proposition~\ref{lem:monotonicity-with-boundary} for some $z\in\R^3\setminus\hat f(\hat\Sigma)$. Then, using \eqref{eq:monotonicity-with-boundary} and \eqref{eq:limits-monotonicity-with-boundary}, we get for every $R>0$,
  \begin{align}
    \mathcal M_{\hat f}(R,z)
    &=
    \int_{\hat f^{-1}(B_R(z))}
    \Big|\frac14\mathbf H_{\hat f}\mathbf n_{\hat f}
    +\frac{(\hat f-z)^\perp}{|\hat f-z|^2}\Big|^2
    \dd\vol_{g_{\hat f}}\\
    &\leq\liminf_{j\to\infty}
    \int_{\Omega_j\cap f_j^{-1}(B_R(z)\setminus\{z\})}
    \Big|\frac14\mathbf H_{f_j}\mathbf n_{f_j}
    +\frac{(f_j-z)^\perp}{|f_j-z|^2}\Big|^2
    \dd\vol_{g_{f_j}}\\
    &\leq
    \liminf_{j\to\infty}\mathcal M_{f_j}(R,z).
  \end{align}
\end{remark}

In the following, fix a smooth compact one-dimensional submanifold $\Gamma\subset\R^3$ and a smooth unit normal field $\mathbf N$ along $\Gamma$. Let $\Sigma$ be a compact oriented surface with boundary, and let $\gamma\colon\partial\Sigma\to\Gamma$ be a smooth diffeomorphism. For immersions $f\colon\Sigma\to\R^3$ with $f\vert_{\partial\Sigma}=\gamma$ and $\partial_{\eta_f}f=\mathbf N\circ\gamma$ on $\partial\Sigma$, the boundary term in \eqref{eq:global-monotonicity-with-boundary} is given by \eqref{eq:def-Boundary-integral}. Moreover, for the proof of \eqref{eq:density-bound-stationary-blow-up} in \Cref{prop:existence-of-blow-ups}, the following quantization of the discontinuity of $\mathbf B_{\Gamma,\mathbf N}$ at $\Gamma$ is crucial.
\begin{lemma}\label{lem:boundary-integral-jump}
    Suppose that $q_j\in\Gamma$, $q_j\to q\in\Gamma$, $0<d_j\to 0$, and fix unit vectors $v_j\in (T_{q_j}\Gamma)^\bot$ converging to $v\in(T_q\Gamma)^\bot$. Then we have that
    \begin{equation}
      \lim_{j\to\infty}\mathbf B_{\Gamma,\mathbf N}(q_j+d_jv_j)=\mathbf B_{\Gamma,\mathbf N}(q)-\pi\langle\mathbf N(q),v\rangle.
    \end{equation}
\end{lemma}
\begin{proof}
  For $j$ sufficiently large, the denominator in \eqref{eq:def-Boundary-integral} only is singular in a neighborhood $\{p\in\Gamma\mid d_{\Gamma}(p,q)<\delta\}$ for some $\delta>0$ suitably small. By dominated convergence, we may thus without loss of generality assume $\Gamma=\gamma((-\delta,\delta))$ where $\gamma\colon(-3\delta,3\delta)\to\R^3$ is smooth with $\vert\gamma'\vert\equiv 1$ and $\gamma(0)=q$. Fix $x_j\in(-\delta,\delta)$ such that $\gamma(x_j)=q_j$. We claim that
  \begin{equation}\label{eq:1-pf-lem:boundary-integral-jump}
    \sup_{s\in(-\delta,\delta)}\Big\vert \frac{\langle \mathbf N(\gamma(x_j+s)), \gamma(x_j+s)-(q_j+d_jv_j)\rangle}{\vert \gamma(x_j+s)-(q_j+d_jv_j) \vert^2} + d_j\frac{\langle\mathbf N(q_j),v_j\rangle}{s^2+d_j^2} \Big\vert \leq C 
  \end{equation} 
  for some $C>0$ independent of $j$. Indeed, by Taylor expansion, we have for a constant $C>0$ which in the following may change from term to term that
  \begin{equation}\label{eq:2-pf-lem:boundary-integral-jump}
    \vert\gamma(x_j+s)-q_j-s\gamma'(x_j)\vert\leq C s^2,\quad\text{and}\quad \vert\mathbf N(\gamma(x_j+s))-\mathbf N(q_j)\vert\leq C\vert s\vert.
  \end{equation}
  Using $\langle\mathbf N(q_j),\gamma'(x_j)\rangle=0$, this yields $\vert \langle\mathbf N(\gamma(x_j+s)),\gamma(x_j+s)-q_j\rangle\vert \leq Cs^2$. Furthermore, since $\langle v_j,\gamma'(x_j)\rangle=0$, \eqref{eq:2-pf-lem:boundary-integral-jump} yields that
  \begin{equation}
    \big\vert \vert \gamma(x_j+s)-(q_j+d_jv_j) \vert^2-(s^2+d_j^2) \big\vert \leq C(\vert s\vert^3+d_js^2).
  \end{equation}
  Combining these estimates and using $d_j\to 0$ yields \eqref{eq:1-pf-lem:boundary-integral-jump} for $\delta>0$ sufficiently small. Now the lemma follows from dominated convergence, using \eqref{eq:1-pf-lem:boundary-integral-jump},  
  \begin{equation}
    \frac{\langle \mathbf N(\gamma(x_j+s)), \gamma(x_j+s)-(q_j+d_jv_j)\rangle}{\vert \gamma(x_j+s)-(q_j+d_jv_j) \vert^2} + d_j\frac{\langle\mathbf N(q_j),v_j\rangle}{s^2+d_j^2} \to \frac{\langle \mathbf N(\gamma(s)),\gamma(s)-q\rangle}{\vert\gamma(s)-q\vert^2} \quad\text{for $s\neq 0$, }
  \end{equation}
  and finally that
  \begin{equation}
    - \int_{-\delta}^{\delta} d_j \frac{\langle \mathbf N(q_j),v_j\rangle}{s^2+d_j^2}\dd s = -2\langle \mathbf N(q_j),v_j\rangle  \arctan(\delta/d_j) \to -\pi \langle \mathbf N(q),v\rangle 
  \end{equation}
  for $j\to\infty$.
\end{proof}

\begin{lemma}\label{lem:area-bound-Dirichlet}
  With $\Gamma$ and $\mathbf N$ as above, $C(\Gamma,\mathbf N)\vcentcolon=\sup_{z\in\R^3}\int_\Gamma\frac{\vert\langle\mathbf N(p),p-z\rangle\vert}{|p-z|^2}\dd\mathcal H^1(p)<\infty$. Let $f\colon\Sigma\to\R^3$ with $f=\gamma$ and $\partial_{\eta_f}f=\mathbf N\circ\gamma$ on $\partial\Sigma$. Then, for all $r>0$ and $z\in\R^3$, we have
  \begin{equation}\label{eq:uniform-area-ratio-Dirichlet}
    \frac{\int_{f^{-1}(B_r(z))}\dd\vol_{g_f}}{r^2}\leq\W(f)+C(\Gamma,\mathbf N).
  \end{equation}
  This estimate, with the same right-hand-side, is invariant with respect to dilation and translation of $f$. Furthermore, for any $z\notin f(\Sigma)$, we have
  \begin{equation}\label{eq:inversion-identity-Dirichlet}
    \int_\Sigma\Big|\frac14\mathbf H_f\mathbf n_f+\frac{(f-z)^\perp}{|f-z|^2}\Big|^2\dd\vol_{g_f}=\frac{\W(f)}4+\frac12\int_\Gamma\frac{\langle\mathbf N(p),p-z\rangle}{|p-z|^2}\dd\mathcal H^1(p).
  \end{equation}
\end{lemma}
\begin{proof}
  Arguing similarly as in the proof of \Cref{lem:boundary-integral-jump}, we get uniformly for $z$ in a fixed tubular neighborhood of $\Gamma$,
  \begin{equation}
    \int_\Gamma\frac{\vert\langle\mathbf N(p),p-z\rangle\vert}{|p-z|^2}\dd\mathcal H^1(p) \leq C+C\int_{-1}^1\frac{s^2+\operatorname{dist}(z,\Gamma)}{s^2+\operatorname{dist}(z,\Gamma)^2}\dd s\leq C.
  \end{equation}
  So $C(\Gamma,\mathbf N)<\infty$ follows since the integral converges to $0$ for $\vert z\vert\to\infty$. Since $\Sigma$ is compact, $\theta^2(f,\infty)=0$. Thus \eqref{eq:monotonicity-with-boundary} and \eqref{eq:limits-monotonicity-with-boundary}, yield
  \begin{align}
    &\frac{\int_{f^{-1}(B_r(z))}\dd\vol_{g_f}}{r^2}+\frac1{16}\int_{f^{-1}(B_r(z))}\mathbf H_f^2\dd\vol_{g_f}+\frac1{2r^2}\int_{f^{-1}(B_r(z))}\mathbf H_f\langle\mathbf n_f,f-z\rangle\dd\vol_{g_f}\\
    &\qquad\leq\frac{\W(f)}4+\frac12\int_\Gamma\frac{\vert\langle\mathbf N(p),p-z\rangle\vert}{\max\{r^2,|p-z|^2\}}\dd\mathcal H^1(p).
  \end{align}
  Young's inequality applied to the mixed term yields \eqref{eq:uniform-area-ratio-Dirichlet}, as in \cite[Equations~(A.5), (A.6)]{kuwertschaetzle2004}, and the claimed invariance follows by scaling and translation. Finally, \eqref{eq:inversion-identity-Dirichlet} follows from \eqref{eq:global-monotonicity-with-boundary}, since $\theta^2(f,z)=\theta^2(f,\infty)=0$ for $z\notin f(\Sigma)$, see also \cite[Lemma~2.1]{novagapozzetta2020} and, for the corresponding point-density estimate, \cite[Proposition~2.3]{schlierf2025acv}.
\end{proof}


\printbibliography

\end{document}